\documentclass[11pt,a4paper,reqno]{amsart}%
\usepackage{amsthm,amsmath,amsfonts,amssymb,xcolor,amsxtra,appendix,bookmark,dsfont,bm,mathrsfs}
\usepackage[latin1]{inputenc}
\usepackage{color} 
\usepackage{amssymb}\usepackage{graphicx}
\usepackage{tikz}
\usepackage{enumitem}
\usepackage{hyperref}
\theoremstyle{plain}
\newtheorem{theorem}{Theorem}%[section]
\newtheorem{lemma}[theorem]{Lemma}
\newtheorem{corollary}[theorem]{Corollary}

\newtheorem{definition}{Definition}

\theoremstyle{remark}
\newtheorem{remark}[theorem]{Remark}

\newcommand{\R}{\mathbb R}
\newcommand{\N}{\mathbb N}
\newcommand{\C}{\mathbb C}
\newcommand{\Z}{\mathbb Z}
\newcommand{\T}{\mathbb T}

\newcommand{\E}{\mathbb E}

\newcommand{\sA}{\mathscr A}
\newcommand{\sN}{\mathscr N}

\newcommand{\sS}{\mathscr S}

\newcommand{\sW}{\mathscr W}

\newcommand{\cL}{\mathcal L}
\newcommand{\cF}{\mathcal F}

\newcommand{\cI}{\mathcal I}

\DeclareMathOperator{\Imagine}{Im}
\DeclareMathOperator{\Real}{Re}

\DeclareMathOperator{\sign}{sign}
\DeclareMathOperator{\dist}{dist}
\DeclareMathOperator{\arcsinh}{arcsinh}

\newcommand{\bbZ}{{\mathbb Z}}

\newcommand{\bbR}{{\mathbb R}}

\newcommand{\bbT}{{\mathbb T}}

\newcommand{\om}{\omega}
\newcommand{\omm}{\omega_{\rm min}}
\newcommand{\omM}{\omega_{\rm max}}

\allowdisplaybreaks

\title[Controlling the rates of a chain of harmonic oscillators]{Controlling the rates of a chain of harmonic oscillators with a point Langevin thermostat}
\author[A. Hannani]{Amirali Hannani
}
\address{Institute for Theoretical Physics, KU Leuven, 3001 Leuven, Belgium}
\email{amirali.hannani@kuleuven.be
}

\author[M.-N. Phung]{Minh-Nhat Phung}
\address{Department of Mathematics, Texas A\&M University, College Station, TX 77843, USA}
\email{pmnt1114@tamu.com}

\author[M.-B. Tran]{Minh-Binh Tran}
\address{Department of Mathematics, Texas A\&M University, College Station, TX 77843, USA}
\email{minhbinh@tamu.edu} 
\thanks{M.-B. T is  funded in part by  a   Humboldt Fellowship,   NSF CAREER  DMS-2303146, and NSF Grants DMS-2204795, DMS-2305523, DMS-2306379.}
\author[E. Tr\'elat]{Emmanuel Tr\'elat}
\address{Sorbonne Universit\'e, CNRS, Universit\' e de Paris, Inria, Laboratoire Jacques-Louis Lions (LJLL),F-75005 Paris, France}
\email{emmanuel.trelat@sorbonne-universite.fr} 
\begin{document}

	\date{\today}

	\maketitle

	\begin{center}
		``Dedicated to the 60th Birthday of Professor Dang Duc Trong''
	\end{center}

	\begin{abstract} 
		We consider the control problem of controlling the rates of an infinite chain of coupled harmonic oscillators with a Langevin thermostat at the origin. We study the effect of  two types of open-loop boundary controls, impulsive control and linear memory-feedback control, in the high frequency limit. We investigate their action on the reflection-transmission coefficients for the
		wave energy for the scattering of the thermostat.
		Our study shows that the impulsive boundary controls have no impact on the rates and are thus not appropriate to act on the system, despite their physical meaning and relevance. In contrast, the second kind of control that we propose, which is less standard and uses the past of the state solution of the system, is adequate and relevant. We prove that any triple of rates satisfying appropriate assumptions is asymptotically reachable thanks to the linear memory-feedback controls that we design explicitly.
	\end{abstract}

	\tableofcontents

\section{Introduction}

Heat reservoirs at temperature $T$ are commonly modeled by the Langevin stochastic dynamics. When the bulk evolution is governed by a discrete wave equation, a small parameter $\varepsilon>0$ is introduced to dictate the ratio between microscopic and macroscopic space-time units. In general, the noise is chosen  so that by the stochastic mechanism, there is only  a finite amount of momentum exchanged  in a   finite interval of time. As thus, each particle undergoes only a finite number of collisions in any finite interval of time. 
It is common to let  $\varepsilon$ tend to $0$, which is often referred to as taking the kinetic limit for the system (for instance, see \cite{LS2007}, \cite{Mielke2006}, \cite{RPK96}). When a  chain has no microscopic boundary, the energy density evolution is often described by a linear kinetic equation.

A  useful tool to localize in space the energy per frequency mode is the 
Wigner distribution, introduced in \cite{PhysRev.40.749}.
In the absence of the thermostat, {by adding a small conservative noise exchanging velocities}, the authors of \cite{basile2010energy} prove that, in the kinetic limit $\varepsilon\to 0$, the Wigner distribution
converges to the solution of the kinetic transport equation
\begin{equation}
	\label{eq:bos1}
	\partial_t W(t,x,k) + v_g(k) \partial_x W(t,x,k) =
	2\gamma_0 \int_{\mathbb{T}} \mathcal R(k,k') \left(W(t,x,k') - W(t,x,k)\right)dk',
\end{equation}
for all $(t,x,k)\in [0,+\infty)\times \mathbb{T}\times \mathbb R$. The  explicit   scattering kernel $\mathcal R(k,k') \ge 0$ is given by 
\begin{equation*}
	%\label{Rk}
\mathcal	R(k):={\int}_{\mathbb{T}}\mathcal R(k,k')dk' \sim |k|^2\quad\mbox{ for }|k| \ll 1.
\end{equation*} Here and in the sequel $\mathbb{T}$ is the unit torus, identified with the interval $[-1/2,1/2]$ with periodic 
endpoints. The parameter $\gamma_0>0$ is the scattering rate for the microscopic chain. The group velocity is defined by $v_{g}(k) := \omega(k)/2\pi$, 
$\omega(k)$ is the dispersion relation of the chain.

When a heat bath at temperature $T$ is applied to one particle, which is labeled $\alpha_0$,
with a coupling strength $\gamma_1>0$, the action of the heat bath is not affected by the scale of the small parameter 
	$\varepsilon$.  As a consequence, 
	 when a thermostat is included in the system, its presence  can be regarded as a singular
	 perturbation of the dynamics of the  system. Mathematically speaking, when $\varepsilon\to 0$, in  \cite{komorowski2020kinetic,KORH2020}, it has been proved that the thermostat enforces that phonons of wave number $k$ are generated with rate $r_{a}(k)T$, incoming $k$-phonons can be  transmitted with probability $r_{t}(k)$ and reflected with probability $r_{r}(k)$, which means that one needs to introduce the boundary  conditions
at $x=0$ on \eqref{eq:bos1}:
\begin{equation}\label{eq:bos2}
	\begin{split}
		W(t,0^+,k) &=r_{r}(k) W(t,0^+, -k) + r_{t}(k)W(t, 0^-,k)+r_{a}(k)T, \quad\hbox{ for $0< k\le 1/2$},\\
		W(t,0^-, k) &=r_{r}(k)W(t,0^-,-k) + r_{t}(k) W(t,0^+, k) +r_{a}(k)T,\quad \hbox{ for $-1/2< k< 0$.}
	\end{split}
\end{equation}
Those quantities are properly normalized according to $ r_{a}(k) + r_{t}(k) + r_{r}(k) = 1$,
so that $W(t,x,k) = T$ is a  thermal equilibrium.

 In the recent years, there have been significant progresses on the control theory for kinetic models \cite{carrillo2022controlling,piccoli2015control}. The goal of our work is to initiate the study of the Wigner distributions for   stochastic discrete wave equations under the point of view of control theory. 
 To be more precise, in the setting of the stochastic discrete wave equations considered in \cite{komorowski2020kinetic,KORH2020}, the three important parameters $r_{a}(k),r_{t}(k),r_{r}(k)$ are respectively the probabilities for absorption, transmission and reflection of a phonon of mode $k$.   Within the control theory viewpoint, we investigate the following question:

 \vspace{0.75em}
{\it {\bf Question A - Controllability at the kinetic limit: }%\footnote{We want to distinguish between ``limit controllability'' and ``exact controllability'': Given a time interval $[0,T_0]$, initial and final states, the exact controllability problem is   to determine a control such that the solution under the influence of the control satisfies both the initial  and final state. In order words, in ``exact controllability'', there is no limit taken.}
if we add a control function, having a relevant physical meaning, to the wave system/equation, can we control  the above three important rates $r_a(k),r_t(k),r_r(k)$: absorption, transmission, and reflection of a phonon of mode $k$ in \eqref{eq:bos2}?}
\vspace{0.75em}

 %At the first sight, this is a little bit different from the standard question in control theory:  Can we add some control term in the PDE model, having a relevant physical meaning, so that   we can drive the solution from the initial condition at time $0$ to a targeted function at time $\tau$. However, both questions serve a common goal: to control the solutions of the underline physical equations.  
 To address Question A, we consider two possible types of open-loop boundary controls. 
  
 The first boundary control is impulsive (see Section \ref{Impulsivecontrol} or further reading at \cite{1130000794184306176,komorowski2023heat}) and consists of adding the term $F(t)\delta_{0,n}$ to the system. At the kinetic limit such a control has an impact only at the boundary, causing a shoot-up explained by the friction $\nu$ at the boundary. Like the thermostat, the control force $F(t)$ can be seen as a wave at all frequencies but only the frequency of the oscillators are kept; all other frequencies are damped by oscillations at the macroscopic limit.
 {Physically, this type of control corresponds to adding a force $F(t)$ at the boundary. Similar physical phenomena (subjecting a chain of oscillators to a point force) have been considered recently in \cite{komorowski2023heat}. The difference with our work lies in the nature of the force: $L^1$ in time in our case vs. periodic in time in theirs. Another difference is the scaling regime that we consider for the control of the kinetic limit, in contrast to their setup that corresponds to the diffusive scaling. } Our control is impulsive in the sense that it creates a new term in the kinetic limit equation (see \eqref{kineticlimit} further), namely, a Dirac delta measure along the characteristic $\delta(x-v_g(k)t)$. This means that the control is too strong and the whole dynamics of the chain  follows the control and not its own dynamics. %A well-known example for an impulsive behavior  can be observed in the following simple 1-dimensional control system
%$$ x'(t) = u(t), \ \   x(0)=0,$$
%in which $x(t)$ represents the system under consideration and $u(t)$ is the added control function. The impulsive behavior happens when we take the control to be a Dirac delta measure $u=\delta_{t=0}$ and the solution $x(t)$ of the system becomes a delta measure. In our case, 
Even when we take a smooth function $F$,  an extra Dirac delta measure along the characteristic $\delta(x-v_g(k)t)$ pops up in the final equation, as proved in \eqref{kineticlimit}. 
It follows that the probabilities for absorption, transmission and reflection of a phonon of mode $k$ are not controllable using this control (see Remark  \ref{controlcoro0}).
  
The second type of open-loop boundary control is a linear memory-feedback control (see Section \ref{LMFcontrol}) on the stochastic process of the wave consisting of adding the convolution term $(F+\nu\delta(t))\star \alpha_{0}\delta_{0,n}$. In control theory, the time-convolution integral for a continuous-time system calculates the output of a   system to a given input using the  response of the system. 
  Therefore it is natural to utilize time-convolution  as a tool to control the system at any time $t$ given that the anterior states of the system are known.  In this respect, the control uses feedback from previous times to control outputs at the current time  and is often called memory-feedback. Since  the control in this case follows closely the dynamics of the system,   the ``impulsive'' behavior in which a delta function pops up in the final result can be avoided.
  Noticing that the three quantities $r_{a},r_{t},r_{r}$ depend on the friction $\nu$, a time-convolution control can be imposed on this parameter $\nu$, which sticks to the stochastic process $\alpha_{0}$. Thanks to this feedback-type boundary control, eventually the effect of the control on the kinetic limit is much better: the three rates $(r_{a},r_{t},r_{r})$ can be changed at the kinetic limit. We will prove that for a given class of function rates $(r_{a},r_{t},r_{r})$, enjoying appropriate conditions, it is possible to find a control $F$ asymptotically steering the rates to desired target functions (see Corollary \ref{controlcoro}).

Another type of control, which seems to be technically harder is the so-called ``internal control'', for which, the control acts on several points of the chain, will be studied in future work. 
{Let us mention that instead of a Langevin thermostat, one may also consider a Poisson scattering mechanism at the boundary. Such mechanisms are studied in \cite{komorowski2022asymptotic}. Controlling the rate of such problems seems doable via our methods.}

{\bf Acknowledgments:} The authors would like to express their gratitude to Tomasz Komorowski, Stefano Olla,    Herbert Spohn and Enrique Zuazua for several useful remarks, guidance, and instructions on the topic.  
\section{Setting and Notations}

We consider the evolution of an infinite particle system governed by the Hamiltonian 
\begin{equation*}
	{\mathcal H}({\alpha},{\beta}):=\frac12\sum_{n\in\bbZ}{\alpha}_n^2+\frac{1}{2}\sum_{n,n'\in\bbZ}\sigma_{n-n'}{\beta}_n{\beta}_{n'}=\frac12\sum_{n\in\Z}\|\alpha\|^2+\langle \beta,\sigma\star \beta\rangle,
\end{equation*}
where $n\in\bbZ$, $({\alpha}_n,{\beta}_n)$ is the position and momentum of the particle  $n$, and $({\alpha},{\beta})=\{({\alpha}_n,{\beta}_n),\,n\in\bbZ\}$. The assumptions on $\sigma$ will be specified in Section \ref{Assumdisper}.

The convolution of two functions on $\bbZ$ is 
\begin{equation}
	(f\star g)_n:=\sum_{n'\in\bbZ}f_{n-n'}g_{n'}.\label{convolutiony}
\end{equation}
We denote by $\delta_{0,n}$ and $\delta(x)$ the Dirac delta functions on $\Z$ and $\R$ respectively.
The Hamiltonian dynamics with stochastic source without control reads
\begin{equation*}
	\begin{split}
	\dot{\beta}_n(t)&= {\alpha}_n(t),
	 \\
	 d{\alpha}_n(t) &=  -(\sigma\star {\beta})_n(t)dt
	+\big(-\nu {\alpha}_0(t)dt+\sqrt{2\nu T}dR(t)\big) \delta_{0,n},\quad n\in\bbZ,
	\end{split}
\end{equation*}
where $\{R(t),\,t\ge0\}$   is a Wiener process on a  probability space with proper filtration $(\Omega, \mathfrak{F}_t, \mathbb{P})$ and an initial probability measure $\mu_{\varepsilon}$ on $\ell^2(\mathbb{Z})$. In this setting, we couple the particle whose label is $ 0$ to a Langevin thermostat at temperature $T$ and we assume that the friction is  $\nu>0$. 

We denote by $\E$ the expectation with respect to $\mathbb{P}$ and the expectation with respect to $\mu_{\varepsilon}$ is denoted by $\langle\cdot\rangle_{\mu_{\varepsilon}}$.
The expectation for the two processes, the Wiener process and the initial measure $\mu_{\varepsilon}$, is denoted by $\E_{\varepsilon}$.

The Fourier transform of  $f_n\in \ell^2(\bbZ)$ and the inverse Fourier transform of $\hat f\in L^2(\bbT)$ are
\begin{equation*}
	%\label{fourier}
	\hat f(k):=\sum_{x\in\bbZ} f_n \exp\{-2\pi ink\}, ~~
	%\end{equation}
	%The  inverse transform is given by
	%\begin{equation}
	%\label{inv-fourier}
	f_n:=\int_{\bbT} \hat f(k) \exp\{2\pi ink\} dk, \quad n\in \bbZ,~~k\in\bbT.
\end{equation*}
The wave function is given by
\begin{equation}
	\label{011307}
	\psi_n(t) := (\tilde{\om} \star{\beta})_n(t) + i{\alpha_n}(t)
\end{equation}
%where  $({\frak q}(t),{\frak p}(t))$ satisfies \eqref{eq:bas2}. 
in which $\{\tilde{\omega}_n,n\in\Z\}$ is the inverse Fourier transform of the dispersion relation defined by
$ 
	\om(k):=\sqrt{\hat \sigma (k)}.
$
The Fourier transform of the wave function is
\begin{equation*}
	%\label{011307a}
	\hat\psi(t,k) := \om(k) \hat {\beta}(t,k) + i\hat{\alpha}(t,k)\qquad\forall k\in\T.\end{equation*}
 
We also have
 
$$
\hat{\alpha}\left(t,k\right)=\frac{1}{2i}[\hat\psi(t,k)-\hat\psi^*(t,-k)],
~~{\alpha}_0(t)=\int_{\bbT} {\rm Im}\,\hat\psi(t,k) dk.
$$

For a function $O(x,k)$ in the Schwartz space $\sS(\R\times\T)$, we denote by $\hat O:\bbR\times\bbT\to\mathbb C$, $\tilde O:\bbR\times\bbZ\to\mathbb C$ the Fourier transforms of $O$ in the $x$ and $k$ variables, respectively,
\begin{align*}
	\hat{O}(\xi,k):=\int_{\R}e^{-2\pi i \xi x}O(x,k)dx,%\label{FOfirst}
	\\
	\tilde{O}(x,n):=\int_{\T}e^{-2\pi i n k}O(x,k)dk.%\label{FOsecond}
\end{align*}

The Laplace transform for the time variable is
\begin{align*}
	\cL(f)(Z):=\int_{0}^{\infty}e^{-Z t}f(t)dt
	%\label{Laplace}
\end{align*}
Also, for now, if we don't say anything further, we assume the domain for convergence of the Laplace transform is $\C_{+}=\{z\in\C:\Real(z)>0\}$.
We keep the notations $\Real(z),\Imagine(z),z^{*}$ for the real part, imaginary part and complex conjugate of $z$.

The Laplace transform of the Wiener process is denoted by
\begin{align*}
	\tilde{R}:=\cL(R).
	%\label{wienerprocess}
\end{align*}
One can observe that $\tilde{R}$ is a Gaussian process, determined by its covariance which is given by: 
\begin{equation*}
    \mathbb{E}(\tilde{R}(Z_1)\tilde{R}(Z_2))= \frac{1}{Z_1+Z_2}, \qquad
    \text{Re} Z_1, \text{Re} Z_2 >0.
\end{equation*}

Recall the convolution notation in \eqref{convolutiony}. We also use the notation $\star$ for convolutions with respect to the time variable $t$, as
\begin{equation*}
	f\star g(t) := \int_{0}^{t}f(t-s)g(s)ds.
	%\label{convolutiont}
\end{equation*}
We use this for convenience. To avoid confusion, we have a convention that we only use this notation when dealing with functions $F,\Theta,\Psi,C_\omega$. We define those functions in Section \ref{Impulsivecontrol} and Section \ref{LMFcontrol}.

The $\varepsilon$ time scaling of a function $f$ is defined by
\begin{align}
	f^{(\varepsilon)}(t):=f(t/\varepsilon).
	\label{timescaling}
\end{align}

For the estimates, we use the following symbols: for $f,g:D\to\R$, we write
\begin{align*}
	&f\lesssim g\text{ if there exists }C>0:f(x)\le Cg(x),x\in D.%\label{lesssimsymbol}
	\\
	&\text{We write }f\approx g\text{ if }f\lesssim g\text{ and }g\lesssim f.
	%\label{estimatesymbol}
\end{align*}
Our estimation usually involves $\varepsilon\to0^{+}$, hence when we use those symbols we mean that $D$ is a small positive neighborhood of $0$ for the variable $\varepsilon$ and $C$ does not depend on $\varepsilon$.

We now state a few basic assumptions.
\subsection{Assumptions on the initial wave}

We define the space $\sN$ as the completion of $\sS(\R\times\T)$ for the norm
\begin{align*}
	\|O\|_{\sN}:=\int_{\R}\sup_{k\in\T}|\hat{O}(\xi,k)|d\xi
	%\label{normO}
\end{align*}
and we consider it dual space $\sN'$, where the dual is defined in \eqref{Wigner}.

The Wigner distribution is defined by
\begin{align}
	\langle O,W_\varepsilon(t)\rangle_{\sN}:=\frac{\varepsilon}{2}\sum_{n,n'\in\Z}\E_{\varepsilon}\left[ \hat{\psi}^{(\varepsilon)}_{n}(t)\hat{\psi}^{(\varepsilon)*}_{n'}(t) \right]\tilde{O}^{*}(\varepsilon\frac{n+n'}{2},n-n'),
	\label{Wigner}
\end{align}
for every test function $O\in\sS(\R\times\T)$.

The Wigner distribution is also defined by its Fourier transform, as
\begin{align}
	\hat{W}_{\varepsilon}(t,\xi,k)&:= \frac{\varepsilon}{2}\E_{\varepsilon}\left[ \hat{\psi}^{(\varepsilon)*}(t,k-\varepsilon\xi/2)\hat{\psi}^{(\varepsilon)}(t,k+\varepsilon\xi/2) \right],\label{FWigner}\\
	\langle O,W_\varepsilon(t)\rangle_{\sN}&= \int_{\T\times\R}\hat{W}_{\varepsilon}(t,\xi,k)\hat{O}^{*}(\xi,k)d\xi dk.\nonumber%\label{FWignerwithO}
\end{align}
For the Wigner distribution's Laplace transform, we use the notation
\begin{align*}
	\sW_{\varepsilon}(Z,\xi,k):=\cL(\hat{W}_{\varepsilon}(\cdot,\xi,k))(Z).
	%\label{LaplaceW}
\end{align*}
For most parts, we will find the limit of the Laplace transform $\sW_{\varepsilon}$ instead of directly find the limit of $\hat{W}_{\varepsilon}$.
When we mention about the initial condition, we omit the time variable and we write
\begin{align*}
	\hat{\psi}(k)=\hat{\psi}(0,k)\quad\text{and}\quad\hat{W}_{\varepsilon}(\xi,k)=\hat{W}_{\varepsilon}(0,\xi,k).
	%\label{omittime}
\end{align*}

The assumptions for the initial state are the following
\begin{enumerate}[label=(I\arabic*)]
	\item \label{init0}$\langle\hat{\psi}(k),\hat{\psi}(h)\rangle_{\mu_{\varepsilon}}= 0,\quad \forall k,h\in\T$.
	\item\label{energyassumption} $\sup_{\varepsilon\in(0,1]}\sum_{n\in\Z}\varepsilon\langle|\psi_{n}|^2\rangle_{\mu_{\varepsilon}}= \sup_{\varepsilon\in(0,1]}\varepsilon\langle\|\hat{\psi}\|^{2}_{L^2(\T)}\rangle_{\mu_{\varepsilon}}<\infty$.
	\item\label{WbracketR} There exist constants $C_1, \kappa>0$ such that $|\hat{W}_{\varepsilon}(\xi,k)|\le \frac{C_1}{\langle \xi\rangle_{\R}^{3+\kappa}}$,
for $(\xi,k)\in \T_{2/\varepsilon}\times\T,\varepsilon\in(0,1]$ ($\T_{2/\varepsilon}$ is the periodic torus $[-\varepsilon^{-1},\varepsilon^{-1}]$). Here, the Japanese bracket is defined by
\begin{align*}
	\langle x\rangle_{\R}:= \sqrt{1+x^2}.
\end{align*}
Note that, for any $\kappa>0$, $\int_{\R}\langle x\rangle_{\R}^{-1-\kappa} dx<\infty$.
\item $\hat{W}_{\varepsilon}(\xi,k)$ converges weakly in $\sN'$ to $\hat{W}_0\in L^1(\R\times\T)$.\label{W_0}\\
The energy of the wave grows at most linearly in time (Lemma \ref{ligrowimpul} or Lemma \ref{linearfeed}), hence
\begin{align*}
	\sup_{t\in[0,\tau]}\|W_{\varepsilon}(t)\|_{\sN'}<\infty,
	%\label{boundW}
\end{align*}
for every $\tau>0$. Therefore, the family of all $W_{\varepsilon}$ is sequentially weak-star compact in $(L^{1}([0,\tau];\sN))^{\star}$. We obtain this assumption by taking a subsequence of $\varepsilon$ if it is neccessary.
\end{enumerate}

\subsection{Assumptions on the dispersion relation}
\label{Assumdisper}
The coupling between two points $n,n'\in\Z$ is denoted by
\begin{align*}
	\sigma_{n-n'}\in\C.
	%\label{coupling}
\end{align*}
The dispersion relation is defined using this coupling
\begin{align*}
	\omega(k):=\sqrt{\hat{\sigma}(k)}.
	%\label{dispersion}
\end{align*}
The function $k\mapsto\omega(k)$ is smooth, even, and is positive on $\T$ (following \ref{expdecay} - \ref{positivesigma}).
Furthermore, we have $\omega'(k)>0$ for $k\in(0,1/2)$ (following \ref{monoome}).
The minimum and maximum values of the dispersion relation are
\begin{equation*}
	\omm= \omega(0)\quad\text{and}\quad\omM= \omega(-1/2)=\omega(1/2).
\end{equation*}
Since $\omega$ is continuous and strictly increases on $[0,1/2]$, given any $\overline{\omega}\in[\omm,\omM]$ there exists a unique $k\in[0,1/2]$ such that $\omega(k)=\overline{\omega}$. We define two inverse functions $\varphi_{+}:[\omm,\omM]\to[0,1/2]$ and $\varphi_{-}:[\omm,\omM]\to[-1/2,0]$, the value of $\varphi(\overline{\omega})$ is set as $k$. We means that
\begin{equation*}
	\varphi_{+}(\omega(k)) := |k|\quad\text{and}\quad \varphi_{-}(\omega(k)) := -|k|.
\end{equation*}

The following assumptions are imposed.
\begin{enumerate}[label=($\sigma$\arabic*)]
	\item\label{expdecay} There exist $C_2,C_3>0$ such that
		\begin{align*}
			|\sigma_{n}|\le C_2e^{-C_3|n|},\qquad\forall n\in\Z.	
		\end{align*}
	\item $\hat{\sigma}$ is even.%\label{evensigma}
	\item $\hat{\sigma}(k)>0$ for $k\in\T$.\label{positivesigma}
\end{enumerate}

\begin{enumerate}[label=($\omega$\arabic*)]
	\item\label{monoome} $\omega$ is decreasing on the negative branch, or equivalently, increasing on the positive branch.
\item \label{inverome} There are smooth positive functions $\varphi_{1}$ and $\varphi_{2}$ such that
		\begin{align}
			\varphi'_{\pm}(\overline{\omega})&= \pm(\overline{\omega}-\omm)^{-1/2}\varphi_{1}(\overline{\omega}),\quad \overline{\omega}-\omm\ll 1,\label{minapprox}\\
			\varphi'_{\pm}(\overline{\omega})&= \pm(\omM-\overline{\omega})^{-1/2}\varphi_{2}(\overline{\omega}),\quad \omM-\overline{\omega}\ll 1.\label{maxapprox}
		\end{align}
\end{enumerate}

\begin{remark}
	We consider one example for the dispersion relation (see \cite{basile2010energy}):
\begin{align*}
	\omega(k)=\sqrt{\omega_0^{2}+\gamma(1-\cos(2\pi k))},
	%\label{exampledispersion}
\end{align*}
where $\omega_0,\gamma>0$. For $\overline{\omega}\in[\omega_0,\sqrt{\omega_0^2+2\gamma}]$, we have
\[\varphi_\pm(\overline{\omega})=\pm\frac{1}{\pi}\arcsin\sqrt{\frac{\overline{\omega}^2-\omega_0^2}{2\gamma}}\]
and
\[\varphi_\pm'(\overline{\omega})=\pm\frac{\overline{\omega}}{\pi\sqrt{2\gamma+\omega_0^2-\overline{\omega}^2}\sqrt{\overline{\omega}^2-\omega_0^2}}.\]
When $\overline{\omega}-\omega_0\ll1$, we define
\[\varphi_1(\overline{\omega})=\frac{\overline{\omega}}{\pi\sqrt{2\gamma+\omega_0^2-\overline{\omega}^2}\sqrt{\overline{\omega}+\omega_0}},\]
and when $\sqrt{\omega_0^2+2\gamma}-\overline{\omega}\ll 1$, we define
\[\varphi_2(\overline{\omega})=\frac{\overline{\omega}}{\pi\sqrt{\sqrt{\omega_0^2+2\gamma}+\overline{\omega}}\sqrt{\overline{\omega}^2-\omega_0^2}}.\]
\end{remark}

Finally, we introduce some notations related to the dispersion relation
\begin{equation}
	\begin{split}
			D_{\varepsilon}(k,\xi)&:= \varepsilon^{-1}\left[ \omega(k+\varepsilon\xi/2)-\omega(k-\varepsilon\xi/2) \right]\qquad\text{for }k\in\T,\xi\in\T_{2/\varepsilon},\\
	D_{\varepsilon}^{+}(k,\xi)&:= \varepsilon^{-1}\left[ \omega(k+\varepsilon\xi)-\omega(k) \right]\qquad\text{for }k\in\T,\xi\in\T_{1/\varepsilon},%\label{plusDepsilon}
	\\
	D_{\varepsilon}^{-}(k,\xi)&:= \varepsilon^{-1}\left[ \omega(k)-\omega(k-\varepsilon\xi) \right] \qquad\text{for }k\in\T,\xi\in\T_{1/\varepsilon}.%\label{minusDepsilon}
	\end{split}
\label{Depsilon}
\end{equation}
When $\varepsilon\to0$, $D_{\varepsilon}(k,\xi), D_{\varepsilon}^{+}(k,\xi), D_{\varepsilon}^{-}(k,\xi)$ all converge to $\omega'(k)\xi$.

\subsection{Asumptions on target rates}
Given target rates $(r_a,r_t,r_r)$, we want to reach them using our controls.
\begin{definition}
	For a triple of functions $(r_a,r_t,r_r)$ on $\T$, it is called ``asymptotically reachable'' if there exist a control function $F$ and a sequence $\{F_N\}_{N=1}^{\infty}$ such that $F_N$ is equal to $F$ in finite time and $r_a^{F_N},r_t^{F_N},r_r^{F_N}$ converge almost everywhere to $r_a,r_t,r_r$ on $\T$ as $N\to\infty$.
	\label{asymreach}
\end{definition}

We state some assumptions on the target rates to make a sufficient condition for ``asymptotically reachable''
		
	\begin{enumerate}[label=(H\arabic*)]
		\item\label{evenT} The three functions $(r_a,r_t,r_r)$ are even on $\T$;
		\item\label{sum=1} $r_a(k)+r_t(k)+r_r(k)=1$ for all $k\in\T$;
		\item \label{posrates}$r_t(k),r_r(k)>0$ and there exists $c_1>0$ such that $r_a(k)\ge c_1$ for $k\in(0,1/2)$;
		\item \label{conrates}The three functions $(r_a,r_t,r_r)$ are continuous on $(0,1/2)$;
			\item\label{sumsquareroot} $\sqrt{r_t(k)}+\sqrt{r_r(k)}\ge 1$ for all $k\in\T$;
			\item \label{intcossin}$\int_{0}^{\infty}\Real(\bar{F}(u))\cos(ut)du=\int_{0}^{\infty}\Imagine(\bar{F}(u))\sin(ut)du$ for all $t>0$.
		\end{enumerate}
		In \ref{intcossin}, $\bar{F}$ is defined by
		\begin{itemize}
			\item $RE(k):=|v_g(k)|(r_t(k)-r_r(k)-1)$ for $k\in(0,1/2)$;
			\item $IM(k):=\sqrt{4|v_g(k)|^2r_r(k)-RE(k)^2}$ for $k\in(0,1/2)$;
			\item $FT(k):=RE(k)+iIM(k)$ for $k\in(0,1/2)$;
			\item $TH(k):= 1+FT(k)\lim_{Z\to0}\cL(C_\omega)(Z-i\omega(k))$ for $k\in(0,1/2)$;
			\item $\bar{F}(u):=FT(\varphi_+(u))/TH(\varphi_+(u))$ for $u\in(\omm,\omM)$.
		\end{itemize}
We have the definition of the group velocity
\begin{align*}
	v_{g}(k):=\frac{\omega'(k)}{2\pi}.
	%\label{groupvelo}
\end{align*}

\section{Main results}
In this section, we will first show in Theorem \ref{Main} the transport equation with impulsive control to see that the rates cannot be controlled in this way. In contrast, we show in Theorem \ref{mainnew} the equation with the rates having control using linear memory-feedback control. Finally, we state the result on a sufficient condition to get asymptotically reachable rates.

We define the rates without control
\begin{align*}
	r_{a}(k)&:= \frac{\nu|\theta(k)|^2}{|v_{g}(k)|},%\label{absorbrate}
	\\
	r_{t}(k)&:= 1-\Real(\theta(k))\frac{\nu}{|v_{g}(k)|}+\frac{\nu r_a(k)}{4|v_{g}(k)|},%\label{transmission}
	\\
	r_{r}(k)&:= \frac{\nu r_a(k)}{4|v_{g}(k)|},%\label{reflection}
\end{align*}
where $\theta$ is defined further in \eqref{thetadef}.
Recall that the terms $r_a,r_t$ and $r_r$ represent the absorption rate, transmission rate and reflection rate. 

\subsection{Impulsive control}
The impulsive control term is associated with a real-valued function $F(t)$. We consider the system of equations
\begin{align}
	\dot{\beta}_n(t)&= \alpha_n(t),\label{dbeta}\\
	d\alpha_n(t)&= -(\sigma\star \beta(t))_ndt+(-\nu \alpha_0dt+F(t)dt+\sqrt{2\nu T}dR(t))\delta_{0,n}.\label{dalpha}
\end{align}
In this system, we assume some additional assumptions.
\begin{enumerate}[label=(C\arabic*)]
	\item \label{Passumption}$\varepsilon^{d}\left|P(\varepsilon,k-\varepsilon\xi/2)\right|$ is bounded for each fixed $k$, $\xi$, where $P(\varepsilon,k)=\langle \hat{\psi}(k)\rangle_{\mu_{\varepsilon}}$ and $d<1/6$.
	\item \label{Fassumption} There exists $\cF \in C(\T)$ such that $\varepsilon^{1/2}\cL(F)(\varepsilon Z+ i\omega(k-\varepsilon\xi/2))\to\cF(k)$ for any fixed $\xi$ as $\varepsilon \to 0$.
	\item \label{integralF^2bound}There exists $C_4>0$ such that $\sup_{\varepsilon\in(0,1]}\int_{0}^{t/\varepsilon}\varepsilon F^{2}(s)ds\le C_4t,\forall t>0$.
\end{enumerate}

	\begin{theorem}
	Consider the wave system governed by \eqref{dbeta} and \eqref{dalpha} along with all the assumptions \ref{init0} - \ref{W_0}, \ref{expdecay} - \ref{positivesigma}, \ref{monoome} - \ref{inverome} and \ref{Passumption} - \ref{integralF^2bound}. For any $\tau>0$ and $O\in L^{1}([0,\tau];\sN)$ we have
	\begin{align}
		\lim_{\varepsilon\to0}\int_{0}^{\tau}\langle O(t),W_{\varepsilon}(t)\rangle_{\sN} dt=\int_{0}^{\tau}\int_{\R\times\T}O^{*}(t,x,k)W(t,x,k)dxdkdt,
		\label{kineticlimit}
	\end{align}
where
\begin{align}
	W(t,x,k)&= W_0(x-v_g(k)t,k)1_{[0,v_g(k)t]^{c}}(x)+r_{a}(k)T1_{[0,v_g(k)t]}(x)\nonumber\\
	&+r_{t}(k)W_0(x-v_g(k)t,k)1_{[0,v_g(k)t]}(x)+r_{r}(k)W_0(-x+v_g(k)t,-k)1_{[0,v_g(k)t]}(x)\nonumber\\
	&+\frac{|v_g(k)|r_a(k)|\cF(k)|^2}{\nu }\delta(x-v_g(k)t).
	\label{Wresult}
\end{align}
	\label{Main}
\end{theorem}

\begin{remark}\label{controlcoro0}
	In a related context, the authors of \cite{komorowski2023heat} enforced impulsive forces on two ends of a finite harmonic chain to study a stationary state relating to the chain. Therefore, the question of impulsive control is quite natural and has already considered in \cite{komorowski2023heat} for a different model.

\end{remark}

\begin{remark}
	In view of Question A, posed in our introduction, we observe that the three rates $r_{a},r_{t},r_{r}$ are not modified under the influence of the impulsive control $F(t)\delta_{0,n}$. The control appears in the term $\frac{|v_g(k)|r_a(k)|\cF(k)|^2}{\nu}\delta(x-\bar{\omega}'(k)t)$ in \eqref{Wresult}($\cF$ is given by \ref{Fassumption}: this is where the control $F(t)\delta_{0,n}$ appears in \eqref{Wresult}). With respect to the main result of \cite{KORH2020}, the latter term is the new additional one.

	Therefore, despite its physical meaning and relevance, such an impulsive control is not appropriate to act on the thermostat system. In Section \ref{LMFcontrol}, we will propose another type of control called memory-feedback because it uses the information provided by the past of the state solution. This alternative kind of control enables us to act on the triple rates.
\end{remark}

\subsection{Linear memory-feedback control}
The linear memory-feedback control is the convolution control term $(F+\nu \delta(t))\star \alpha_{0}(t)\delta_{0,n}$. Hence, we will consider the system of equations:
\begin{align}
	\dot{\beta}_n(t)&= \alpha_n(t),\label{dbetanew}\\
	d\alpha_n(t)&= -(\sigma\star \beta(t))_ndt+(F\star \alpha_0(t)dt+\sqrt{2\nu T}dR(t))\delta_{0,n}.\label{dalphanew}
\end{align}
In the system, we use the assumption
\begin{enumerate}[label=(L\arabic*)]
	\item \label{inequalitycondition1}$F\in C^{\infty}((0,+\infty))\cap L^m((0,+\infty)),m\in[1,\infty)$ and $\Real\left( \hat{F}\left( \frac{\omega(k)}{2\pi} \right) \right)<0$ for $k\in\T$.
\end{enumerate}
Because $F \in L^1$ we have the equality \eqref{LaplaceFourier} and the Fourier transform $\hat{F}$ is a bounded and continuous function.

We define the three new rates and use them for our main theorem 
\begin{align}
	r_a^{F}(k)&:= -\Real(\hat{F}(\omega(k)/2\pi))\frac{|\theta_{F}(k)|^{2}}{|v_g(k)|},\label{absorbratenew1}\\
	r_t^F(k)&:= 1+\frac{\Real(\hat{F}^{*}(\omega(k)/2\pi)\theta_{F}(k))}{|v_g(k)|}+\frac{|\hat{F}(\omega(k)/2\pi)|^{2}|\theta_{F}(k)|^{2}}{4|v_g(k)|^{2}},\label{transmissionnew1}\\
	r_r^F(k)&:= \frac{|\hat{F}(\omega(k)/2\pi)|^{2}|\theta_{F}(k)|^{2}}{4|v_g(k)|^{2}}.\label{reflectionnew1}
	\end{align}
	We define the function $\theta_F$ later in \eqref{tildeThetaFlimit}.

\begin{theorem}
	Consider the system of equations \eqref{dbetanew} and \eqref{dalphanew} with all the assumptions \ref{init0} - \ref{W_0}, \ref{expdecay} - \ref{positivesigma}, \ref{monoome} - \ref{inverome} and \ref{inequalitycondition1}. Then, for any $\tau>0$ and $G\in L^{1}([0,\tau];\sA)$ we have
	\begin{align*}
		\int_{0}^{\tau}\left\langle G(t),W_{\varepsilon}(t)\right\rangle dt=\int_{0}^{\tau}\int_{\R\times\T}G^{*}(t,x,k)W(t,x,k)dxdkdt,
		%\label{kineticlimit1}
	\end{align*}
	where
	\begin{align*}
		W(t,x,k)&= W_{0}(x-v_g(k)t,k)1_{[0,v_g(k)t]^{c}}(x)-\frac{\nu Tr_a^{F}(k)}{\Real(\hat{F}(\omega(k)/2\pi))}1_{[0,v_g(k)t]}(x)\\
		&+r_t^{F}(k)W_{0}(x-v_g(k)t,k)1_{[0,v_g(k)t]}(x)+r_r^F(k)W_{0}(-x+v_g(k)t,-k)1_{[0,v_g(k)t]}(x).
		%\label{Wresult1}
	\end{align*}
	\label{mainnew}
\end{theorem}
On the well-posedness of \eqref{dbetanew}-\eqref{dalphanew}, we recall \eqref{Psi} in which the solution $\hat{\psi}$ is uniquely determined by \eqref{psiexpandnew}.

\begin{remark}
	In view of addressing Question A, posed in the introduction, we observe that the three rates $(r_a^F,r_t^F,r_r^F)$ can now be controlled under the influence of the linear memory-feedback control. An ``inverse problem'' question is: given  a triple of functions $(r_a,r_t,r_r)$ as target rates, can we determine if $(r_a,r_t,r_r)$ can be asymptotically reached following \eqref{absorbratenew1}-\eqref{transmissionnew1}-\eqref{reflectionnew1} in sense of Definition \ref{asymreach}?
\end{remark}

\begin{corollary}	\label{controlcoro}
	If a triple $(r_a,r_t,r_r)$ satisfies \ref{evenT} - \ref{intcossin}, then it is asymptotically reachable.
\end{corollary}

\begin{remark}
	We emphasize on the convergence type used in asymptotically reachable is only in the a.e. sense, which is weak. We desire to improve on this convergence type. Another highly important issue is the rate at which the control is achieved. For both issues, we hope that a better construction of $F_N$ would allow one to obtain a better convergence type with a good convergence rate. To be more precise, in the design of  $\bar{F}$ of \ref{intcossin}, the function $F_N$ is constructed by multiplying the inverse Fourier of $\bar{F}$ with a smooth cut-off function, whose values are $1$ on $[0,N],$ are in $[0,1]$ on $(N,N+1)$, and are $0$ on $[N+1,+\infty)$. This simple construction of $F_N$ gives the desired achievable control in a.e. sense. However, this construction is not sufficient to obtain a better convergence type nor the rate at which the control is achieved. A better construction of $F_N$ may give a better result. However, such computations would be rather very long and we will devote this issue for an upcoming work. 
\end{remark}
\begin{remark}
	As discussed above, the work \cite{komorowski2023heat} establishes the first impulsive control, in which the force is periodic. Time periodic control is therefore a very interesting direction to explore. For instance, it would be important to see whether or not our linear-feedback control (or even different types of control) can be done in a periodic in time manner. We will investigate this question in a follow-up paper. 
\end{remark}
\section{Impulsive control - Proof of Theorem \ref{Main}}
\label{Impulsivecontrol}

\subsection{Preliminary computations}
In this section, we work with the impulsive control introduced in \eqref{dbeta} and \eqref{dalpha}.
We use the notation $\hat{F}$ for the Fourier transform of $F$ in the time variable. We recall a basic relation between Laplace and Fourier transforms. For $F\in L^{1}([0,+\infty))$ we have
\begin{align}
	\lim_{Z\to0}\cL(F)(Z+i\omega)=\lim_{Z\to0}\int_{0}^{+\infty}e^{-Z t-i\omega t}F(t)dt=\int_{0}^{+\infty}e^{-i\omega t}F(t)dt=\hat{F}(\omega/2\pi).
	\label{LaplaceFourier}
\end{align}
We clarify the reasons for the assumptions in this control. In \ref{Passumption}, we have $d<1/6$ so that $P$ is small enough, when the control acts on the wave, the new term involving $P$ will disappear when taking the limit. Assumption \ref{Fassumption} is a genelization of \eqref{LaplaceFourier}. The bound in \ref{integralF^2bound} ensures the growth of the energy is at most linear.
We prove the following lemma using an It\^{o} analysis shown in the Section \ref{proligrowimpul}.
\begin{lemma}
	Consider the wave system governed by \eqref{dbeta} and \eqref{dalpha} with all the assumptions \ref{init0} - \ref{W_0}, \ref{expdecay} - \ref{positivesigma}, \ref{monoome} - \ref{inverome} and \ref{Passumption} - \ref{integralF^2bound}. There is a constant $C_5$ such that
	\begin{equation*}
		\varepsilon\E_{\varepsilon}\|\hat{\psi}^{(\varepsilon)}(t)\|^{2}_{L^2(\T)}\le \varepsilon\E_{\varepsilon}\|\hat{\psi}^{(\varepsilon)}(0)\|^2_{L^2(\T)}+C_5t.
	\end{equation*}
	\label{ligrowimpul}
\end{lemma}

In the proof of Theorem \ref{Main}, to analyze the kinetic limit, the Wigner distribution is splitted into $13$ terms that we categorize into $4$ types.
\begin{enumerate}
	\item Terms not involving the control: one of the thermal term \eqref{thermal4}, the ballistic term \eqref{ballistic}, the first and the second scattering terms \eqref{I}, \eqref{II}.
	\item Terms with a single Wiener process: \eqref{thermal2} and \eqref{thermal3}. Since the expectation of the Wiener process is zero, thus these terms are zero.
	\item Terms involving one occurence of $F$ in the product: those terms are \eqref{term3deriLaplace}, the third, the fourth, the fifth and the sixth scattering terms, \eqref{III}, \eqref{IV}, \eqref{V}, \eqref{VI}.
	\item Terms involving two occurences of $F$ in the product: \eqref{thermal1} and \eqref{VII}.
\end{enumerate}
We will see that only the first and the fourth types have non-zero kinetic limits. All other terms vanish at the limit.

We write $d\psi$ using the system of equations  (see \cite{hannani2022wave,staffilani2021wave} for a related nonlinear problem)
\begin{align}
	d\hat{\psi}(t,k)&= d(\omega(k)\hat{\beta}(t,k))+id\hat{\alpha}(t,k)\nonumber\\
	&= \omega(k)\hat{\alpha}(t,k)dt+i(-\omega^2(k)\hat{\beta}dt-\nu \alpha_0(t)dt+F(t)dt+\sqrt{2\nu T}dR(t))\nonumber\\
	&= -i\omega(k)\hat{\psi}(t,k)dt-i\nu \alpha_0(t)dt+iF(t)dt+\sqrt{2\nu T}idR(t).
	\label{dpsi}
\end{align}
Solving the linear ODE \eqref{dpsi} gives
\begin{align}
	\hat{\psi}(t,k)=e^{-i\omega(k)t}\hat{\psi}(k)-i\nu \int_{0}^{t}e^{-i\omega(k)(t-s)}\alpha_0(s)ds&+i\int_{0}^{t}e^{-i\omega(k)(t-s)}F(s)ds\nonumber\\
	&\quad+i\sqrt{2\nu T}\int_{0}^{t}e^{-i\omega(k)(t-s)}dR(s).\label{solvepsi}
\end{align}
Noting that by definition of inverse Fourier transform,
\begin{align}
	\alpha_0(t)=  \int_{\T}\Imagine(\hat{\psi}(t,k))dk,
	\label{alpha_0}
\end{align}
and that the real part of $e^{-i\omega(k)(t-s)}$ is $\cos(\omega(t-s))$, we infer from \eqref{solvepsi} and \eqref{alpha_0}
\begin{align}
	\alpha_0(t)&= \Psi(t)-\nu \int_{\T}\int_{0}^{t}\cos\left( \omega(k)(t-s) \right)\alpha_0(s)dsdk+\int_{\T}\int_{0}^{t}\cos\left( \omega(k)(t-s) \right)F(s)dsdk\nonumber\\
	&\hspace{15em}+\sqrt{2\nu T}\int_{\T}\int_{0}^{t}\cos\left( \omega(k)(t-s) \right)dR(s)dk\nonumber%\label{p_0cos}
	\\
	&=\Psi(t)-\nu C_\omega\star \alpha_0(t)+C_\omega\star F(t)+\sqrt{2\nu T}C_\omega\star dR(t),
	\label{p_0Comega}
\end{align}
where
\begin{align}
	\Psi(t):=\int_{\T}\Imagine\left( \hat{\psi}(k)e^{-i\omega(k)t} \right)dk,
	\label{Psi}
\end{align}
and
\begin{align}
	C_\omega(t):=\int_{\T}\cos(\omega(k)t)dk.
	\label{Comegadef}
\end{align}
We can compute the Laplace transform
\begin{equation}
	\cL(C_\omega)(Z)=\int_{\T}\frac{Z}{Z^2+\omega^2(k)}dk.
	\label{LComega}
\end{equation}
Recall that $\Real(Z)>0$ so $\Real(\cL(C_\omega)(Z))>0$. For $\cL(C_\omega)$, we also get the following lemma. We put the detailed proof for this lemma in Section \ref{prolimLC}.
\begin{lemma}
	For $k\in(0,1/2)$, we have
	\begin{equation*}
		\lim_{Z\to0}\cL(C_\omega)(Z-i\omega(k))= \frac{\pi}{|\omega'(k)|}+i\left( \int_{0}^{1/2}\frac{dh}{\omega(k)+\omega(h)}+\log\frac{\omega(k)-\omm}{\omM-\omega(k)}+\frac{C_{\omega,0}(\omega(k))}{|\omega'(k)|} \right),
	\end{equation*}
	where $C_{\omega,0}$ is a continuous, bounded and real-valued function.
	\label{limLC}
\end{lemma}

Using the Laplace transform on \eqref{p_0Comega} and the notation
\begin{align*}
	\tilde{\Theta}(Z):=\frac{1}{1+\nu \cL(C_\omega)(Z)}=\sum_{j=0}^{\infty}(-\nu \cL(C_\omega)(Z))^j,
	%\label{tildeTheta}
\end{align*}
we get
\begin{align}
	\tilde{\alpha}_{0}(Z)=\cL(\alpha_0)(Z)&= \tilde{\Theta}(Z)\cL(\Psi)(Z) +\tilde{\Theta}(Z)\cL(C_\omega)(Z)\cL(F)(Z)+\sqrt{2\nu T}\tilde{\Theta}(Z)\cL(C_\omega)(Z)\tilde{R}(Z)\nonumber\\
	&= \tilde{\Theta}(Z)\cL(\Psi)(Z)+\frac{1}{\nu}(1-\tilde{\Theta}(Z))\cL(F)(Z) +\sqrt{2\nu T}\tilde{\Theta}(Z)\cL(C_\omega)(Z)\tilde{R}(Z).
	\label{tildealpha_0}
\end{align}
We will write $\Theta$ as the inverse Laplace transform of $\tilde{\Theta}$ and it can be understood using the infinite sum of convolution terms
\begin{align*}
	\Theta(t)&:= \cL^{-1}(\tilde{\Theta})=\delta(t)+\cL^{-1}(\tilde{\Theta}-1)=\delta(t)+\sum_{n=1}^{\infty}(-\nu )^{n}(\star^{n}C_\omega),
	%\label{Thetadef}
\end{align*}
where $(\star^n C_\omega)$ is the $n$-time convolution of $C_\omega$ with itself.

Applying the Laplace transform to \eqref{solvepsi}, and using \eqref{tildealpha_0} we obtain
\begin{align*}
	\cL(\hat{\psi}(\cdot,k))(Z)&= \frac{\hat{\psi}(k)-i\nu\tilde{\alpha}_0(Z)+i\cL(F)(Z)+i\sqrt{2\nu T}\tilde{R}(Z)}{Z+i\omega(k)}\\
	&= \frac{\hat{\psi}(k)-i\nu \cL(\Psi)(Z)+i\tilde{\Theta}(Z)\cL(F)(Z)+i\sqrt{2\nu T}\tilde{\Theta}(Z)\tilde{R}(Z)}{Z+i\omega(k)}.
	%\label{Lhatpsi}
\end{align*}
Therefore, by inverse Laplace transform,
\begin{align}
	\hat{\psi}(t,k)&= e^{-i\omega(k)t}\hat{\psi}(k)-i\nu \int_{0}^{t}\Theta_\omega(t-s,k)\Psi(s)ds+i\int_{0}^{t}\Theta_\omega(t-s,k)F(s)ds\nonumber\\
	&\hspace{15em}+i\sqrt{2\nu T}\int_{0}^{t}\Theta_\omega(t-s,k)dR(s),
	\label{psiexpand}
\end{align}
where
\begin{align*}
	\Theta_\omega(t,k):= \int_{0}^{t}e^{-i\omega(k)(t-s)}d\Theta(s).
	%\label{Theta_omega}
\end{align*}
Denoting by
\begin{align}
	\Theta^\omega(t,k):= \int_{0}^{t}e^{i\omega(k)s}d\Theta(s)= e^{i\omega(k)t}\Theta_\omega(t,k),
	\label{Theta^omega}
\end{align}
we have the important identity
\begin{align*}
	\cL((\Theta^\omega)^{(\varepsilon)}(\cdot,k))(Z)=\frac{\tilde{\Theta}(\varepsilon Z-i\omega(k))}{Z}.
	%\label{Theta^omegaTheta}
\end{align*}
Because $\Real(1+\nu\cL(C_\omega)(Z))>1$ we get $\tilde{\Theta}(\varepsilon Z-i\omega(k))$ is bounded by $1$. Hence, using  Fatou's Theorem, $\tilde{\Theta}(\varepsilon Z-i\omega(k))$ converges (a.e. and in any $L^{p}$) when $Z\to 0$. We also write
\begin{align}
	\theta(k):=\lim_{\varepsilon\to0^{+}}\tilde{\Theta}(\varepsilon Z-i\omega(k)).
	\label{thetadef}
\end{align}
$\theta$ is used to define the rates, it has an important property that allows us to compute the sum of the rates. Lemma \ref{thetaprop} is an immediate consequence of Lemma \ref{limLC}.
\begin{lemma}
	We have
	\begin{equation*}
		\Real(\theta(k))=\left( 1+\frac{\nu\pi}{|\omega'(k)|} \right)|\theta(k)|^2.
	\end{equation*}
	As a consequence, $r_a(k)+r_t(k)+r_r(k)=1$ for a.e. $k\in\T$.
	\label{thetaprop}
\end{lemma}

We consider the truncation
\begin{align}
	\hat{\psi}^{1}(0,k)&= \hat{\psi}(k)\chi_{\varrho}(k),\label{multiplybump}\\
	d\hat{\psi}^{1}(t,k)&= \left\{ -i\omega(k)\hat{\psi}^{1}(t,k)-\frac{\nu}{2i}\int_{\T}[\hat{\psi}^{1}(t,k')-\hat{\psi}^{1*}(t,k')]dk'+F(t) \right\}dt+i\sqrt{2\nu T}dR(t).\nonumber\\%\label{newpsi}
	\hat{\psi}^{2}(0,k)&= \hat{\psi}(k)(1-\chi_{\varrho}(k)),\nonumber%\label{trimpart}
	\\
	d\hat{\psi}^{2}(t,k)&= \left\{ -i\omega(k)\hat{\psi}^{1}(t,k)-\frac{\nu }{2i}\int_{\T}[\hat{\psi}^{2}(t,k')-\hat{\psi}^{2*}(t,k')]dk' \right\}dt.\nonumber%\label{eliminatedpsi}
\end{align}
Here, $\chi_{\varrho}(k)$ is a smooth function with values in $[0,1]$ such that it is equal to $0$ on $[\T](\varrho)$ and to $1$ on $\T\setminus [\T](2\varrho)$,
where
\begin{equation}
\begin{split}
	[\T](\varrho):=&  \{k:\dist(k,[\T](0))<\varrho\},\\
	[\T](0):=&  \{k:\omega'(k)=0\}=\{0,\pm1/2\}.
\end{split}
\label{Tvarrho}
\end{equation}
In the computations of scattering terms (Section \ref{ScatIm}), we will consider $\psi^1$ in \eqref{multiplybump} as $\psi$. For a sufficiently small $\varrho$, we will show that $\psi^1$ in \eqref{multiplybump} and $\psi$ in \eqref{dpsi} give the same kinetic limit (Section \ref{proofmain}). Thus, we only consider the situation where the initial state satisfies
\begin{equation}
	\text{There exists }\varrho>0\text{ such that the support of }\hat{\psi}(0,\cdot)\text{ is in }\T\setminus[\T](\varrho).
\label{Trunass}
\end{equation}

We now split the wave into several terms to compute the kinetic limit.

\subsection{Thermal terms}
The thermal part of the wave is considered independently from the initial state of the wave system. Therefore, for the thermal part, we put $\hat{\psi}(k)\equiv 0$. With that assumption, \eqref{psiexpand} becomes
\begin{align*}
	\hat{\psi}(t,k)=i\int_{0}^{t}\Theta_\omega(t-s,k)F(s)ds+i\sqrt{2\nu T}\int_{0}^{t}\Theta_\omega(t-s,k)dR(s).
	%\label{thermalwave}
\end{align*}
The Wigner distribution's definition \eqref{FWigner} gives
\begin{align}
	\hat{W}_\varepsilon(t,\xi,k)=\frac{\varepsilon}{2}\E_\varepsilon\left[\int_{0}^{t/\varepsilon}\Theta_\omega(t/\varepsilon-s,k+\varepsilon\xi/2)F(s)ds\int_{0}^{t/\varepsilon}\Theta_\omega^{*}(t/\varepsilon-s,k-\varepsilon\xi/2)F(s)ds\right.\label{thermal1}\\
	+\sqrt{2\nu T}\int_{0}^{t/\varepsilon}\Theta_\omega(t/\varepsilon-s,k+\varepsilon\xi/2)F(s)ds\int_{0}^{t/\varepsilon}\Theta_\omega^{*}(t/\varepsilon-s,k-\varepsilon\xi/2)dR(s)\label{thermal2}\\
	+\sqrt{2\nu T}\int_{0}^{t/\varepsilon}\Theta_\omega^{*}(t/\varepsilon-s,k-\varepsilon\xi/2)F(s)ds\int_{0}^{t/\varepsilon}\Theta_\omega(t/\varepsilon-s,k+\varepsilon\xi/2)dR(s)\label{thermal3}\\
	+\left.2\nu T\int_{0}^{t/\varepsilon}\Theta_\omega(t/\varepsilon-s,k+\varepsilon\xi/2)dR(s)\int_{0}^{t/\varepsilon}\Theta_\omega^{*}(t/\varepsilon-s,k-\varepsilon\xi/2)dR(s)\right].
	\label{thermal4}
\end{align}

\subsubsection{First thermal term \eqref{thermal1}}\label{thermal1re}
Since the term is independent from any random process, its Laplace transform is
\begin{align}
	\frac{\varepsilon}{2}\cL\left( (\Theta_\omega\star F)^{(\varepsilon)}(\cdot,k+\varepsilon\xi/2)(\Theta_\omega^{*}\star F)^{(\varepsilon)}(\cdot,k-\varepsilon\xi/2)\right).
	\label{Laplacethermal1}
\end{align}
We rewrite the convolution $\Theta_\omega\star F$ using \eqref{Theta^omega} as
\begin{align}
	\Theta_\omega\star F(t,k)=\int_{0}^{t}\Theta_\omega(t-s,k)F(s)ds&= \int_{0}^{t}\Theta^\omega(t-s,k)e^{i\omega(k)(s-t)}F(s)ds\nonumber\\
	&= \int_{0}^{t}\Theta^\omega(t-s,k)e^{-i\omega(k)t}F^{*}_{k}(s)ds\nonumber\\
	&= e^{-i\omega(k)t}\Theta^\omega\star F^{*}_{k}(t),\label{rewriteconvolution}
\end{align}
where
\begin{align}	
	F_{k}(t):=F(t)e^{-i\omega(k)t}.
	\label{Fk}
\end{align}
We also rewrite the conjugation,
\begin{align}
	\Theta_\omega^{*}\star F(t,k)=e^{i\omega(k)t}(\Theta^\omega)^{*}\star F_{k}(t).
	\label{rewriteconvolutionstar}
\end{align}
Using \eqref{Depsilon}, \eqref{rewriteconvolution} and \eqref{rewriteconvolutionstar}, the product of the two convolutions is
\begin{align}
	&(\Theta_\omega\star F(t,k+\varepsilon\xi/2))(\Theta_\omega^{*}\star F(t,k-\varepsilon\xi/2))\nonumber\\
	&\quad= e^{-i\varepsilon D_{\varepsilon}(k,\xi)t}(\Theta^\omega\star F^{*}_{k+\varepsilon\xi/2})(t,k+\varepsilon\xi/2)( (\Theta^\omega)^{*}\star F_{k-\varepsilon\xi/2})(t,k-\varepsilon\xi/2).
	\label{productconvolution}
\end{align}
Now, we use \eqref{productconvolution}, \eqref{Laplacethermal1}, and the formula of the Laplace transform of a product to get
\begin{align}
	\frac{\varepsilon}{4i\pi}&\lim_{L\to\infty}\int_{c-iL}^{c+iL}\cL( (\Theta^\omega\star F^{*}_{k+\varepsilon\xi/2})^{(\varepsilon)})(\zeta)\cL( ( (\Theta^\omega)^{*}\star F_{k-\varepsilon\xi/2})^{(\varepsilon)})(Z+iD_{\varepsilon}(k,\xi)-\zeta)d\zeta\nonumber\\
	&= \frac{\varepsilon}{4i\pi}\lim_{L\to\infty}\int_{c-iL}^{c+iL}\cL( (\Theta^\omega)^{(\varepsilon)})(\zeta,k+\varepsilon\xi/2)\cL(F^{*}_{k+\varepsilon\xi/2})(\varepsilon\zeta)\nonumber\\
	&\qquad\cL( (\Theta^\omega)^{(\varepsilon)*})(Z+iD_{\varepsilon}(k,\xi)-\zeta,k-\varepsilon\xi/2)\cL(F_{k-\varepsilon\xi/2})(\varepsilon(Z+iD_{\varepsilon}(k,\xi)-\zeta))d\zeta%\label{expandLaplacethermal1}
	\nonumber\\
	&= \frac{\varepsilon}{4i\pi}\lim_{L\to\infty}\int_{c-iL}^{c+iL}\frac{\tilde{\Theta}(\varepsilon\zeta-i\omega(k+\varepsilon\xi/2))}{\zeta}\cL(F)(\varepsilon\zeta-i\omega(k+\varepsilon\xi/2))\nonumber\\
	&\qquad\frac{\tilde{\Theta}(\varepsilon(Z+iD_{\varepsilon}(k,\xi)-\zeta)+i\omega(k-\varepsilon\xi/2))}{Z+iD_{\varepsilon}(k,\xi)-\zeta}\cL(F)(\varepsilon(Z+iD_{\varepsilon}(k,\xi)-\zeta)+i\omega(k-\varepsilon\xi/2))d\zeta,\nonumber
\end{align}
where we follow the notation of \eqref{timescaling}.
From \eqref{thetadef} and \ref{Fassumption}, taking the limit $\varepsilon\to0$, we obtain
\begin{align}
	\frac{1}{4i\pi}\lim_{L\to\infty}\int_{c-iL}^{c+iL}\frac{|\theta(k)|^2|\cF(k)|^2}{\zeta(Z+i\omega'(k)\xi-\zeta)}d\zeta
	.\label{thermal1limitepsilon}
\end{align}
We have
\begin{align*}
	\frac{1}{\zeta(Z+i\omega'(k)\xi-\zeta)}=\frac{1}{Z+i\omega'(k)\xi}\left( \frac{1}{\zeta}+\frac{1}{Z+i\omega'(k)\xi-\zeta} \right).
	%\label{basicfraction}
\end{align*}
We also have the identity
\begin{align}
	\lim_{L\to\infty}\int_{c-iL}^{c+iL}\frac{1}{\zeta}d\zeta&=i\lim_{L\to\infty}\left(  \int_{-L}^{L}\frac{c}{c^2+x^2}dx-i\int_{-L}^{L}\frac{x}{c^2+x^2}dx\right)\nonumber\\
	&= i\pi.
	\label{integralfraction}
\end{align}
Hence, \eqref{thermal1limitepsilon} equals
\begin{align*}
	\frac{|\theta(k)|^2|\cF(k)|^2}{2(Z+i\omega'(k)\xi)}=\frac{|v_g(k)|r_a(k)|\cF(k)|^{2}}{2\nu (Z+i\omega'(k)\xi)}.
	%\label{thermal1limit}
\end{align*}

\subsubsection{Second and third thermal terms \eqref{thermal2},\eqref{thermal3}}\label{thermal23re}
Each term has one $dR(s)$, and because they follow a Gaussian distribution, those terms vanish.

\subsubsection{Last thermal term \eqref{thermal4}}\label{thermal4re}
Because of the correlation of Wiener process, \eqref{thermal4} is equal to
\begin{equation*}
	\nu T\varepsilon\int_{0}^{t/\varepsilon}\Theta_\omega(t/\varepsilon-s,k+\varepsilon\xi/2)\Theta_\omega^*(t/\varepsilon-s,k-\varepsilon\xi/2)ds=\nu T\int_{0}^{t}\Theta_\omega^{(\varepsilon)}(s,k+\varepsilon\xi/2)(\Theta_{\omega}^*)^{(\varepsilon)}(s,k-\varepsilon\xi/2)ds.
\end{equation*}
Similarly to the computation of \eqref{thermal1}, by using the identity \eqref{Theta^omega} and the formula of the Laplace transform of a product, the last thermal term gives
\begin{align*}
	\sW_{\varepsilon}(Z,\xi,k)= \frac{\nu T}{Z}\int_{0}^{\infty}e^{-(Z+iD_{\varepsilon}(k,\xi))s}\Theta^\omega\left( \varepsilon^{-1}s,k+\varepsilon\xi/2 \right)(\Theta^\omega)^{*}\left( \varepsilon^{-1}s,k-\varepsilon\xi/2 \right)ds\\
	= \frac{\nu T}{Z}\frac{1}{2\pi i}\lim_{L\to\infty}\int_{c-iL}^{c+iL}\frac{\tilde{\Theta}(\varepsilon\zeta-i\omega(k+\varepsilon\xi/2))\tilde{\Theta}^{*}(\varepsilon(Z+iD_{\varepsilon}(k,\xi)-\zeta)-i\omega(k-\varepsilon\xi/2))}{\zeta(Z+iD_{\varepsilon}(k,\xi)-\zeta)}d\zeta.
\end{align*}
The limit of $\sW_{\varepsilon}$ is
\begin{align*}
	\frac{\nu T|\theta(k)|^2}{Z(Z+i\omega'(k)\xi)}=\frac{T|v_g(k)|r_a(k)}{Z(Z+i\omega'(k)\xi)}.
	%\label{thermal4limit}
\end{align*}

\subsection{Scattering terms}\label{ScatIm}
From \eqref{FWigner}, we have 
\begin{align}
	\partial_t\hat{W}_{\varepsilon}(t,\xi,k)&= \frac{\varepsilon}{2}\E_{\varepsilon}\left[ \left( \partial_t\hat{\psi}^{(\varepsilon)}(t,k+\varepsilon\xi/2) \right)\hat{\psi}^{(\varepsilon)*}(t,k-\varepsilon\xi/2)+\left( \partial_t\hat{\psi}^{(\varepsilon)*}(t,k-\varepsilon\xi/2) \right)\hat{\psi}^{(\varepsilon)}(t,k+\varepsilon\xi/2) \right].
	\label{derivativeW}
\end{align}
Then, we compute each term of \eqref{derivativeW} by replacing $T=0$ in \eqref{dpsi}:
\begin{align}
	\left( \partial_t\hat{\psi}^{(\varepsilon)}(t,k+\varepsilon\xi/2) \right)&\hat{\psi}^{(\varepsilon)*}(t,k-\varepsilon\xi/2)\label{deri1}\\
	&= \frac{1}{\varepsilon}\left[-i\omega(k+\varepsilon\xi/2)\hat{\psi}^{(\varepsilon)}(t,k+\varepsilon\xi/2)-i\nu \alpha_0^{(\varepsilon)}(t)+iF^{(\varepsilon)}(t)  \right]\hat{\psi}^{(\varepsilon)*}(t,k-\varepsilon\xi/2),\nonumber\\
	\left( \partial_t\hat{\psi}^{(\varepsilon)*}(t,k-\varepsilon\xi/2) \right)&\hat{\psi}^{(\varepsilon)}(t,k+\varepsilon\xi/2)\label{deri2}\\
	&= \frac{1}{\varepsilon}\left[i\omega(k-\varepsilon\xi/2)\hat{\psi}^{(\varepsilon)*}(t,k-\varepsilon\xi/2)+i\nu \alpha_0^{(\varepsilon)}(t)-iF^{(\varepsilon)}(t)  \right]\hat{\psi}^{(\varepsilon)}(t,k+\varepsilon\xi/2).\nonumber
\end{align}
We infer from \eqref{deri1} and \eqref{deri2} that
\begin{align}
	&\partial_{t}\hat{W}_{\varepsilon}(t,\xi,k)=-iD_{\varepsilon}(k,\xi)\hat{W}_{\varepsilon}(t,\xi,k)\label{term1deri}\\
	&-\frac{\nu}{2}\E_{\varepsilon}\left[ i\alpha_0^{(\varepsilon)}(t)\hat{\psi}^{(\varepsilon)*}(t,k-\varepsilon\xi/2)-i\alpha_0^{(\varepsilon)}(t)\hat{\psi}^{(\varepsilon)}(t,k+\varepsilon\xi/2) \right]\label{term2deri}\\
	&\frac{i}{2}\E_{\varepsilon}\left[ F^{(\varepsilon)}(t)\hat{\psi}^{(\varepsilon)*}(t,k-\varepsilon\xi/2)-F^{(\varepsilon)}(t)\hat{\psi}^{(\varepsilon)}(t,k+\varepsilon\xi/2) \right].\label{term3deri}
\end{align}
Using the Laplace transform on \eqref{term1deri}, \eqref{term2deri} and \eqref{term3deri}, we obtain
\begin{align}
	&(Z+iD_{\varepsilon}(k,\xi))\sW_{\varepsilon}(Z,\xi,k)=\hat{W}_{\varepsilon}(\xi,k)\label{ballistic}\\
	&-\frac{\nu}{2}\left[ S_{\varepsilon}(Z,k-\varepsilon\xi/2)+S_{\varepsilon}^{*}(Z,k+\varepsilon\xi/2) \right]\label{Sepsilon}\\
	&-\cL\left( \frac{iF^{(\varepsilon)}(t)}{2}(\E_{\varepsilon}[ \hat{\psi}^{(\varepsilon)}(t,k+\varepsilon\xi/2)]-\E_{\varepsilon}[\hat{\psi}^{(\varepsilon)*}(t,k-\varepsilon\xi/2) ]) \right)(Z),\label{term3deriLaplace}
\end{align}
where
\begin{align}
	S_{\varepsilon}(Z,k):=i\cL\left( \E_{\varepsilon}\alpha_0^{(\varepsilon)}(t)\hat{\psi}^{(\varepsilon)*}(t,k) \right),\label{Sepsilondef}\\
	S_{\varepsilon}^{*}(Z,k):=-i\cL\left( \E_{\varepsilon}\alpha_0^{(\varepsilon)*}(t)\hat{\psi}^{(\varepsilon)}(t,k) \right).\nonumber%\label{Sepsilonstar}
\end{align}
To handle the term \eqref{Sepsilon} we split \eqref{Sepsilondef} into seven terms. First, rewrite \eqref{psiexpand} and \eqref{tildealpha_0} in convolution form with $T=0$
\begin{align}
	\hat{\psi}(t,k)&= e^{-i\omega(k)t}\hat{\psi}(k)-i\nu\Theta_\omega\star \Psi(t)+i\Theta_\omega\star F(t)\label{psiTzero},\\
	\alpha_0(t)&= \Theta\star \Psi(t)+\frac{1}{\nu}F(t)-\frac{1}{\nu}\Theta\star F(t).\label{alpha_0Tzero}
\end{align}
Then $S_{\varepsilon}$ is the sum of the Laplace transforms of $I_\epsilon^{(\varepsilon)},II_\epsilon^{(\varepsilon)},III_\epsilon^{(\varepsilon)},IV_\epsilon^{(\varepsilon)},V_\epsilon^{(\varepsilon)},VI_\epsilon^{(\varepsilon)},VII_\epsilon^{(\varepsilon)}$:
\begin{enumerate}
	\item $I_{\varepsilon}$ is obtained as the product of the first term of \eqref{alpha_0Tzero} and the first term of \eqref{psiTzero}:
		\begin{align}
			I_{\varepsilon}(t,k):=i\E_\varepsilon\left[ g\star \Psi(t)e^{i\omega(k)t}\hat{\psi}^{*}(k) \right]&= ie^{i\omega(k)t}\int_{0}^{t}\langle \Psi(t-s)\hat{\psi}^{*}(k)\rangle_{\mu_{\varepsilon}}d\Theta(s).\label{I}
		\end{align}
	\item $II_{\varepsilon}$ is obtained as the product of the first term of \eqref{alpha_0Tzero} and the second term of \eqref{psiTzero}:
		\begin{align}
			II_{\varepsilon}(t,k):=i\E_{\varepsilon}\left[ \Theta\star \Psi(t)\left( i\nu\Theta_\omega\star \Psi(t) \right) \right]&= -\nu\int_{0}^{t}\int_{0}^{t}\Theta_\omega^{*}(t-s,k)\langle \Psi(s)\Psi(t-s')\rangle_{\mu_{\varepsilon}}dsd\Theta(s').\label{II}
		\end{align}
	\item $III_{\varepsilon}$ is obtained as the product of the first term of \eqref{alpha_0Tzero} and the third term of \eqref{psiTzero}:
		\begin{align}
			III_{\varepsilon}(t,k):=i\E_{\varepsilon}\left[ \left( \Theta\star \Psi(t) \right)\left( -i\Theta_\omega^{*}\star F(t) \right) \right]&= \int_{0}^{t}\Theta_\omega^{*}(s,k)F(t-s)ds\left\langle\int_{0}^{t}\Psi(t-s')d\Theta(s')\right\rangle_{\mu_{\varepsilon}}.\label{III}
		\end{align}
	\item $IV_{\varepsilon}$ is obtained as the product of the second term of \eqref{alpha_0Tzero} and $\hat{\psi}^{*}$:
		\begin{align}
			IV_{\varepsilon}(t,k):=i\E_{\varepsilon}\left[ \left( \frac{1}{\nu}F(t) \right)\hat{\psi}^{*}(t,k) \right]=\frac{i}{\nu}F(t)\langle\hat{\psi}^{*}(t,k)\rangle_{\mu_{\varepsilon}}.\label{IV}
		\end{align}
	\item $V_{\varepsilon}$ is obtained as the product of the third term of \eqref{alpha_0Tzero} and the first term of \eqref{psiTzero}:
		\begin{align}
			V_{\varepsilon}(t,k):=i\E_{\varepsilon}\left[ \left( -\frac{1}{\nu}\Theta\star F(t) \right)\left( e^{i\omega(k)t}\hat{\psi}^{*}(k) \right) \right]&= \frac{-i}{\nu}\int_{0}^{t}F(t-s)d\Theta(s)e^{i\omega(k)t}\langle\hat{\psi}^{*}(k)\rangle_{\mu_{\varepsilon}}.\label{V}
		\end{align}
	\item $VI_{\varepsilon}$ is obtained as the product of the third term of \eqref{alpha_0Tzero} and the second term of \eqref{psiTzero}:
		\begin{align}
			VI_{\varepsilon}(t,k):=i\E_{\varepsilon}\left[ \left( -\frac{1}{\nu}\Theta\star F(t) \right)\left( i\nu\Theta_\omega^{*}\star \Psi(t) \right) \right]&= \int_{0}^{t}F(t-s)d\Theta(s)\left\langle\int_{0}^{t}\Theta_\omega^{*}(s',k)\Psi(t-s')ds'\right\rangle_{\mu_{\varepsilon}}.\label{VI}
		\end{align}
	\item $VII_{\varepsilon}$ is obtained as the product of the third term of \eqref{alpha_0Tzero} and the third term of \eqref{psiTzero}:
		\begin{align}
			VII_{\varepsilon}(t,k):=i\E_{\varepsilon}\left[ \left( -\frac{1}{\nu}\Theta\star F(t) \right)\left( -i\Theta_\omega^{*}\star F(t) \right) \right]&= -\frac{1}{\nu}\int_{0}^{t}F(t-s)d\Theta(s)\left\langle\int_{0}^{t}\Theta_\omega^{*}(s',k)F(t-s')ds'\right\rangle_{\mu_{\varepsilon}}.\label{VII}
		\end{align}
\end{enumerate}
We next show how to treat the scattering terms.

The ballistic term is the term generated from $\hat{W}_{\varepsilon}(\eta,k)$ in \eqref{ballistic}. The limit is
\begin{align}
	\int_{\R\times\T}\hat{O}(\xi,k)\frac{\hat{W}_{\varepsilon}(\xi,k)}{Z+iD_{\varepsilon}(k,\xi)}d\xi dk\to\int_{\R\times\T}\hat{O}(\xi,k)\frac{\hat{W}_{0}(\xi,k)}{Z+i\omega'(k)\xi}d\xi dk,\text{ as }\varepsilon\to0,
	\label{limitballistic}
\end{align}
where $O$ is the test function used in \eqref{Wigner}.
\subsubsection{The first scattering term \eqref{I}}\label{scat1re}
The first scattering term is obtained in \eqref{Sepsilon} using $I_{\varepsilon}$. The technique and result for the kinetic limit of this term is similar to that of Section \ref{Iscatterpoint}. In this case, the limit is
\begin{align*}
	-\nu\int_{\R\times\T}\Real(\theta(k))\frac{\hat{W}_{0}(\xi',k)}{Z+i\omega'(k)\xi'}\int_{\R}\frac{\hat{O}^{*}(\xi,k)}{Z+i\omega'(k)\xi}d\xi d\xi' dk,
	%\label{Ilimit}
\end{align*}
under the same limit with \eqref{limitballistic}. 
\subsubsection{The second scattering term \eqref{II}}\label{scat2re}
The second scattering term is  obtained in \eqref{Sepsilon} using $II_{\varepsilon}$. Similar to the treatment of first scattering term, the computations in Section \ref{IIscatterpoint} can be used to get the limit
\begin{align*}
	\frac{\nu}{4}&\int_{\R\times\T}\frac{r_a(k)\hat{W}_{0}(\xi',k)}{Z+i\omega'(k)\xi'}\int_{\R}\frac{\hat{O}^{*}(\xi,k)}{Z+i\omega'(k)\xi}d\xi d\xi' dk\\
	&\qquad+\frac{\nu}{4}\int_{\R\times\T}\frac{r_a(k)\hat{W}_{0}(\xi',-k)}{Z-i\omega'(k)\xi'}\int_{\R}\frac{\hat{O}^{*}(\xi,k)}{Z+i\omega'(k)\xi}d\xi d\xi' dk,
	%\label{IIlimit}
\end{align*}
under the same limit with \eqref{limitballistic}. 

The sum of the first and second scattering terms is
\begin{align*}
	&\int_{\R\times\T}(r_t(k)-1)\int_{\R}\frac{\hat{W}_0(\xi',k)}{Z+i\omega'(k)\xi'}\frac{\hat{O}^{*}(\xi,k)|v_g(k)|}{Z+i\omega'(k)\xi}d\xi' d\xi dk%\label{transmissionterm}
	\\
	&+\int_{\R\times\T}(r_r(k))\int_{\R}\frac{\hat{W}_0(\xi',-k)}{Z-i\omega'(k)\xi'}\frac{\hat{O}^{*}(\xi,k)|v_g(k)|}{Z+i\omega'(k)\xi}d\xi' d\xi dk,%\label{reflectionterm}
\end{align*}
under the same limit with \eqref{limitballistic}. 

\subsubsection{The fourth scattering terms \eqref{IV}}\label{scat4re}
The fourth scattering terms are  obtained in \eqref{Sepsilon} using $IV_{\varepsilon}$. We see that $-\frac{\nu}{2}(\cL(IV_{\varepsilon}^{(\varepsilon)}(\cdot,k-\varepsilon\xi/2))+\cL(IV_{\varepsilon}^{(\varepsilon)*}(\cdot,k+\varepsilon\xi/2))$ is the same as the term \eqref{term3deriLaplace} but with a different sign. Therefore, the fourth scattering term and term \eqref{term3deriLaplace} cancel each other. The limit is $0$.

	\subsubsection{The fifth scattering terms \eqref{V}}\label{scat5re}
The fifth scattering terms are  obtained in \eqref{Sepsilon} using $V_{\varepsilon}$. We have
\begin{align*}
	\cL(V_{\varepsilon}^{(\varepsilon)})(Z,k-\varepsilon\xi/2)&= \varepsilon\cL(V_{\varepsilon})(\varepsilon Z,k-\varepsilon\xi/2)\\
	&= \frac{-i\varepsilon}{\nu}\cL(\Theta\star F)(\varepsilon Z-i\omega(k-\varepsilon\xi/2))P(\varepsilon,k-\varepsilon\xi/2)\\
	&= \frac{-i\varepsilon}{\nu}\tilde{\Theta}(\varepsilon Z-i\omega(k-\varepsilon\xi/2))\cL(F)(\varepsilon Z-i\omega(k-\varepsilon\xi/2))P(\varepsilon,k-\varepsilon\xi/2).
	%\label{Vcalculate}
\end{align*}
By \ref{Passumption}, \ref{Fassumption} and \eqref{thetadef}, the limit is $0$.

\subsubsection{The third and the sixth scattering terms \eqref{III}, \eqref{VI}}\label{scat36re}
We compute $III_{\varepsilon}$ and $VI_{\varepsilon}$ in \eqref{Sepsilon}
\begin{align*}
	&-\frac{\nu}{2}\int_{\R\times\T}\cL\left( III_{\varepsilon}^{(\varepsilon)}(\cdot,k-\varepsilon\xi/2)+III_{\varepsilon}^{(\varepsilon)*}(\cdot,k+\varepsilon\xi/2) \right)(Z)\frac{\hat{O}^{*}(\xi,k)}{Z+iD_{\varepsilon}(k,\xi)}d\xi dk\\
	&-\frac{\nu}{2}\int_{\R\times\T}\cL\left( VI_{\varepsilon}^{(\varepsilon)}(\cdot,k-\varepsilon\xi/2)+VI_{\varepsilon}^{(\varepsilon)*}(\cdot,k+\varepsilon\xi/2) \right)(Z)\frac{\hat{O}^{*}(\xi,k)}{Z+iD_{\varepsilon}(k,\xi)}d\xi dk\\
	&= -\frac{\nu}{2}\left[\int_{\R\times\T}\cL\left( III_{\varepsilon}^{(\varepsilon)}(\cdot,k)\right)(Z)\frac{\hat{O}^{*}(\xi,k+\varepsilon\xi/2)}{Z+iD_{\varepsilon}^{+}(k,\xi)}d\xi dk+\int_{\R\times\T}\cL\left(III_{\varepsilon}^{(\varepsilon)*}(\cdot,k) \right)(Z)\frac{\hat{O}^{*}(\xi,k-\varepsilon\xi/2)}{Z+iD_{\varepsilon}^{-}(k,\xi)}d\xi dk\right]\\
	&\quad-\frac{\nu}{2}\left[\int_{\R\times\T}\cL\left( VI_{\varepsilon}^{(\varepsilon)}(\cdot,k)\right)(Z)\frac{\hat{O}^{*}(\xi,k+\varepsilon\xi/2)}{Z+iD_{\varepsilon}^{+}(k,\xi)}d\xi dk+\int_{\R\times\T}\cL\left(VI_{\varepsilon}^{(\varepsilon)*}(\cdot,k) \right)(Z)\frac{\hat{O}^{*}(\xi,k-\varepsilon\xi/2)}{Z+iD_{\varepsilon}^{-}(k,\xi)}d\xi dk\right]\\
	&= -\frac{\nu}{2}\left[\int_{\R\times\T}S_{\varepsilon}^{III}(Z,k)\frac{\hat{O}^{*}(\xi,k+\varepsilon\xi/2)}{Z+iD_{\varepsilon}^{+}(k,\xi)}d\xi dk+\int_{\R\times\T}S_{\varepsilon}^{III*}(Z,k)\frac{\hat{O}^{*}(\xi,k-\varepsilon\xi/2)}{Z+iD_{\varepsilon}^{-}(k,\xi)}d\xi dk\right]%\label{S3}
	\\
	&\quad -\frac{\nu}{2}\left[\int_{\R\times\T}S_{\varepsilon}^{VI}(Z,k)\frac{\hat{O}^{*}(\xi,k+\varepsilon\xi/2)}{Z+iD_{\varepsilon}^{+}(k,\xi)}d\xi dk+\int_{\R\times\T}S_{\varepsilon}^{VI*}(Z,k)\frac{\hat{O}^{*}(\xi,k-\varepsilon\xi/2)}{Z+iD_{\varepsilon}^{-}S(k,\xi)}d\xi dk\right].%\label{S6}
\end{align*}
Here, we use the notations
\begin{align}
	S_{\varepsilon}^{III}(Z,k):=\int_{0}^{\infty}\varepsilon e^{-Z\varepsilon t}III_{\varepsilon}(t,k)dt,\label{S3def}\\
	S_{\varepsilon}^{III*}(Z,k):=\int_{0}^{\infty}\varepsilon e^{-Z\varepsilon t}III_{\varepsilon}^{*}(t,k)dt,\label{S3star}\\
	S_{\varepsilon}^{VI}(Z,k):=\int_{0}^{\infty}\varepsilon e^{-Z\varepsilon t}VI_{\varepsilon}(t,k)dt,\label{S6def}\\
	S_{\varepsilon}^{VI*}(Z,k):=\int_{0}^{\infty}\varepsilon e^{-Z\varepsilon t}VI_{\varepsilon}^{*}(t,k)dt.\label{S6star}
\end{align}
Let us compute \eqref{S3def} $+$ \eqref{S6def}. The sum \eqref{S3star} $+$ \eqref{S6star} is similar.

To this end, we have 
\begin{align*}
	S^{III}_{\varepsilon}(Z,k)&= \varepsilon\int_{0}^{\infty}e^{-Z\varepsilon t}III_{\varepsilon}(t,k)dt\\
	&= \varepsilon\int_{0}^{\infty}e^{-Z\varepsilon t}\int_{0}^{t}e^{i\omega(k)(t-s)}\Theta\star F(s)\langle \Theta\star \Psi(t)\rangle_{\mu_{\varepsilon}}dsdt.%\label{S3expand}
\end{align*}
Hence
\begin{align*}
	\Real S^{III}_{\varepsilon}(Z,k)&= \varepsilon\int_{0}^{\infty}e^{-Z\varepsilon t}\int_{0}^{t}\cos(\omega(k)s)\Theta\star F(s)\cos(\omega(k)t)\langle \Theta\star \Psi(t)\rangle_{\mu_{\varepsilon}}dsdt%\label{S3realcos}
	\\
	&\quad+\varepsilon\int_{0}^{\infty}e^{-Z\varepsilon t}\int_{0}^{t}\sin(\omega(k)s)\Theta\star F(s)\sin(\omega(k)t)\langle \Theta\star \Psi(t)\rangle_{\mu_{\varepsilon}}dsdt.%\label{S3realsin}
\end{align*}
Proceeding similarly with $VI_{\varepsilon}$ we get
\begin{align*}
	\Real S^{VI}_{\varepsilon}(Z,k)&= \varepsilon\int_{0}^{\infty}e^{-Z\varepsilon t}\cos(\omega(k)t)\Theta\star F(t)\left\langle\int_{0}^{t}\cos(\omega(k)s) \Theta\star \Psi(s)ds\right\rangle_{\mu_{\varepsilon}}dt%\label{S6realcos}
	\\
	&\quad+\varepsilon\int_{0}^{\infty}e^{-Z\varepsilon t}\sin(\omega(k)t)\Theta\star F(t)\left\langle\int_{0}^{t}\sin(\omega(k)s) \Theta\star \Psi(s)ds\right\rangle_{\mu_{\varepsilon}}dt.%\label{S6realsin}
\end{align*}
Therefore, using integrations by parts, we have
\begin{align}
	&\Real S^{III}_{\varepsilon}(Z,k)+\Real S^{VI}_{\varepsilon}(Z,k)\nonumber\\
	&= \varepsilon^2Z\int_{0}^{\infty}e^{-Z\varepsilon t}\left( \int_{0}^{t}\cos(\omega(k)s)\Theta\star F(s)ds \right)\left\langle \int_{0}^{t}\cos(\omega(k)s)\Theta\star \Psi(s)ds \right\rangle_{\mu_{\varepsilon}}dt\label{sumrealcos}\\
	&+\varepsilon^2Z\int_{0}^{\infty}e^{-Z\varepsilon t}\left( \int_{0}^{t}\sin(\omega(k)s)\Theta\star F(s)ds \right)\left\langle \int_{0}^{t}\sin(\omega(k)s)\Theta\star \Psi(s)ds \right\rangle_{\mu_{\varepsilon}}dt.\nonumber%\label{sumrealsin}
\end{align}
Using the identity
\begin{align*}
	\Psi(t)=\frac{1}{2i}\left( \int_{\T}e^{-i\omega(k)t}\hat{\psi}(0,k)-\int_{\T}e^{i\omega(k)t}\hat{\psi}^{*}(0,k) \right),%\label{Psiexpand}
\end{align*}
we expand the term \eqref{sumrealcos} as
\begin{align}
	&\frac{\varepsilon^{2}Z}{2i}\int_{0}^{\infty}e^{-Z\varepsilon t}\int_{0}^{t}\int_{0}^{t}\cos(\omega(k)s)\cos(\omega(k)s')\nonumber\\
	&\qquad \times \int_{0}^{s}\int_{0}^{s'}F(s-u)\left\langle \int_{\T}e^{-i\omega(h)(s'-u')}\hat{\psi}(0,h)dh\right\rangle_{\mu_{\varepsilon}}d\Theta(u)d\Theta(u')dsds'dt\label{sumrealcos1}\\
	&-\frac{\varepsilon^{2}Z}{2i}\int_{0}^{\infty}e^{-Z\varepsilon t}\int_{0}^{t}\int_{0}^{t}\cos(\omega(k)s)\cos(\omega(k)s')\nonumber\\
	&\qquad \times \int_{0}^{s}\int_{0}^{s'}F(s-u)\left\langle \int_{\T} e^{i\omega(h)(s'-u')}\hat{\psi}^{*}(0,h)dh\right\rangle_{\mu_{\varepsilon}}d\Theta(u)d\Theta(u')dsds'dt.\label{sumrealcos2}
\end{align}
Let us estimate \eqref{sumrealcos1}, the term \eqref{sumrealcos2} can be estimated similarly. After that, we subtract \eqref{sumrealcos1} $ and $ \eqref{sumrealcos2} to get \eqref{sumrealcos}. Using the identity 
\begin{align*}
	\delta(t-t')=\frac{1}{2\pi}\int_{\R} e^{ia(t'-t)}da,%\label{diracdelta}
\end{align*}
and changing $e^{-\lambda\varepsilon t}$ into $e^{-\lambda\varepsilon(t+t')/2}$, \eqref{sumrealcos1} becomes
\begin{align*}
	\frac{\varepsilon^{2}Z}{2i}\int_{0}^{\infty}e^{-Z\varepsilon (t+t')/2}&\int_{0}^{t}\int_{0}^{t'}\cos(\omega(k)s)\cos(\omega(k)s')\\
	&\qquad\int_{0}^{s}\int_{0}^{s'}F(s-u)\left\langle \int_{\T}e^{-i\omega(h)(s'-u')}\hat{\psi}(0,h)dh\right\rangle_{\mu_{\varepsilon}}d\Theta(u)d\Theta(u')dsds'\delta(t-t')dtdt'\\
	=\frac{\varepsilon^{2}Z}{4i\pi}\int_{\R} &\int_{0}^{\infty}\int_{0}^{\infty}e^{ia(t'-t)}e^{-Z\varepsilon (t+t')/2}\int_{0}^{t}\int_{0}^{t'}\cos(\omega(k)s)\cos(\omega(k)s')\\
	&\int_{0}^{s}\int_{0}^{s'}F(s-u)\left\langle \int_{\T} e^{-i\omega(h)(s'-u')}\hat{\psi}(0,h)dh\right\rangle_{\mu_{\varepsilon}}d\Theta(u)d\Theta(u')dsds'dtdt'da.%\label{Applydiracdelta}
\end{align*}

By rewriting the domain of integration, 
\begin{align*}
	(u,s,t)\in[0,s]\times[0,t]\times[0,\infty]\to(t,s,u)\in[s,\infty]\times[u,\infty]\times[0,\infty],%\label{newintegration}
\end{align*}
this term is equal to
\begin{align*}
	\frac{\varepsilon Z}{4i\pi}\int_{\R}\int_{\T}\varepsilon P(\varepsilon,h)A(h,k,Z)A_{F}(k,Z)dhda,%\label{Aterm}
\end{align*}
where
\begin{align*}
	A(h,k,Z):=\frac{1}{2(Z\varepsilon/2-ia)}\left( \frac{\tilde{\Theta}(Z\varepsilon/2-ia-i\omega(k))}{Z\varepsilon/2-ia-i(\omega(k)-\omega(h))}+\frac{\tilde{\Theta}(Z\varepsilon/2-ia+i\omega(k))}{Z\varepsilon/2-ia+i(\omega(k)+\omega(h))} \right),%\label{A}
\end{align*}
and
\begin{align*}
	A_F(k,Z)&:= \int_{0}^{\infty}\int_{u}^{\infty}\cos(\omega(k)s)F(s-u)\int_{s}^{\infty}e^{-(Z\varepsilon/2+ia)t}dtdsd\Theta(u)\\
	&= \frac{1}{Z\varepsilon/2+ia}\int_{0}^{\infty}\int_{u}^{\infty}\cos(\omega(k)s)F(s-u)e^{-(Z\varepsilon/2+ia)s}dsd\Theta(u)\\
	&= \frac{1}{2(Z\varepsilon/2+ia)}\int_{0}^{\infty}\left( e^{-(Z\varepsilon/2+ia-i\omega(k))u}\cL(F)(Z\varepsilon/2+ia-i\omega(k))d\Theta(u) \right.\\
	&\qquad \left.+e^{-(Z\varepsilon/2+ia+i\omega(k))u}\cL(F)(Z\varepsilon/2+ia+i\omega(k))\right)\\
	&= \frac{1}{2(Z\varepsilon/2+ia)}\left( \tilde{\Theta}(Z\varepsilon/2+ia-i\omega(k))\cL(F)(Z\varepsilon/2+ia-i\omega(k)) \right.%\label{AFminus}
	\\
	&\qquad \left.+\tilde{\Theta}(Z\varepsilon/2+ia+i\omega(k))\cL(F)(Z\varepsilon/2+ia+i\omega(k))\right).%\label{AFplus}
\end{align*}

Making the change variable $a\to\varepsilon a$, \eqref{sumrealcos1} is equal to
\begin{align}
	&\frac{Z}{16i\pi}\int_{\R}\int_{\T}P(\varepsilon,h)\frac{1}{(\lambda/2)^2+a^2}\nonumber\\
	&\times \left\{ \tilde{\Theta}(\varepsilon(Z/2+ia)-i\omega(k))\cL(F)(\varepsilon(Z/2+ia)-i\omega(k))\right.\nonumber\\
	&+\left.\tilde{\Theta}(\varepsilon(Z/2+ia)+i\omega(k))\cL(F)(\varepsilon(Z/2+ia)+i\omega(k)) \right\}\frac{\tilde{\Theta}(\varepsilon(Z/2-ia)-i\omega(k))}{Z/2-ia-i\varepsilon^{-1}(\omega(k)-\omega(h))}dhda\label{minusdenom}\\
	&+\frac{Z}{16i\pi}\int_{\R}\int_{\T} P(\varepsilon,h)\frac{1}{(Z/2)^2+a^2}\nonumber\\
	&\times \left\{ \tilde{\Theta}(\varepsilon(Z/2+ia)-i\omega(k))\cL(F)(\varepsilon(Z/2+ia)-i\omega(k))\right.\nonumber\\
	&+\left.\tilde{\Theta}(\varepsilon(Z/2+ia)+i\omega(k))\cL(F)(\varepsilon(Z/2+ia)+i\omega(k)) \right\}\frac{\varepsilon\tilde{\Theta}(\varepsilon(Z/2-ia)+i\omega(k))}{\varepsilon(Z/2-ia)+i(\omega(k)+\omega(h))}dhda.\label{plusdenom}
\end{align}
We expect that most of the contribution at the limit comes from \eqref{minusdenom}. We will prove that the limit of \eqref{plusdenom} is zero. From \ref{Passumption}, \ref{Fassumption}, \eqref{thetadef}, using the bound of test function $\hat{O}^{*}$, and the bound $|Z+iD_{\varepsilon}^{+}(k,\eta)|\ge\Real(Z)$, it suffices to estimate
\begin{align}
	\left|\frac{\varepsilon^{1-d}}{\varepsilon(Z/2-ia)+i\omega(k)+i\omega(h)}\right|\lesssim \varepsilon^{1-d}.
	\label{estimationofplus}
\end{align}
For $\varepsilon$ sufficient small, $|\varepsilon(Z/2-ia)+i\omega(k)+i\omega(h)|\gtrsim 2\omega_{min}>0$, and \eqref{estimationofplus} follows.

Next, we handle the term \eqref{minusdenom}.
If $|\omega(k)-\omega(h)|\gtrsim\varepsilon^{1/3}$, then we do the same as \eqref{estimationofplus} to show that the integral on this domain is small. We consider the domain $|\omega(k)-\omega(h)|\le \varrho_{1}\varepsilon^{1/3}$ for some $\varrho_{1}$ to be specified later.

Recall that we consider the \eqref{Trunass} where $\varrho$ to be chosen later. In the domain, $k\in\left[ -\frac{1}{2},-\frac{1}{2}+\frac{\varrho}{2} \right]\cup\left[ -\frac{\varrho}{2},\frac{\varrho}{2} \right]\cup\left[ \frac{1}{2}-\frac{\varrho}{2},\frac{1}{2} \right]$, we choose $\varrho_{1}\ll\frac{\varrho^2}{4}$. Then, $h\in [\T](\varrho)$ according to \eqref{minapprox} and \eqref{maxapprox}. As a result, $P(\varepsilon,h)=0$ on this domain.

We will focus on: $k\in\left[ -\frac{1}{2}+\frac{\varrho}{2},-\frac{\varrho}{2} \right]\cup \left[ \frac{\varrho}{2},\frac{1}{2}-\frac{\varrho}{2} \right]$ and $h$ such that $|\omega(k)-\omega(h)|\le \varrho_{1}\varepsilon^{1/3}$ for some small $\varrho_{1}$. Because $|\omega(k)-\omega(h)|=|\omega'(k')||k-h|\gtrsim |k-h|$, we consider $h-k\in[-\varrho_{1}\varepsilon^{1/3},\varrho_{1}\varepsilon^{1/3}]$ for some small $\varrho_{1}$.

We perform the change of variables $h\to k-\varepsilon\xi'$ where $\xi'\in [-\varrho_{1}\varepsilon^{-2/3},\varrho_{1}\varepsilon^{-2/3}]$. \eqref{minusdenom} becomes
\begin{align}
	&\frac{-Z\varepsilon}{16i\pi}\int_{\R}\int_{-\varrho_{1}\varepsilon^{-2/3}}^{\varrho_{1}\varepsilon^{-2/3}} P(\varepsilon,k-\varepsilon\xi')\frac{1}{(Z/2)^2+a^2}\nonumber\\
	&\times \left\{ \tilde{\Theta}(\varepsilon(Z/2+ia)-i\omega(k))\cL(F)(\varepsilon(Z/2+ia)-i\omega(k))\right.\nonumber\\
	&+\left.\tilde{\Theta}(\varepsilon(Z/2+ia)+i\omega(k))\cL(F)(\varepsilon(Z/2+ia)+i\omega(k)) \right\}\frac{\tilde{\Theta}(\varepsilon(Z/2-ia)-i\omega(k))}{Z/2-ia-i\varepsilon^{-1}(\omega(k)-\omega(k-\varepsilon\xi'))}d\xi'da.\label{minusdenometa}
\end{align}

We approximate \eqref{minusdenometa} using the derivative of $\omega$, that is
\begin{align}
	&\frac{-Z\varepsilon}{16i\pi}\int_{\R}\int_{-\varrho_{1}\varepsilon^{-2/3}}^{\varrho_{1}\varepsilon^{-2/3}} P(\varepsilon,k-\varepsilon\xi')\frac{1}{(Z/2)^2+a^2}\nonumber\\
	&\times \left\{ \tilde{\Theta}(\varepsilon(Z/2+ia)-i\omega(k))\cL(F)(\varepsilon(Z/2+ia)-i\omega(k))\right.\nonumber\\
	&+\left.\tilde{\Theta}(\varepsilon(Z/2+ia)+i\omega(k))\cL(F)(\varepsilon(Z/2+ia)+i\omega(k)) \right\}\frac{\tilde{\Theta}(\varepsilon(Z/2-ia)-i\omega(k))}{Z/2-ia-i\omega'(k)\xi'}d\xi'da.\label{minusdenomderi}
\end{align}
We have 
\begin{align*}
	&\left|\frac{\varepsilon}{Z/2-ia-i\varepsilon^{-1}(\omega(k)-\omega(k-\varepsilon\xi'))}-\frac{\varepsilon}{Z/2-ia-\omega'(k)\xi'}\right|\\
	&= \left|\frac{\omega(k)-\omega(k-\varepsilon\xi')+\omega'(k)\varepsilon\xi'}{(Z/2-ia-i\varepsilon^{-1}(\omega(k)-\omega(k-\varepsilon\xi')))(Z/2-ia-i\omega'(k)\xi')}\right|,
	%\label{differenceofdenom}
\end{align*}
and $|\omega(k)-\omega(k-\varepsilon\xi')+\omega'(k)\varepsilon\xi'|\approx|\omega''(k')(\varepsilon\xi')^{2}|\lesssim\varrho_{1}^{2}\varepsilon^{2/3}$. Besides, by Cauchy-Schwarz's inequality, we have
\begin{align*}
	&\left|\frac{1}{(Z/2-ia-i\varepsilon^{-1}(\omega(k)-\omega(k-\varepsilon \xi')))(Z/2-ia-i\omega'(k)\xi')}\right|\\
	&\lesssim\frac{1}{|Z/2-ia-i\varepsilon^{-1}(\omega(k)-\omega(k-\varepsilon\xi'))|^{2}}+\frac{1}{|Z/2-ia-i\omega'(k)\xi'|^{2}}\\
	&\approx \frac{1}{1+|a+\varepsilon^{-1}(\omega(k)-\omega(k-\varepsilon\xi'))|^{2}}+\frac{1}{1+|a+\omega'(k)\xi'|^{2}}.
	%\label{CSestimate}
\end{align*}
Next, by the estimate
\begin{align}
	\int_{\R}\frac{da}{(1+a^2)(1+(a+b)^2)}\lesssim\frac{1}{1+b^2},
	\label{asquareestimate}
\end{align}
the difference \eqref{minusdenometa} $-$ \eqref{minusdenomderi} is estimated by
\begin{align*}
	\frac{\varrho_{1}^{2}\varepsilon^{1/6-d}}{\inf_{k\in[\varrho/2-\varrho_{1},1/2-\varrho/2+\varrho_1]}|\omega'(k)|}\int_{\R}\frac{dx}{1+x^2}\lesssim\varrho_{1}^{2}\varepsilon^{1/6-d}.
	%\label{estimatedelta1}
\end{align*}
In the latter estimate, we recall that $d\in[0,1/6)$ so $\varepsilon^{1/6-d}$ tends to $0$ as $\varepsilon\to0$.

We now estimate \eqref{minusdenomderi} for some small $\varrho_{1}$.
Using \ref{Passumption}, \ref{Fassumption} and \eqref{thetadef}, we only need to estimate
\begin{align}
	\varepsilon^{1/2-d}\int_{-\varrho_{1}/\varepsilon^{1/3}}^{\varrho_{1}/\varepsilon^{1/3}}\frac{d\xi'}{|Z/2-ia-i\omega'(k)\xi'|}&= \varepsilon^{1/2-d}\int_{-\varrho_1\varepsilon^{-2/3}}^{\varrho_1\varepsilon^{-2/3}}\frac{d\xi'}{\sqrt{\Real(Z/2)^2+(\omega'(k)\xi'+a-\Imagine(Z/2))^2}}\nonumber\\
	&\lesssim\frac{\varepsilon^{1/2-d}}{\omega'(k)}\arcsinh(x)\bigg|_{-\omega'(k)\varrho_1\varepsilon^{-2/3}+a-\Imagine(Z/2)}^{\omega'(k)\varrho_1\varepsilon^{-2/3}+a-\Imagine(Z/2)}\approx \varepsilon^{1/2-d}\log(\varepsilon^{-2/3}).
	\label{estimatederi}
\end{align}
The limit of term \eqref{estimatederi} is $0$. This implies that the term \eqref{minusdenom} also has limit $0$.

\subsubsection{The seventh scattering terms \eqref{VII}}\label{scat7re}
The seventh scattering terms are obtained in \eqref{Sepsilon} using $VII_{\varepsilon}$. We have a transformation for $VII_{\varepsilon}(t,k)$ using \eqref{Theta^omega} and \eqref{Fk}
\begin{align*}
	VII_{\varepsilon}(t,k)&= \frac{-1}{\nu}\int_{0}^{t}F(t-s')d\Theta(s')\int_{0}^{t}\Theta_\omega^{*}(t-s,k)F(s)ds\\
	&= \frac{-1}{\nu}\int_{0}^{t}e^{i\omega(k)s'}e^{i\omega(k)(t-s')}F(t-s')d\Theta(s')\int_{0}^{t}(\Theta^\omega)^{*}(t-s,k)e^{-i\omega(k)s}F(s)ds\\
	&= \frac{-1}{\nu}((e^{i\omega(k)\cdot}\Theta)\star F_{k}^{*}(t))( (\Theta^\omega)^{*}\star F_{k}(t)).%\label{VIItransform}
\end{align*}
Next, we compute $\cL(VII_{\varepsilon}^{(\varepsilon)})(Z,k-\varepsilon\xi/2)$ using the formula for the Laplace transform of a product
\begin{align}
	&\cL(VII_{\varepsilon}^{(\varepsilon)})(Z,k-\varepsilon\xi/2)\nonumber\\
	&= \frac{-1}{2\nu i\pi}\lim_{L\to\infty}\int_{c-iL}^{c+iL}\cL((e^{i\omega(\cdot)\cdot}\Theta)\star F_{\cdot}^{*}{}^{(\varepsilon)})(\sigma,k-\varepsilon\xi/2)\cL( (\Theta^\omega)^{*}\star F_{\cdot}{}^{(\varepsilon)})(Z-\zeta,k-\varepsilon\xi/2)d\zeta\nonumber\\
	&= \frac{-\varepsilon^2}{2\nu i\pi}\lim_{L\to\infty}\int_{c-iL}^{c+iL}\tilde{\Theta}(\varepsilon\zeta-i\omega(k-\varepsilon\xi/2))\cL(F)(\varepsilon\zeta-i\omega(k-\varepsilon\xi/2))\nonumber\\
	&\qquad\cL( (\Theta^\omega)^{*})(\varepsilon(Z-\zeta),k-\varepsilon\xi/2)\cL(F)(\varepsilon(Z-\zeta)+i\omega(k-\varepsilon\xi/2))d\zeta\nonumber\\
	&= \frac{-\varepsilon}{2\nu i\pi}\lim_{L\to\infty}\int_{c-iL}^{c+iL}\tilde{\Theta}(\varepsilon\zeta-i\omega(k-\varepsilon\xi/2))\cL(F)(\varepsilon\zeta-i\omega(k-\varepsilon\xi/2))\nonumber\\
	&\qquad\cL( (\Theta^\omega)^{(\varepsilon)*})(Z-\zeta,k-\varepsilon\xi/2)\cL(F)(\varepsilon(Z-\zeta)+i\omega(k-\varepsilon\xi/2))d\zeta\nonumber\\
	&= \frac{-\varepsilon}{2\nu i\pi}\lim_{L\to\infty}\int_{c-iL}^{c+iL}\tilde{\Theta}(\varepsilon\zeta-i\omega(k-\varepsilon\xi/2))\cL(F)(\varepsilon\zeta-i\omega(k-\varepsilon\xi/2))\nonumber\\
	&\qquad\frac{\tilde{\Theta}(\varepsilon(Z-\zeta)+i\omega(k-\varepsilon\xi/2))}{Z-\zeta}\cL(F)(\varepsilon(Z-\zeta)+i\omega(k-\varepsilon\xi/2))d\zeta.\label{VIIcalculate}
\end{align}
At the limit $\varepsilon\to0$, the term \eqref{VIIcalculate} gives  
\begin{align*}
	\frac{-1}{2\nu i\pi}\lim_{L\to\infty}\int_{c-iL}^{c+iL}\frac{|\theta(k)|^2|\cF(k)|^2}{Z-\zeta}d\zeta.
	%\label{VIIlimit}
\end{align*}
Now, the integral is estimated using \eqref{integralfraction}. We get
\begin{align}
	\frac{-|\theta(k)|^2|\cF(k)|^2}{2\nu}.
	\label{VIIlimitresult}
\end{align}
Expanding \eqref{VIIlimitresult} in \eqref{Sepsilon}, the result is
\begin{align*}
	\frac{|\theta(k)|^2|\cF(k)|^2}{2(Z+i\omega'(k)\xi)}.
	%\label{sumVII}
\end{align*}

\subsection{Proof of Theorem \ref{Main}}\label{proofmain}

To prove the theorem, we first show that the kinetic limit of $\psi^1$ in \eqref{multiplybump} and the kinetic limit of $\psi$ in \eqref{dpsi} give the same result. The growth of energy of $\hat{\psi}^{1}$ is at most linear, which is similar to $\hat{\psi}$. Meanwhile, the energy of $\hat{\psi}^{2}$ is decreasing.  We use the bound on the initial condition \ref{energyassumption} and get the following estimate
\begin{align*}
	&\limsup_{\varepsilon\to0^{+}}\sup_{\xi\in\T_{2/\varepsilon}}\int_{\T}|\sW_{\varepsilon}(Z,\xi,k)-\sW_{\varepsilon}^{1}(Z,\xi,k)|dk\\
&= \limsup_{\varepsilon\to0^{+}}\sup_{\xi\in\T_{2/\varepsilon}}\frac{\varepsilon}{2}\E_{\varepsilon}\left[ \int_{\T}\left|\int_{0}^{\infty}e^{-Zt}\left( \hat{\psi}^{(\varepsilon)}(t,k+\frac{\varepsilon\xi}{2})\hat{\psi}^{(\varepsilon)*}(t,k-\frac{\varepsilon\xi}{2})\right.\right.\right.\\
&\hspace{15em}\left.\left.\left.-\hat{\psi}^{1(\varepsilon)}(t,k+\frac{\varepsilon\xi}{2})\hat{\psi}^{1(\varepsilon)*}(t,k-\frac{\varepsilon\xi}{2}) \right)dt\right|dk \right]\\
&\le \limsup_{\varepsilon\to0^{+}}\sup_{\xi\in\T_{2/\varepsilon}}\frac{\varepsilon}{2}\E_{\varepsilon}\left[ \int_{0}^{\infty}e^{-\Real Zt}\int_{\T}\left|\left( \hat{\psi}^{(\varepsilon)}(t,k+\frac{\varepsilon\xi}{2})\hat{\psi}^{(\varepsilon)*}(t,k-\frac{\varepsilon\xi}{2})\right.\right.\right.\\
&\hspace{15em}\left.\left.\left.-\hat{\psi}^{1(\varepsilon)}(t,k+\frac{\varepsilon\xi}{2})\hat{\psi}^{1(\varepsilon)*}(t,k-\frac{\varepsilon\xi}{2}) \right)dt\right|dk \right]\\	
	&\le \limsup_{\varepsilon\to0^{+}}\sup_{\xi\in\T_{2/\varepsilon}}\frac{\varepsilon}{2}\E_{\varepsilon}\left[ \int_{0}^{\infty}e^{-\Real Zt}\int_{\T}\left| \hat{\psi}^{2(\varepsilon)}(t,k+\frac{\varepsilon\xi}{2})\hat{\psi}^{2(\varepsilon)*}(t,k-\frac{\varepsilon\xi}{2})dt \right|dk \right]\\
	&+\limsup_{\varepsilon\to0^{+}}\sup_{\xi\in\T_{2/\varepsilon}}\varepsilon\E_{\varepsilon}\left[ \int_{0}^{\infty}e^{-\Real Zt}\int_{\T}\left| \hat{\psi}^{2(\varepsilon)}(t,k+\frac{\varepsilon\xi}{2})\hat{\psi}^{1(\varepsilon)*}(t,k-\frac{\varepsilon\xi}{2})dt \right|dk \right]\\
	&\lesssim \int_{0}^{\infty}e^{-\Real Zt}\limsup_{\varepsilon\to0^{+}}\varepsilon\E_{\varepsilon}\left\|\hat{\psi}^{2}(t/\varepsilon)\right\|_{L^{2}(\T)}\left( \left\|\hat{\psi}^{2}(t/\varepsilon)\right\|_{L^{2}(\T)}+\left\|\hat{\psi}^{1}(t/\varepsilon)\right\|_{L^{2}(\T)} \right)dt\\
	&\lesssim \rho(\varrho)\left( \rho(\varrho)+\int_{0}^{t}e^{-\Real Zt}tdt \right)\ll 1.
	%\label{prooftrimming}
\end{align*}
In this estimate, $\rho(\varrho)^{2}$ is a bound for $\sup_{\varepsilon}\varepsilon\E_{\varepsilon}\|\hat{\psi}^{2}\|_{L^2}^{2}$. We see that the limit $\rho(\varrho)$ is $0$ as $\varrho$ tends to $0$ by Dominated Convergence Theorem. If we denote $\sW^{1}$ to be the limit of $\sW^{1}_{\varepsilon}$, then for $\hat{O}^{*}\in\sS(\R\times\T)$ we have
\begin{align*}
	\lim\sup_{\varepsilon\to0^{+}}&\left|\int_{\R\times\T}\hat{O}^{*}(\xi,k)\sW^{1}(Z,\xi,k)d\xi dk-\int_{\R\times\T}\hat{O}^{*}(\xi,k)\sW_{\varepsilon}(Z,\xi,k)d\xi dk\right|\\
	&\le\left|\int_{\R\times\T}\hat{O}^{*}(\xi,k)\sW_{\varepsilon}^{1}(Z,\xi,k)d\xi dk-\int_{\R\times\T}\hat{O}^{*}(\xi,k)\sW^{1}(Z,\xi,k)d\xi dk\right|\\
	&\quad
	+\left|\int_{\R\times\T}\hat{O}^{*}(\xi,k)\sW_{\varepsilon}^{1}(Z,\xi,k)d\xi dk-\int_{\R\times\T}\hat{O}^{*}(\xi,k)\sW_{\varepsilon}(Z,\xi,k)d\xi dk\right|\\
	&\lesssim \rho(\varrho).
	%\label{w1tow}
\end{align*}

With \ref{Trunass}, we combine all results in Section \ref{thermal1re}, \ref{thermal23re}, \ref{thermal4re}, \ref{scat1re}, \ref{scat2re}, \ref{scat4re}, \ref{scat5re}, \ref{scat36re} and \ref{scat7re} to get
\begin{equation}
	\begin{split}
		\sW^1(Z,\xi,k)=&  T|v_g(k)|r_a(k)\underbrace{\frac{1}{Z(Z+i\omega'(k))\xi}}_{(a)}+\underbrace{\frac{\hat{W}_{0}(\xi,k)}{Z+i\omega'(k)\xi}}_{(b)}+(r_t(k)-1)\underbrace{\int_{\R}\frac{|v_g(k)|\hat{W}_{0}(\xi',k)d\xi'}{(Z+i\omega'(k)\xi)(Z+i\omega'(k)\xi')}}_{(c)}\\
		&+r_r(k)\underbrace{\int_{\R}\frac{|v_g(k)|\hat{W}_{0}(\xi',-k)d\xi'}{(Z+i\omega'(k)\xi)(Z-i\omega'(k)\xi')}}_{(d)}+|v_g(k)|r_a(k)|\cF(k)|^2\underbrace{\frac{1}{\nu(Z+i\omega'(k)\xi)}}_{(e)}.
	\end{split}
	\label{sumallw}
\end{equation}
Hence, $\sW(Z,\xi,k)=\sW^1(Z,\xi,k)$ since $\sW^1$ does not depend on $\varrho$.

To get $W$ in the main theorem, we take the inverse Laplace and Fourier transform of each term. Let us make precise the technical computations as follows.

We first compute the inverse Laplace of term $(a)$ of \eqref{sumallw} with test function $O$
\begin{align}
	\int_{\R} \hat{O}^{*}(\xi,k)&\cL^{-1}\left( \frac{1}{Z(Z+i\omega'(k)\xi)} \right)(t)d\xi= \int_{\R}\int_{\R}\frac{1-e^{-i2\pi v_g(k)t\xi}}{i2\pi v_g(k)\eta}e^{i2\pi x\xi}O^{*}(x,k)dxd\xi\nonumber\\
	&=\int_{\R}O^{*}(x,k)\left( 1_{[0,\infty)}(x)-1_{[v_g(k)t,\infty)}(x) \right)dx=\int_{\R}O^{*}(x)1_{[0,v_g(k)t]}(x)dx.
	\label{thermal4limitshowing}
\end{align}

Term $(b)$ of \eqref{sumallw} can be computed as
\begin{align}
	\int_{\R} \hat{O}^{*}(\xi,k)&\cL^{-1}\left( \frac{\hat{W}_{0}(\xi,k)}{Z+i\omega'(k)\xi} \right)(t)d\xi= \int_{\R}\int_{\R} e^{-i2\pi v_g(k)t\xi}e^{i2\pi x\xi}O^{*}(x,k)\hat{W}_{0}(\xi,k)dxd\xi\nonumber\\
	&=\int_{\R}O^{*}(x,k)W_{0}(x-v_g(k)t,k)dx.
	\label{ballisticlimitshowing}
\end{align}

Term $(c)$ of \eqref{sumallw} is now computed
\begin{align}
	\int_{\R}&\int_{\R}\cL^{-1}\left(\frac{\hat{W}_{0}(\xi',k)\hat{O}^{*}(\xi,k)|v_g(k)|}{(Z+i\omega'(k)\xi')(Z+i\omega'(k)\xi)}\right)(t)d\xi'd\xi\nonumber\\
	&= \int_{\R}\int_{\R}e^{i2\pi x \xi}O^{*}(x,k)|v_g(k)|\int_{\R}\left( e^{-i2\pi v_g(k)\xi\cdot}\star(e^{-i2\pi v_g(k)\xi'}\hat{W}_{0}(\xi',k))(t) \right)d\xi'd\xi dx\nonumber\\
	&= \int_{\R}\int_{0}^{t}|v_g(k)|O^{*}(x,k)W_{0}(-v_g(k)(t-s),k)\delta(x-v_g(k)s)dsdx\nonumber\\
	&= \int_{\R}O^{*}(x,k)W_{0}(x-v_g(k)t,k)1_{[0,v_g(k)t]}(x)dx.
	\label{transmittinglimitshowing}
\end{align}

Now, we calculate term $(d)$
\begin{align}
	\int_{\R} &\int_{\R}\cL^{-1}\left(\frac{\hat{W}_{0}(\xi',-k)\hat{O}^{*}(\xi,k)|v_g(k)|}{(Z-i\omega'(k)\xi')(Z+i\omega'(k)\xi)}\right)(t)d\xi'd\xi\nonumber\\
	&= \int_{\R}\int_{\R} e^{i2\pi x \xi}O^{*}(x,k)|v_g(k)|\int_{\R}\left( e^{-i2\pi v_g(k)\xi\cdot}\star(e^{i2\pi v_g(k)\xi'}\hat{W}_{0}(\xi',-k))(t) \right)d\xi'd\xi dx\nonumber\\
	&= \int_{\R}\int_{0}^{t}|v_g(k)|O^{*}(x,k)W_{0}(v_g(k)(t-s),-k)\delta(x-v_g(k)t)dsdx\nonumber\\
	&= \int_{\R}O^{*}(x,k)W_{0}(v_g(k)t-x,-k)1_{[0,v_g(k)t]}(x)dx.
	\label{reflectinglimitshowing}
\end{align}

Finally, we develop term $(e)$ of \eqref{sumallw}
\begin{align*}
	\int_{\R}\hat{O}^{*}(\xi,k)&\cL^{-1}\left( \frac{1}{Z+i\omega'(k)\xi} \right)(t)d\xi= \int_{\R}\int_{\R} e^{i2\pi(x-v_g(k)t)\xi}O^{*}(x,k)d\xi dx\\
	&=\int_{\R}O^{*}(x,k)\delta(x-v_g(k)t)dx.
	%\label{thermal1limitshowing}
\end{align*}
We have proved that
\begin{align*}
	\lim_{\varepsilon}\left\langle O,W_{\varepsilon}(t)\right\rangle_{\sN}&= \int_{\R\times\T} O^{*}(x,k)W(t,x,k)dkdx
\end{align*}
for any $t$ and any $O\in\sS(\R\times\T)$. The conclusion of the theorem follows.

\section{Linear memory-feedback control - Proof of Theorem \ref{mainnew} \& Corollary \ref{controlcoro}}
\label{LMFcontrol}
\subsection{Preliminary computations}

In this Section, we focus on the system defined by \eqref{dbetanew} and \eqref{dalphanew}. Assumption \ref{inequalitycondition1} allows us to keep the control over the energy of the wave. We have the lemma.
\begin{lemma}
	Consider the wave system governed by \eqref{dbetanew} and \eqref{dalphanew} with all the assumptions \ref{init0} - \ref{W_0}, \ref{expdecay} - \ref{positivesigma}, \ref{monoome} - \ref{inverome} and \ref{inequalitycondition1}. There is a constant $C_6$ such that
	\begin{equation*}
		\varepsilon\E_{\varepsilon}\|\hat{\psi}^{(\varepsilon)}(t)\|_{L^2(\T)}^2\le\varepsilon\E_{\varepsilon}\|\hat{\psi}^{(\varepsilon)}(0)\|_{L^2(\T)}^2+C_6t.
	\end{equation*}
	\label{linearfeed}
\end{lemma}
Similar to Lemma \ref{ligrowimpul}, we use an It\^{o} analysis for the proof of Lemma \ref{linearfeed}, the details will be given in Section \ref{prolinearfeed}.

Recalling \eqref{011307}, we write \eqref{dbetanew}, \eqref{dalphanew} as:
\begin{equation}
	d\hat{\psi}(t,k)=-i\omega(k)\hat{\psi}(t,k)dt+iF\star \alpha_{0}(t)dt+i\sqrt{2\nu T}dR(t).
	\label{dpsinew}
\end{equation}
Solving \eqref{dpsinew}, we obtain
\begin{equation}
	\hat{\psi}(t,k)=e^{-i\omega(k)t}\hat{\psi}(k)+ie^{-i\omega(k)\cdot}\star F\star \alpha_{0}(t)+i\sqrt{2\nu T}\int_{0}^{t}e^{-i\omega(k)(t-s)}dR(s).
	\label{solvepsinew}
\end{equation}
Hence,
\begin{equation}
	\alpha_{0}(t)=\Psi(t)+C_\omega\star F\star \alpha_0(t)+\sqrt{2\nu T}C_\omega\star dR(t),
	\label{alpha_0Comeganew}
\end{equation}
where we have used the notation of \eqref{Comegadef}.
To solve \eqref{alpha_0Comeganew}, we define, formally, $$\tilde{\Theta}_{F}(Z):=\frac{1}{1-\cL(F)(Z)\cL(C_\omega)(Z)}\quad\text{for }Z\in\C_+.$$

We are going to see that $\tilde{\Theta}_{F}$ is well defined and is bounded, similarly to $\tilde{\Theta}$. First, we infer from \eqref{LaplaceFourier} that
\begin{align}
	\theta_F(k):=\lim_{Z\to0}\tilde{\Theta}_{F}(Z-i\omega(k))=\frac{1}{1-\hat{F}^{*}(\omega(k)/2\pi)\lim_{Z\to0}\cL(C_\omega)(Z-i\omega(k))}.
	\label{tildeThetaFlimit}
\end{align}
Note that, $\lim_{Z\to0}\cL(C_\omega)(Z-i\omega(k))=\frac{1}{\nu}\left( \frac{1}{\theta(k)}-1 \right)$ if $\theta(k)\ne0$. By \ref{inequalitycondition1}, the limit $\lim_{Z\to0}\tilde{\Theta}_{F}(Z-i\omega(k))$ exists. When $\theta(k)=0$, we have $|\cL(C_\omega)(Z-i\omega(k))|\to\infty$, and we also get $\theta_{F}(k)=0$. Since $\theta(k)$ is defined almost everywhere, so is $\theta_{F}(k)$.

The boundedness of $\tilde{\Theta}_{F}(Z-i\omega(k))$ comes from \ref{inequalitycondition1} and it is satisfied almost everywhere. This can be proved as follows.
By \ref{inequalitycondition1}, we get $\left|1-\frac{\hat{F}(\omega(k)/(2\pi))}{\nu}\left( \frac{1}{\theta(k)}-1 \right)\right|\ge C>0$, for a constant $C$. Indeed, the function $\Real(\hat{F})$ is bounded continuous and $\Real(\hat{F}(\omega(k)/(2\pi)))<0$ so there exist $C',C''>0$ such that $C'\ge|\hat{F}(\omega(k)/(2\pi))|\ge|\Real(\hat{F}(\omega(k)/(2\pi)))|\ge C''$. Thus
	\begin{align*}
	\left|1-\frac{\hat{F}(\omega(k)/(2\pi))}{\nu}\left( \frac{1}{\theta(k)}-1 \right)\right|&= \frac{|\hat{F}(\omega(k)/(2\pi))|}{\nu}\left|\frac{\nu}{\hat{F}(\omega(k)/(2\pi))}-\left( \frac{1}{\theta(k)}-1 \right)\right|\\
	&\ge\frac{C''}{\nu}\left|\Real\left(  \frac{\nu}{\hat{F}(\omega(k)/(2\pi))}\right)-\Real\left( \frac{1}{\theta(k)}-1 \right)\right|\\
	&\ge\frac{C''}{\nu}\frac{\nu(-\Real(\hat{F}))}{|\hat{F}|^2}\ge\left( \frac{C''}{C'} \right)^{2}>0.
				\end{align*}

				Similar to Lemma \ref{thetaprop}, we have another lemma for linear memory-feedback controls.
				\begin{lemma}
					\begin{equation*}
						\Real(\hat{F}^{*}(\omega(k)/2\pi)\theta_{F}(k))=\left( \Real(\hat{F}(\omega(k)/2\pi))-|\hat{F}(\omega(k)/2\pi)|^{2}\frac{\pi}{|\omega'(k)|} \right)|\theta_{F}(k)|^2.
					\end{equation*}
As a consequence,
\begin{align}
	r_a^{F}(k)+r_t^{F}(k)+r_r^{F}(k)=1
	\label{sumrates}
\end{align}
for $k\in\T$ almost everywhere.

					\label{thetaFprop}
				\end{lemma}
				The lemma is an immediate consequence of Lemma \ref{limLC}, we show the proof here.
				\begin{proof}
								Using \eqref{tildeThetaFlimit}, we get
								\begin{align*}
									\hat{F}^{*}(\omega(k)/2\pi)\theta_{F}(k)&= \frac{\hat{F}^{*}(\omega(k)/2\pi)}{1-\hat{F}^{*}(\omega(k)/2\pi)\lim_{\varepsilon\to0}\cL(C_\omega)(\varepsilon-i\omega(k))}\\
									&= |\hat{F}(\omega(k)/2\pi)|^{2}\frac{1}{\hat{F}(\omega(k)/2\pi)-|\hat{F}(\omega(k)/2\pi)|^{2}\lim_{\varepsilon\to0}\cL(C_\omega)(\varepsilon-i\omega(k))}.
							\end{align*}
						From Lemma \ref{limLC}, we have $\Real(\lim_{\varepsilon\to0}\cL(C_\omega)(\varepsilon-i\omega(k)))=\frac{\pi}{|\omega'(k)|}$. Thus, using $\Real z=|z|^{2}\Real \frac{1}{z}$, the conclusion of the lemma follows.
					\end{proof}

We introduce more notations. We write \eqref{solvepsinew} in the form
\begin{align}
	\hat{\psi}(t,k)=e^{-i\omega(k)t}\hat{\psi}(k)+i\Theta_\omega^F\star F\star\Psi(t)+i\sqrt{2\nu T}\Theta_\omega^F\star R(t),
	\label{psiexpandnew}
\end{align}
where
\begin{align*}
	\Theta_\omega^F(t,k):=\int_{0}^{t}e^{-i\omega(k)(t-s)}\Theta_{F}(ds).
	%\label{Theta_omegaFdef}
\end{align*}
We will also consider
\begin{align*}
	\Theta^\omega_{F}(t,k)&:= \int_{0}^{t}e^{i\omega(k)s}\Theta_{F}(ds)= e^{i\omega(k)t}\Theta_\omega^F(t,k),\\
	\cL( (\Theta^\omega_{F})^{(\varepsilon)})(Z,k)&= \frac{\tilde{\Theta}_{F}(\varepsilon Z-i\omega(k))}{Z}.
	%\label{Theta^omegaF}
\end{align*}

Using $[\T](\varrho)$ defined by \eqref{Tvarrho}, we define a truncation similar to \eqref{multiplybump}:
\begin{align}
	\hat{\psi}^{1}(0,k)&= \hat{\psi}(k)\chi_{\varrho}(k),\label{multiplybump1}
	\\
	d\hat{\psi}^{1}(t,k)&= \left\{ -i\omega(k)\hat{\psi}^{1}(t,k)+\frac{1}{2i}\int_{\T}\int_{0}^{t}F(s)[\hat{\psi}^{1}(t-s,k')-\hat{\psi}^{1*}(t-s,k')]dsdk' \right\}dt\nonumber\\
	&\hspace{20em}+i\sqrt{2\nu T}dR(t),\nonumber%\label{newpsi1}
	\\
	\hat{\psi}^{2}(0,k)&= \hat{\psi}(k)(1-\chi_{\varrho}(k)),\nonumber%\label{trimpart1}
	\\
	d\hat{\psi}^{2}(t,k)&= \left\{ -i\omega(k)\hat{\psi}^{1}(t,k)+\frac{1}{2i}\int_{\T}\int_{0}^{t}F(s)[\hat{\psi}^{2}(t-s,k')-\hat{\psi}^{2*}(t-s,k')]dsdk' \right\}dt.\nonumber%\label{eliminatedpsi1}
\end{align}
Similar to \eqref{multiplybump}, for sufficiently small $\varrho$, $\psi^1$ in \eqref{multiplybump1} and $\psi$ in \eqref{dpsinew} give the same kinetic limit. Therefore, in the scattering terms (Section \ref{ScatLM}), we only consider the situation where the initial state satisfies \eqref{Trunass}.

\subsection{Thermal term}\label{therF}
Due to the independence between the initial wave and the thermostat, we simplify the computation of the thermal term by setting the initial wave to be $0$.
The thermal term in linear memory-feedback control is
\begin{equation*}
	\nu T\int_{0}^{t}(\Theta_{\omega}^{F})^{(\varepsilon)}(s,k+\varepsilon\xi/2)(\Theta_{\omega}^{F})^{*(\varepsilon)}(s,k-\varepsilon\xi/2)ds.
\end{equation*}
In that case, we compute 
\begin{align*}
	\sW_{\varepsilon}(Z,\xi,k)= \frac{\nu T}{Z}\int_{0}^{\infty}e^{-(Z+iD_{\varepsilon}(k,\xi))s}\Theta^\omega_{F}\left( \varepsilon^{-1}s,k+\varepsilon\xi/2 \right)(\Theta^\omega_{F})^{*}\left( \varepsilon^{-1}s,k-\varepsilon\xi/2 \right)ds\\
	= \frac{\nu T}{Z}\frac{1}{2\pi i}\lim_{L\to\infty}\int_{c-iL}^{c+iL}\frac{\tilde{\Theta}_{F}(\varepsilon\zeta-i\omega(k+\varepsilon\xi/2))\tilde{\Theta}_{F}^{*}(\varepsilon(Z+iD_{\varepsilon}(k,\xi)-\zeta)-i\omega(k-\varepsilon\xi/2))}{\zeta(Z+iD_{\varepsilon}(k,\xi)-\zeta)}d\zeta.
	%\label{newthermal}
\end{align*}

Taking the limit $\varepsilon\to0$ and using \eqref{integralfraction}, we obtain
\begin{align*}
	\frac{\nu T|\theta_{F}(k)|^2}{Z(Z+i\omega'(k)\xi)}=\frac{\nu T|v_g(k)|r_a^{F}(k)}{\Real(\hat{F}(\omega(k)/(2\pi)))Z(Z+i\omega'(k)\xi)}.
	%\label{newthermallimit}
\end{align*}

\subsection{Scattering terms}\label{ScatLM}

For scattering terms, we consider $T=0$. The derivative of the Wigner transform is given by
\begin{align*}
	\partial_{t}\hat{W}_{\varepsilon}(t,\xi,k)&= -iD_{\varepsilon}(k,\xi)\hat{W}_{\varepsilon}(t,\xi,k)\\
	&+\E_{\varepsilon}\left[ i(F\star\Psi)^{(\varepsilon)}(t)\hat{\psi}^{(\varepsilon)*}(t,k-\varepsilon\xi/2)-i(F\star\Psi)^{(\varepsilon)}(t)\hat{\psi}^{(\varepsilon)}(t,k+\varepsilon\xi/2) \right].
	%\label{derivativeWnew}
\end{align*}

Using the Laplace transform we obtain
\begin{align}
	(Z+iD_{\varepsilon}(k,\xi))&\sW_{\varepsilon}(Z,\xi,k)=\hat{W}_{\varepsilon}(\xi,k) \label{ballisticnew}\\
	&+ S_{\varepsilon}^{I}(Z,k-\varepsilon,\xi/2,F)+S_{\varepsilon}^{I*}(Z,k+\varepsilon,\xi/2,F)\label{Iscatnew}\\
	&+ S_{\varepsilon}^{II}(Z,k-\varepsilon,\xi/2,F)+S_{\varepsilon}^{II*}(Z,k+\varepsilon,\xi/2,F)\label{IIscatnew},
\end{align}
where $S_{\varepsilon}^{I}$ and $S_{\varepsilon}^{II}$ will be defined in \eqref{S1new} and \eqref{S2new}.

\subsubsection{The ballistic term \eqref{ballisticnew}}
The ballistic term is unchanged with respect to \eqref{limitballistic}. The term only involves the initial condition of the wave, neither the control nor the thermostat affects this term.

\subsubsection{The first scattering term \eqref{Iscatnew}}
\label{Iscatterpoint}
The first scattering term is given by
\begin{align*}
	I_{\varepsilon}(t,k,F)&:= ie^{i\omega(k)t}\int_{0}^{t}\langle F\star\Psi(t-s)\hat{\psi}^{*}(k)\rangle_{\mu_{\varepsilon}}d\Theta_{F}(s)\\
	&= \frac{1}{2}\int_{0}^{t}\int_{\T}\int_{0}^{t-s}F(s')\langle\hat{\psi}^{*}(k)\hat{\psi}(h)\rangle_{\mu_{\varepsilon}} e^{i(\omega(k)-\omega(h))t+i\omega(h)(s+s')}ds'dhd\Theta(s),%\label{Inew}
\end{align*}
and
\begin{equation}
	S_{\varepsilon}^{I}(Z,k,F):= \cL(I_{\varepsilon}^{(\varepsilon)})(Z,k,F).\label{S1new}
\end{equation}
We have 
\begin{align*}
	S_{\varepsilon}^{I}(Z,k,F)&= \frac{\varepsilon}{2}\int_{0}^{\infty}e^{-\varepsilon Zt}\int_{0}^{t}\int_{\T}\int_{0}^{t-s}F(s')\langle\hat{\psi}^{*}(k)\hat{\psi}(h)\rangle_{\mu_{\varepsilon}} e^{i(\omega(k)-\omega(h))t+i\omega(h)(s+s')}ds'dhd\Theta(s)dt\\
	&= \frac{\varepsilon}{2}\int_{\T}\langle\hat{\psi}^{*}(k)\hat{\psi}(h)\rangle_{\mu_{\varepsilon}}\int_{0}^{\infty}F(s')\int_{0}^{\infty}\int_{s+s'}^{\infty}e^{-\varepsilon Zt} e^{i(\omega(k)-\omega(h))t+i\omega(h)(s+s')}dtd\Theta(s)ds'dh\\
	&= \frac{\varepsilon}{2}\int_{\T}\frac{\langle\hat{\psi}^{*}(k)\hat{\psi}(h)\rangle_{\mu_{\varepsilon}}}{\lambda\varepsilon+i(\omega(h)-\omega(k))}\int_{0}^{\infty}F(s')\int_{0}^{\infty} e^{-\varepsilon Z(s+s')}e^{i\omega(k)(s+s')}d\Theta(s)ds'dh\\
	&= \frac{\varepsilon}{2}\int_{\T}\frac{\langle\hat{\psi}^{*}(k)\hat{\psi}(h)\rangle_{\mu_{\varepsilon}}}{\lambda\varepsilon+i(\omega(h)-\omega(k))}\cL(F)(\varepsilon Z-i\omega(k))\tilde{\Theta}_{F}(\varepsilon Z-i\omega(k))dh.
\end{align*}

Changing the variables from $k$ to $k'-\varepsilon\xi'/2$ and $h$ to $k'+\varepsilon\xi'/2$, we call the new domain $T_{\varepsilon}\subset\T_{2/\varepsilon}\times\T$. We get
\begin{align*}
	&\int_{\R\times\T}\frac{\hat{O}^{*}(\xi,k+\varepsilon\xi/2)S^{1}_{\varepsilon}(Z,k,F)}{Z+iD_{\varepsilon}^{+}(k,\xi)}d\xi dk\\
	&= \int_{\R}\int_{T_{\varepsilon}}\frac{\hat{W}_{\varepsilon}(\xi',k')\cL(F)(\varepsilon Z-i\omega(k'-\varepsilon\xi/2))\tilde{\Theta}_{F}(\varepsilon Z-i\omega(k'-\varepsilon\xi'/2))\hat{O}^{*}(\xi,k'+\varepsilon\xi/2-\varepsilon\xi'/2)}{(\varepsilon Z+i(\omega(k'+\varepsilon\xi'/2)-\omega(k'-\varepsilon\xi'/2)))(Z+iD_{\varepsilon}^{+}(k'-\varepsilon\xi'/2,\xi))}d\xi'd\xi dk\\
	&\overset{\varepsilon\to0}{\longrightarrow} \int_{\R^2\times\T}\frac{\hat{W}_{0}(\xi',k')\hat{F}^{*}(\omega(k')/(2\pi))\theta_{F}(k)\hat{O}^{*}(\xi,k')}{(Z+i\omega'(k')\xi)(Z+i\omega'(k')\xi')}d\xi'd\xi dk.
	%\label{S1newchangevari}
\end{align*}
Similarly, we have
\begin{align*}
	\int_{\R\times\T}\frac{\hat{O}^{*}(\xi,k-\varepsilon\xi/2)S^{I*}_{\varepsilon}(Z,k,F)}{Z+iD_{\varepsilon}^{-}(k,\xi)}d\xi dk\overset{\varepsilon\to0}{\longrightarrow} \int_{\R^2\times\T}\frac{\hat{W}^{*}_{0}(\xi',k')\hat{F}(\omega(k')/(2\pi))\theta^{*}_{F}(k)\hat{O}^{*}(\xi,k')}{(Z+i\omega'(k')\xi)(Z-i\omega'(k')\xi')}d\xi'd\xi dk.
	%\label{S1starnew}
\end{align*}
Taking the sum and using the fact that $\hat{W}^{*}_{0}(\eta,k)=\hat{W}_{0}(-\eta,k)$, we get the limit of the first scattering term  
\begin{align*}
	\int_{\R^{2}\times\T} \Real\left( \hat{F}^{*}(\omega(k)/(2\pi))\theta_{F}(k) \right)\frac{\hat{W}_{0}(\xi',k)O^{*}(\xi,k)}{(Z+i\omega'(k)\xi')(Z+i\omega'(k)\xi)}d\xi'd\xi dk.
	%\label{firstscatterlim}
\end{align*}

\subsubsection{The second scattering term \eqref{IIscatnew}}
\label{IIscatterpoint}

The second scattering term is given by
\begin{align*}
	II_{\varepsilon}(t,k,F):=\int_{0}^{t}\int_{0}^{t}\Theta_\omega^{F}(t-s,k)\langle F\star\Psi(s)F\star\Psi(t-s')\rangle_{\mu_{\varepsilon}}dsd\Theta(s'),
	%\label{IInew}
\end{align*}
and
\begin{equation}
	S_{\varepsilon}^{II}(Z,k,F):= \cL(II_{\varepsilon}^{(\varepsilon)})(Z,k,F).\label{S2new}
\end{equation}
We compute
\begin{align}
	&\Real S^{II}_{\varepsilon}(Z,k,F)=\Real \cL(II_\varepsilon^{(\varepsilon)})(Z,k,F)\nonumber\\
	&= \varepsilon^2Z\int_{0}^{t}e^{-\varepsilon Zt}\left\langle\left( \int_{0}^{t}\cos(\omega(k)s)(\Theta_{F}\star F\star \Psi)(s)ds \right) \left(\int_{0}^{t}\cos(\omega(k)s)(\Theta_{F}\star F\star \Psi)(s)ds \right)\right\rangle_{\mu_{\varepsilon}}dt\label{sumreal2cosnew}\\
	&+\varepsilon^2Z\int_{0}^{t}e^{-\varepsilon Zt}\left\langle\left( \int_{0}^{t}\sin(\omega(k)s)(\Theta_{F}\star F\star \Psi)(s)ds \right) \left(\int_{0}^{t}\sin(\omega(k)s)(\Theta_{F}\star F\star \Psi)(s)ds \right)\right\rangle_{\mu_{\varepsilon}}dt.\label{sumreal2sinnew}
\end{align}
Let us analyze \eqref{sumreal2cosnew}. \eqref{sumreal2sinnew} can be treated similarly.

We write  
\begin{align}
	&\frac{\varepsilon^{2}Z}{2}\int_{0}^{\infty}e^{-\varepsilon Zt}\int_{0}^{t}\int_{0}^{t}\cos(\omega(k)s)\cos(\omega(k)s')\int_{0}^{s}\int_{0}^{s'} \nonumber\\
	&\qquad \times \int_{\T^2}e^{-i\omega(h)(s-u)}e^{i\omega(h')(s'-u')}\left\langle \hat{\psi}(h)\hat{\psi}^{*}(h')\right\rangle_{\mu_{\varepsilon}}dh'dhd\Theta_{F}\star F(u') d\Theta_{F}\star F(u)ds'dsdt%\label{sumreal2cos1}
	\nonumber\\
	&= \frac{\varepsilon Z}{4\pi}\int_{\R} \int_{\T^2}\varepsilon\left\langle \hat{\psi}(h)\hat{\psi}^{*}(h')\right\rangle_{\mu_{\varepsilon}}\int_{0}^{\infty}\int_{0}^{\infty}e^{ia(t-t')}e^{-\varepsilon Z(t+t')/2}\int_{0}^{t}\int_{0}^{t}\cos(\omega(k)s)\cos(\omega(k)s')\nonumber\\
	&\qquad \times \int_{0}^{s}\int_{0}^{s'} e^{-i\omega(h)(s-u)}e^{i\omega(h')(s'-u')}d\Theta_{F}\star F(u') d\Theta_{F}\star F(u)ds'dsdt'dtdh'dhda\nonumber\\
	&= \frac{\varepsilon Z}{4\pi}\int_{\R} \int_{\T^2} \varepsilon\left\langle \hat{\psi}(h)\hat{\psi}^{*}(h')\right\rangle_{\mu_{\varepsilon}}A(h,k,\lambda,F)A^{*}(h',k,\lambda,F)dh'dhda,\label{sumreal2cosA}
\end{align}
in which
\begin{align}
	A(h,k,Z,F)&:= \int_{0}^{\infty}\cos(\omega(k)s)\int_{0}^{s}e^{-i\omega(h)(s-u)}\int_{s}^{\infty}e^{(-\varepsilon Z/2+ia)t}dtd\Theta_{F}\star F(u)ds\nonumber\\
	&= \frac{1}{2(\varepsilon Z/2-ia)} \frac{\cL(F)(\varepsilon Z/2-ia-i\omega(k))\tilde{\Theta}_{F}(\varepsilon Z/2-ia-i\omega(k))}{\varepsilon Z/2-ia-i(\omega(k)-\omega(h))}\label{Anewdefminus}\\ 
	&\qquad+\frac{1}{2(\varepsilon Z/2-ia)}\frac{\cL(F)(\varepsilon Z/2-ia+i\omega(k))\tilde{\Theta}_{F}(\varepsilon Z/2-ia+i\omega(k))}{\varepsilon Z/2-ia+i(\omega(k)+\omega(h))},
	\label{Anewdefplus}
\end{align}
and $A^{*}$ is defined in a similar manner, that is
\begin{align}
	A^{*}(h',k,Z,F)
	&:= \frac{1}{2(\varepsilon Z/2+ia)} \frac{\cL(F)(\varepsilon Z/2+ia+i\omega(k))\tilde{\Theta}_{F}(\varepsilon Z/2+ia+i\omega(k))}{\varepsilon Z/2+ia+i(\omega(k)-\omega(h'))}\label{Astardefminus}\\ 
	&\qquad+\frac{1}{2(\varepsilon Z/2+ia)}\frac{\cL(F)(\varepsilon Z/2+ia-i\omega(k))\tilde{\Theta}_{F}(\varepsilon Z/2+ia-i\omega(k))}{\varepsilon Z/2+ia-i(\omega(k)+\omega(h'))}.\label{Astardefplus}
\end{align}

We will see that among the four terms in the expansion of \eqref{sumreal2cosA}, these are the products between \eqref{Anewdefminus} or \eqref{Anewdefplus} and \eqref{Astardefminus} or \eqref{Astardefplus}, only the term generated by \eqref{Anewdefminus} and \eqref{Astardefminus} contributes to the high frequency limit. To derive the limit, we need some lemmas.

\begin{lemma} The following limit holds true
	\begin{align*}
	\lim_{\varepsilon\to0}&\frac{1}{\varepsilon}\int_{\R}\frac{1}{(Z/2)^2+a^2}\int_{\T^3}\left|\left\langle\hat{\psi}(h)\hat{\psi}^{*}(h')\right\rangle_{\mu_{\varepsilon}}\right|\\
	&\times\left|\frac{\cL(F)(\varepsilon Z/2-i\varepsilon a-i\omega(k))\tilde{\Theta}_{F}(\varepsilon Z/2-i\varepsilon a-i\omega(k))}{Z/2-ia-i\varepsilon^{-1}(\omega(k)-\omega(h))}\right|\\
	&\times\left|\frac{\cL(F)(\varepsilon Z/2+i\varepsilon a-i\omega(k))\tilde{\Theta}_{F}(\varepsilon Z/2+i\varepsilon a-i\omega(k))}{Z/2+ia-i\varepsilon^{-1}(\omega(k)+\omega(h'))}\right|dh'dhdkda=0
		%\label{limitplusdenom}
	\end{align*}

	\label{plusdenomlimitzero}
\end{lemma}

The proof of the lemma will be provided later in Section \ref{prodenomlimitzero}. To use Lemma \ref{plusdenomlimitzero}, we change $a$ into $\varepsilon a$ in \eqref{sumreal2cosA}. Therefore,
\begin{align*}
	\frac{Z}{16\pi\varepsilon}&\int_{\R}\frac{1}{(Z/2)^2+a^2}\int_{\T^2}\varepsilon\left\langle\hat{\psi}(h)\hat{\psi}^{*}(h')\right\rangle_{\mu_{\varepsilon}}\\
	&\times\frac{\cL(F)(\varepsilon Z/2-i\varepsilon a-i\omega(k))\tilde{\Theta}_{F}(\varepsilon Z/2-i\varepsilon a-i\omega(k))}{Z/2-i a-i\varepsilon^{-1}(\omega(k)-\omega(h))}\\
	&\times\frac{\cL(F)(\varepsilon Z/2+i\varepsilon a+i\omega(k))\tilde{\Theta}_{F}(\varepsilon Z/2+i\varepsilon a+i\omega(k))}{Z/2+i a+i\varepsilon^{-1}(\omega(k)-\omega(h'))}dh'dhda
	%\label{2denomminus}
\end{align*}
is the only term contributing in the limit. By changing the variables from $h$ into $k'+\varepsilon\xi'/2$ and $h'$ into $k'-\varepsilon\xi'/2$, we need to evaluate
\begin{align*}
	&\frac{Z}{16\pi}\int_{\R\times\T_{2/\varepsilon}}\frac{1}{(Z/2)^2+a^2}\int_{\T\times T^{2}_{\varepsilon}}\frac{\hat{W}_{\varepsilon}(\xi',k')}{Z/2-ia-i\varepsilon^{-1}(\omega(k)-\omega(k'+\varepsilon\xi'/2))}\\
	&\times\frac{\left|\cL(F)(\varepsilon Z/2-i\varepsilon a-i\omega(k))\tilde{\Theta}_{F}(\varepsilon Z/2-i\varepsilon a-i\omega(k))\right|^2\hat{O}^{*}(\xi,k+\varepsilon\xi/2)}{(Z/2+ia+i\varepsilon^{-1}(\omega(k)-\omega(k'-\varepsilon\xi'/2)))(Z+iD_{\varepsilon}^{+}(k,\xi))}dk'd\xi'dkd\xi da.
	%\label{changevar2denomminus}
\end{align*}
In the integral, we have the set $T^{2}_{\varepsilon}=\left\{(\xi',k'):|\xi'|\le\frac{\rho}{2^{100}\varepsilon},|k'|\le\frac{1-\varepsilon|\xi'|}{2}\right\}\subset \T_{2/\varepsilon}\times\T$, where $\rho$ is in \eqref{Trunass}. By Assumption \ref{WbracketR}, if $|\xi'|>\frac{\rho}{2^{100}\varepsilon}$, then $|\hat{W}_{\varepsilon}|\lesssim \varepsilon^{3/2}$; it makes the whole term tend  to $0$.

We define
\begin{align*}
	\cI_{\pm}(Z,\varepsilon,F)&= \frac{Z}{32\pi}\int_{\R\times\T_{2/\varepsilon}}\frac{1}{(Z/2)^2+a^2}\int_{T^{3}_{\varepsilon,\pm}}\frac{\hat{W}_{\varepsilon}(\xi',\pm k+\varepsilon\xi''/2)}{Z/2-ia-i\varepsilon^{-1}(\omega(k)-\omega(k\pm\varepsilon(\xi'+\xi'')/2))}\\
	&\times\frac{\left|\cL(F)(\varepsilon Z/2-i\varepsilon a-i\omega(k))\tilde{\Theta}_{F}(\varepsilon Z/2-i\varepsilon a-i\omega(k))\right|^2\hat{O}^{*}(\xi,k+\varepsilon\xi/2)}{(Z/2+ia+i\varepsilon^{-1}(\omega(k)-\omega(k\mp\varepsilon(\xi'-\xi'')/2)))(Z+iD_{\varepsilon}^{+}(k,\xi))}d\xi''d\xi'dkd\xi da.
	%\label{cI}
\end{align*}
In this integral,
\[T^{3}_{\varepsilon,\pm}=\left\{(k,\xi',\xi''):k\in\T,|\xi'|\le\frac{\varrho}{2^{50}\varepsilon},|k\pm\varepsilon\xi''/2|\le\frac{1-\varepsilon|\xi'|}{2},\sign k=\sign(k\pm\varepsilon(\xi'+\xi'')/2)\right\},\]
where we have used the change of variables $k'\to\pm k+\varepsilon\xi''/2$.

We also define the following integrals
	\begin{align*}
		\cI_{\pm}^{(1)}(Z,\varepsilon,F)&:= \frac{Z}{32\pi}\int_{\R\times\T_{2/\varepsilon}}\frac{1}{(Z/2)^2+a^2}\int_{T^{3}_{\varepsilon,\pm}}\frac{\hat{W}_{\varepsilon}(\xi',\pm k+\varepsilon\xi''/2)}{Z/2-ia\pm i\omega'(k)(\xi'+\xi'')/2}\\
		&\times\frac{\left|\cL(F)(\varepsilon Z/2-i\varepsilon a-i\omega(k))\tilde{\Theta}_{F}(Z\varepsilon/2-i\varepsilon a-i\omega(k))\right|^2\hat{O}^{*}(\xi,k+\varepsilon\xi/2)}{(Z/2+i a+i\varepsilon^{-1}(\omega(k)-\omega(k\mp\varepsilon(\xi'-\xi'')/2)))(Z+iD_{\varepsilon}^{+}(k,\xi))}d\xi''d\xi'dkd\xi da;
		%\label{cI1}
	\end{align*}

	\begin{align*}
		\cI_{\pm}^{(2)}(Z,\varepsilon,F)&:= \frac{Z}{32\pi}\int_{\R\times\T_{2/\varepsilon}}\frac{1}{(Z/2)^2+a^2}\int_{T^{3}_{\varepsilon,\pm}}\frac{\hat{W}_{\varepsilon}(\xi',\pm k+\varepsilon\xi''/2)}{Z/2-ia\pm i\omega'(k)(\xi'+\xi'')/2}\\
		&\times\frac{\left|\cL(F)(\varepsilon Z/2-i\varepsilon a-i\omega(k))\tilde{\Theta}_{F}(\varepsilon Z/2-i\varepsilon a-i\omega(k))\right|^2\hat{O}^{*}(\xi,k+\varepsilon\xi/2)}{(Z/2+ia\pm i\omega'(k)(\xi'-\xi'')/2)(Z+iD_{\varepsilon}^{+}(k,\xi))}d\xi''d\xi'dkd\xi da;
		%\label{cI2}
	\end{align*}

	\begin{align*}
		\cI_{\pm}^{(3)}(Z,\varepsilon,F)&:= \frac{Z}{32\pi}\int_{\R\times\T_{2/\varepsilon}}\frac{1}{(Z/2)^2+a^2}\int_{T^{3}_{\varepsilon,\pm}}\frac{\hat{W}_{\varepsilon}(\xi',\pm k+\varepsilon\xi''/2)}{Z/2-ia\pm i\omega'(k)(\xi'+\xi'')/2}\\
		&\times\frac{\left|\hat{F}(\omega(k)/(2\pi))\theta_{F}(k)\right|^2\hat{O}^{*}(\xi,k+\varepsilon\xi/2)}{(Z/2+ia\pm i\omega'(k)(\xi'-\xi'')/2)(Z+iD_{\varepsilon}^{+}(k,\xi))}d\xi''d\xi'dkd\xi da;
		%\label{cI3}
	\end{align*}

	\begin{align*}
		\cI_{\pm}^{(4)}(Z,F)&:= \frac{Z}{32\pi}\int_{\R^2}\frac{1}{(Z/2)^2+a^2}\int_{\T\times\R^2}\frac{\hat{W}_{0}(\xi',\pm k)}{Z/2-ia\pm i\omega'(k)(\xi'+\xi'')/2}\\
		&\times\frac{\left|\hat{F}(\omega(k)/(2\pi))\theta_{F}(k)\right|^2\hat{O}^{*}(\xi,k)}{(Z/2+ia\pm i\omega'(k)(\xi'-\xi'')/2)(Z+i\omega'(k)\xi)}d\xi''d\xi'dkd\xi da.
		%\label{cI4}
	\end{align*}
\begin{lemma}
	We have the following approximations
	\begin{align}
		\lim_{\varepsilon\to0}\left|\cI_{\pm}(Z,\varepsilon,F)-\cI_{\pm}^{(1)}(Z,\varepsilon,F)\right|=0;
		\label{approxI1}
	\end{align}

	\begin{align}
		\lim_{\varepsilon\to0}\left|\cI_{\pm}^{(1)}(Z,\varepsilon,F)-\cI_{\pm}^{(2)}(Z,\varepsilon,F)\right|=0;
		\label{approxI2}
	\end{align}

	\begin{align}
		\lim_{\varepsilon\to0}\left|\cI_{\pm}^{(2)}(Z,\varepsilon,F)-\cI_{\pm}^{(3)}(Z,\varepsilon,F)\right|=0;
		\label{approxI3}
	\end{align}

	\begin{align}
		\lim_{\varepsilon\to0}\left|\cI_{\pm}^{(3)}(Z,\varepsilon,F)-\cI_{\pm}^{(4)}(Z,F)\right|=0.
		\label{approxI4}
	\end{align}
	As a consequence,
	\[\lim_{\varepsilon\to0}\cI_{\pm}(Z,\varepsilon,F)=\cI_{\pm}^{(4)}(Z,F).\]
	\label{lemI1-4}
\end{lemma}

\begin{lemma}
	\begin{align*}
		\lim_{\varepsilon\to0}\int_{\R\times\T}\Imagine S^{II}_{\varepsilon}(Z,k,F)\left( \frac{\hat{O}^{*}(\xi,k+\varepsilon\xi/2)}{Z+iD_{\varepsilon}^{+}(k,\xi)}-\frac{\hat{O}^{*}(\xi,k-\varepsilon\xi/2)}{Z+iD_{\varepsilon}^{-}(k,\xi)} \right)d\xi dk=0.
		%\label{limim}
	\end{align*}
	\label{limim0}
\end{lemma}
The proofs of Lemma \eqref{lemI1-4} and Lemma \eqref{limim0} will be given in Section \ref{proapp}.

We apply to $\cI^{(4)}_{\pm}$ the following identity
\begin{align*}
	\int_{\R}\frac{dz}{(z-z_{+})(z-z_{-})}=\frac{2\pi i}{z_{+}-z_{-}} \text{ when }\Imagine(z_{+})>0>\Imagine(z_{-}).
	%\label{identityz}
\end{align*}
The variable $z$ is $a\mp\omega'(k)\xi''$ in our case, $d\xi''$   is integrated out first. Then, we also integrate with respect to $a$:
\begin{align*}
	\cI^{(4)}_{\pm}(\lambda,F)&= \frac{Z}{16}\int_{\R^2}\frac{1}{(Z/2)^2+a^2}\int_{\T\times\R}\frac{|\hat{F}(\omega(k)/(2\pi))|^{2}|\theta_{F}(k)|^{2}\hat{W}_{0}(\xi',\pm k)\hat{O}(\xi,k)}{|\omega'(k)|(Z\pm i\omega(k) \xi')(Z+i\omega'(k)\xi)}d\xi'dkd\xi da\\
	&= \int_{\T\times\R^{2}}\frac{|\hat{F}(\omega(k)/(2\pi))|^{2}|\theta_{F}(k)|^{2}\hat{W}_{0}(\xi',\pm k)\hat{O}(\xi,k)}{16|v_g(k)|(Z\pm i\omega(k)\eta')(Z+i\omega'(k)\xi)}d\xi'dk d\xi.
	%\label{cI4result}
\end{align*}
By repeating these steps, the  term \eqref{sumreal2sinnew} will produce the same result as \eqref{sumreal2cosnew}. That concludes the computation for $S^{II}_{\varepsilon}(Z,k,F)$. We also have the same result for $\frac{\hat{O}^{*}(k-\varepsilon\xi/2)}{Z+iD_{\varepsilon}^{-}(k,\xi)}$, so the result get doubled. The second scattering term limit is
\begin{align*}
	&\int_{\T\times\R^{2}}\frac{|\hat{F}(\omega(k)/(2\pi))|^{2}|\theta_{F}(k)|^{2}\hat{W}_{0}(\xi',k)\hat{O}(\xi,k)}{4|v_g(k)|(Z+i\omega(k)\xi')(Z+i\omega'(k)\xi)}d\xi'd\xi dk%\label{secondscatplus}
	\\
	&+\int_{\T\times\R^{2}}\frac{|\hat{F}(\omega(k)/(2\pi))|^{2}|\theta_{F}(k)|^{2}\hat{W}_{0}(\xi',-k)\hat{O}(\xi,k)}{4|v_g(k)|(Z-i\omega(k)\xi')(Z+i\omega'(k)\xi)}d\xi' d\xi dk.%\label{secondscatminus}
\end{align*}

\subsection{Proof of Theorem \ref{mainnew}}
The technical detail for this part is similar to that of the impulsive control. We see that \eqref{thermal4limitshowing}, \eqref{ballisticlimitshowing}, \eqref{transmittinglimitshowing}, and \eqref{reflectinglimitshowing} do not depend on the rates. Thus, we do not repeat the same arguments to obtain the conclusion of Theorem \ref{mainnew} even though the rates are different.

\subsection{Proof of Corollary \ref{controlcoro}}
%\label{Controlresult}
In this section, we show how to design controls that asymptotically steer the system to some desired rates $(r_a,r_t,r_r)$ in the kinetic limit. First, we state some basic properties that the rates must satisfy and the implications in finding a control.

\begin{enumerate}
	\item By \eqref{absorbratenew1}, \eqref{transmissionnew1}, and \eqref{reflectionnew1}, the target rates must be defined when $\omega'(k)\ne0$. Hence the domain of definition is $(-1/2,0)\cup(0,1/2)$. Since  $\omega(\cdot)$ is an even function,     the three rates are even functions. We reduce the domain of definition to the positive branch $(0,1/2)$. This explains Assumption \ref{evenT}.
	\item By \eqref{sumrates}, the sum of the three rates is $1$, that is why we have Assumption \ref{sum=1} for the target rates. Hence, we only need to know  two among the three rates $(r_a,r_t,r_r)$. We focus on $r_t$ and $r_r$.
	\item All the rates are positive functions, i.e. $r_a^F,r_t^F,r_r^F>0$. For the target rates, we assume a stronger positivity, which is Assumption \ref{posrates}.
	\item Recall the Lemma \ref{limLC}: for $k\in(0,1/2)$, we have
\begin{align*}
	\lim_{Z\to0}\cL(C_\omega)(Z-i\omega(k))&= \frac{\pi}{|\omega'(k)|}+i\left( \int_{0}^{1/2}\frac{dh}{\omega(k)+\omega(h)}+\log\frac{\omega(k)-\omega_{min}}{\omega_{max}-\omega(k)}+\frac{C_{\omega,0}(\omega(k))}{|\omega'(k)|} \right),
	%\label{limComega}
\end{align*}
where $C_{\omega,0}$ is a continuous, bounded and real-valued function. Because of the continuity of $\hat{F}$, we deduce that the rates are continuous on $(0,1/2)$. This means the target rates need to be continuous on $(0,1/2)$, or target rates need to meet Assumption \ref{conrates}.

	\item By the triangular inequality, we have
\begin{align}
	\sqrt{r_t^F(k)}+\sqrt{r_r^F(k)}&= \left|1+\frac{\hat{F}^{*}(\omega(k)/(2\pi))\theta_{F}(k)}{2|v_g(k)|}\right|+\left|\frac{\hat{F}^{*}(\omega(k)/(2\pi))\theta_{F}(k)}{2|v_g(k)|}\right|\ge1.
	\label{sumsqrtrate}
\end{align}
This corresponds to Assumption \ref{sumsquareroot} on the target rates.

\end{enumerate}
To sum up, we consider two continuous functions $r_t,r_r$ on $(0,1/2)$ such that $r_t(k)+r_r(k)\le1-c_1$ and $\sqrt{r_t(k)}+\sqrt{r_r(k)}\ge1$. These are Assumptions \ref{evenT} - \ref{sumsquareroot}.

We define
\begin{align*}
	RE(k)&:= |v_g(k)|\left( r_t(k)-r_r(k)-1 \right)\qquad\forall k\in(0,1/2).%\label{RealFthetaFfind}
\end{align*}
It follows from \eqref{transmissionnew1} and \eqref{reflectionnew1} that $RE(k)=\Real(\hat{F}(\omega(k)/(2\pi))\theta_{F}(k))$.

We also define
\begin{align*}
	IM(k)&:= \sqrt{4|v_g(k)|^{2}r_r(k)-RE(k)^{2}}\qquad\forall k\in(0,1/2),
	%\label{ImagineFthetaFfind}
\end{align*}
which is well-defined because, using \eqref{sumsqrtrate},
\begin{align*}
	4|v_g(k)|^{2}r_r(k)-RE(k)^{2}&= |v_g(k)|^{2}(\sqrt{r_t(k)}+\sqrt{r_r(k)}-1)(\sqrt{r_t(k)}-\sqrt{r_r(k)}+1)\\
	&\qquad\times(-\sqrt{r_t(k)}+\sqrt{r_r(k)}+1)(\sqrt{r_t(k)}+\sqrt{r_r(k)}+1)\ge 0.
	%\label{welldefineImagineFtheta}
\end{align*}

Finally, we define a complex-valued function $FT(k)$, which corresponds to $\hat{F}(\omega(k)/(2\pi))\theta_{F}(k)$ by
\begin{align*}
	FT(k)&:= RE(k)+iIM(k)\qquad\forall k\in(0,1/2).
	%\label{FthetaFfind}
\end{align*}

By \eqref{tildeThetaFlimit}, the function corresponding to $\theta_{F}(k)$ is defined by
\begin{align*}
	TH(k)&:= 1+FT(k)\lim_{Z\to0}\cL(C_\omega)(Z-i\omega(k)).
	%\label{thetaFfind}
\end{align*}
We claim that $|TH(k)|\ge c_1/4$ for any $k\in(0,1/2)$. To prove this claim, we consider two cases: either $\Imagine(\lim_{\lambda\to0}\cL(C_\omega)(\lambda-i\omega(k)))\ge0$ or $\Imagine(\lim_{\lambda\to0}\cL(C_\omega)(\lambda-i\omega(k)))<0$. In the first case,
\begin{align*}
	\Real(TH(k))&= 1+RE(k)\Real(\lim_{Z\to0}\cL(C_\omega)(Z-i\omega(k)))+IM(k)\Imagine(\lim_{Z\to0}\cL(C_\omega)(Z-i\omega(k)))\\
	&\ge1-\frac{1+r_r(k)-r_t(k)}{2}=\frac{1+r_t(k)-r_r(k)}{2}\ge r_t(k)+c_{1}/2.
	%\label{RealthetaFbound}
\end{align*}
In the second case, we have either
\begin{align*}
	|IM(k)\Imagine(\lim_{Z\to0}\cL(C_\omega)(Z-i\omega(k)))|&\le\frac{1+r_t(k)-r_r(k)}{4}\\
	\Rightarrow\Real(TH(k))&\ge\frac{1+r_t(k)-r_r(k)}{4}\ge\frac{r_t(k)}{2}+\frac{c_{1}}{4},%\label{smallcase1}
\end{align*}
or
\begin{align*}
	|IM(k)\Imagine(\lim_{Z\to0}\cL(C_\omega)(Z-i\omega(k)))|&>\frac{1+r_t(k)-r_r(k)}{4},%\label{smallcase2}
\end{align*}
and then
\begin{align*}
	\Imagine(TH(k))&= RE(k)\Imagine(\lim_{Z\to0}\cL(C_\omega)(Z-i\omega(k)))+IM(k)\Real(\lim_{Z\to0}\cL(C_\omega)(Z-i\omega(k)))\\
	&> \sqrt{(1+r_r(k)-r_t(k))(1+r_t(k)-r_r(k))}/\sqrt{2}\ge\frac{\sqrt{c_{1}+2r_t(k)}\sqrt{c_{1}+2r_r(k)}}{\sqrt{2}}.
	%\label{ImaginethetaFbound}
\end{align*}
The claim is proved. In particular, $TH(k)\ne0$ for every $k\in(0,1/2)$, so we can define $\bar{F}(u):=FT(\varphi_{+}(u))/TH(\varphi_{+}(u)),u\in(\omm,\omM)$, which corresponds to $\hat{F}^{*}(\omega(k)/2\pi)$.

The next step consists of computing the inverse  Fourier   of $\bar{F}$. We recall that the control function $F$ is real-valued and is only defined on $[0,+\infty)$. It is more convenient to use the half-line Fourier transform defined by
\begin{align*}
	f^{c}(s)&= \int_{0}^{\infty}f(t)\cos(st)dt,%\label{halflinecos}
	\\
	f^{s}(s)&= \int_{0}^{\infty}f(t)\sin(st)dt%\label{halflinesin}
	,
\end{align*}
for all $f\in L^1([0,+\infty))\cap L^2([0,+\infty))$. We obtain $f^c,f^s\in L^2([0,+\infty))$.
It is also possible to define the inverses in the case the transformed functions are also in $L^{1}$ by
\begin{align*}
	f(x)&= \frac{2}{\pi}\int_{0}^{\infty}\int_{0}^{\infty}f(t)\cos(st)dt\cos(sx)ds,%\label{inversehalflinecos}
	\\
	f(x)&= \frac{2}{\pi}\int_{0}^{\infty}\int_{0}^{\infty}f(t)\sin(st)dt\sin(sx)ds.%\label{inversehalflinesin}
\end{align*}
By \eqref{LaplaceFourier}, if there is a control $F\in L^1$, we have
\begin{align*}
	\lim_{Z\to0}\cL(F)(Z+i\omega)&= F^{c}(\omega)-iF^{s}(\omega).
	%\label{LaplacehalflineFourier}
\end{align*}
We see that $\Real(\bar{F})$ and $\Imagine(\bar{F})$ correspond to $F^{c}$ and $F^{s}$, respectively. Thus, the existence of $F$ depends on whether we have
\begin{align*}
	\int_{0}^{\infty}\Real(\bar{F}(u))\cos(ut)du&= \int_{0}^{\infty}\Imagine(\bar{F}(u))\sin(ut)du,
	%\label{cosequalsin}
\end{align*}
for all $t>0$. This happens because if $F$ exists, then both sides equal $\frac{\pi}{2}F(t)$. Therefore, on the target rates, we expect them to satisfy Assumptions \ref{evenT} - \ref{intcossin}.

For a triple of rates that satisfy Assumptions \ref{evenT} - \ref{intcossin}, let us design  an explicit control function $F$ by 
\[F(t)=\frac{2}{\pi}\int_{0}^{\infty}\Real(\bar{F}(u))\cos(u t)du.\] It is easy to see that $\bar{F}$ is continuous. Therefore,
\begin{align*}
	\frac{d^{n}}{dt^{n}}F(t)=\frac{2}{\pi}\int_{0}^{\infty}\Real(\bar{F}(u))\frac{d^{n}}{dt^{n}}\cos(ut)du
	%\label{dF}
\end{align*}
are continuous for all $n$. This means that $F$ must belong to $C^{\infty}(0,+\infty)$. We also note that $|\bar{F}|\le\frac{8|v_g(k)|r_r(k)}{c_{1}}$ and $\bar{F}$ is supported on $[\omm,\omM]$, in particular $\bar{F}\in L^2(0,+\infty)$, and thus also $F\in L^{2}(0,+\infty)$.

We define a family of function $F_N\in C^{\infty}_c\subset C^\infty\cap L^m,m\ge1$. For each $N\in\N$, $F_N$ is the product of $F$ with a smooth function with values in $[0,1]$ such that it is equal to $1$ on $[0,N]$ and to $0$ on $[N+1,+\infty)$. It is clear that $F_N$ converges to $F$ in $L^2$. Thus, $F_{N}^{c},F_{N}^{s}$ converge to $F^{c}=\Real(\bar{F}),F^{s}=\Imagine(\bar{F})$ in $L^2$, respectively. Taking a subsequence if necessary, we obtain a family $F_N$ such that $\hat{F}^{*}_{N}(\cdot/2\pi)$ converges almost everywhere to $\bar{F}$. We have 
\begin{equation*}
	FT(k)=\frac{\bar{F}(\omega(k))}{1-\bar{F}(\omega(k))\lim_{Z\to0}\cL(C_\omega)(Z-i\omega(k))}
\end{equation*}
and
\begin{equation*}
	\hat{F}_{N}^{*}(\omega(k)/2\pi)\theta_{F}(k)=\frac{\hat{F}_{N}^{*}(\omega(k)/2\pi)}{1-\hat{F}_{N}^{*}(\omega(k)/2\pi)\lim_{Z\to0}\cL(C_\omega)(Z-i\omega(k))}.
\end{equation*}
Hence, $\hat{F}_{N}^{*}(\omega(k)/2\pi)\theta_{F}(k)$ converges almost everywhere to $FT(k)$. We also have
\begin{align*}
	r_r(k)&= \frac{RE(k)^2+IM(k)^2}{4|v_g(k)|^2}=\frac{|FT(k)|^2}{4|v_g(k)|^2},\\
	r_t(k)&= 1+\frac{RE(k)}{|v_g(k)|}+r_r(k)=1+\frac{\Real(FT(k))}{|v_g(k)|}+\frac{|FT(k)|^2}{4|v_g(k)|^2}.
\end{align*}
By \eqref{transmissionnew1} and \eqref{reflectionnew1}, we get that $r_t^{F_N},r_r^{F_N}$ converge almost everywhere to $r_t,r_r$ on $(0,1/2)$. By \ref{evenT} and \ref{sum=1}, the result is extended to all three rates on the whole $\T$, i.e. $r_a^{F_N},r_{t}^{F_N},r_{r}^{F_N}$ converge almost everywhere to $r_a,r_t,r_r$ on $\T$.

\appendix
\section{Proofs of the lemmas}
\subsection{Proof of lemma \ref{limLC}}
\label{prolimLC}
	Recalling \eqref{LComega}, we have
	\begin{equation*}
		\begin{split}
			\lim_{Z\to0}\cL(C_\omega)(Z-i\omega(k))=& \lim_{Z\to0}\int_{\T}\frac{Z-i\omega(k)}{(Z-i\omega(k))^2+\omega(h)}dh\\
			=& \frac{1}{2}\left( \lim_{Z\to0}\int_{\T}\frac{dh}{Z-i\omega(k)+i\omega(h)}+\lim_{Z\to0}\int_{\T}\frac{dh}{Z-i\omega(k)-i\omega(h)} \right)\\
			=& \frac{i}{2}\left( \lim_{Z\to0}\int_{\T}\frac{dh}{iZ+\omega(k)-\omega(h)}+\int_\T\frac{dh}{\omega(k)+\omega(h)} \right).
		\end{split}
	\end{equation*}
	Due to the symmetry of $\T$ in both integrals, we get 
	\begin{equation*}
		i\left( \lim_{Z\to0}\int_{0}^{1/2}\frac{dh}{iZ+\omega(k)-\omega(h)}+\int_{0}^{1/2}\frac{dh}{\omega(k)+\omega(h)} \right).
	\end{equation*}
	It is therefore enough to compute the limit of the first integral. We change the variable $\omega(h)\to w$, the integral becomes
	\begin{equation*}
		\int_{\omm}^{\omM}\frac{dw}{\omega'(\varphi_+(w))(iZ+\omega(k)-w)},
	\end{equation*}
	where $\varphi_+$ is defined in \ref{inverome}. For $u\in (\omm,\omM)$, we define
	\begin{equation*}
		C_{\omega,0}(u):=\lim_{Z\to0}\int_{\omm}^{\omM} \frac{(\omega'(\varphi_+(u))-\omega'(\varphi_+(w)))dw}{\omega'(\varphi_+(w))(iZ+u-w)}=\int_{\omm}^{\omM} \frac{(\omega'(\varphi_+(u))-\omega'(\varphi_+(w)))dw}{\omega'(\varphi_+(w))(u-w)}.
	\end{equation*}
	When both $u,w$ are close to each other and close to $\omm$ or $\omM$, then Assumption \ref{inverome} ensures that
	\begin{equation*}
		|\omega'(\varphi_+(u))-\omega'(\varphi_+(w))|=\left|\frac{1}{\varphi_+'(u)}-\frac{1}{\varphi_+'(w)}\right|\approx |\sqrt{u}-\sqrt{w}|\le |u-w|/(2\sqrt{\omm}).
	\end{equation*}
	In the case both $u,w$ are not close to the $\omm$ and $\omM$, then for $Z:|Z|<|u-w|/2$
	\begin{equation*}
		\left|\frac{\omega'(\varphi_+(u))-\omega'(\varphi_+(w))}{iZ+u-w}\right|\lesssim\left|\varphi'_+(w')\omega''(\varphi_+(w'))\right|\lesssim 1,
	\end{equation*}
	which is due to the fact that $w'$ lies between $u,w$ which are far from $\omm,\omM$. Since $\int_{\omm}^{\omM}dw/\omega'(\varphi_+(w))=\int_{0}^{1/2}dh=1/2$ we apply Dominated Convergence Theorem to see that $C_{\omega,0}$ is well-defined and bounded. Again, by Dominated Convergence Theorem, $C_{\omega,0}$ is a continuous function.

	To conclude, we compute
	\begin{align*}
		\lim_{Z\to0}\int_{\omm}^{\omM}\frac{d\omega}{iZ+\omega(k)-w}=&\lim_{Z\to0}\int_{(\omega(k)+\Imagine(Z)-\omm)/\Real(Z)}^{(\omega(k)+\Imagine(Z)-\omM)/\Real(Z)}\frac{(-w+i)d\omega}{1+w^2}\\
		=&\lim_{Z\to0}\left(\log \sqrt{1+w^2}\big|_{(\omega(k)+\Imagine(Z)-\omM)/\Real(Z)}^{(\omega(k)+\Imagine(Z)-\omm)/\Real(Z)}-i\arctan w\big|_{(\omm-\omega(k)-\Imagine(Z))/\Real(Z)}^{(\omM-\omega(k)-\Imagine(Z))/\Real(Z)} \right)\\
		=& \log\frac{\omega(k)-\omm}{\omM-\omega(k)}-i\pi.
	\end{align*}

\subsection{Proof of Lemma \ref{plusdenomlimitzero}}
\label{prodenomlimitzero}
	Changing the variables from $h$ into $k'+\varepsilon\xi'/2$ and $h'$ into $k'-\varepsilon\xi'/2$, we estimate
	\begin{align*}
		\lim_{\varepsilon\to0}&\frac{\|\hat{F}\|_{\infty}^{2}\|\tilde{\Theta}_{F}\|_{\infty}^{2}}{\varepsilon}\int_{\R}\frac{1}{(Z/2)^2+a^2}\int_{\T\times T^2_\varepsilon}\left|\frac{\hat{W}_{\varepsilon}(\xi',k')}{Z/2-ia-i\varepsilon^{-1}(\omega(k)-\omega(k'+\varepsilon\xi'/2))}\right|\\
		&\times\frac{1}{|Z/2+ia-i\varepsilon^{-1}(\omega(k)+\omega(k'-\varepsilon\xi'/2))|}d\xi'dk'dkda.
		%\label{limitplusdenomchangevari}
	\end{align*}
	Using the identity
	\begin{align*}
		\frac{1}{(Z/2-iu)(Z/2+iv)}=\left( \frac{1}{Z/2-iu}+\frac{1}{Z/2+iv} \right)\frac{1}{Z-i(u-v)}
		%\label{basicidentity}
	\end{align*}
	we only have to estimate
	\begin{align}
		\lim_{\varepsilon\to0}&\frac{\|\hat{F}\|_{\infty}^{2}\|\tilde{\Theta}_{F}\|_{\infty}^{2}}{\varepsilon}\int_{\R}\frac{1}{(Z/2)^2+a^2}\int_{\T^2\times \T_{2/\varepsilon}}\left|\frac{\hat{W}_{\varepsilon}(\xi',k')}{Z-i\varepsilon^{-1}(2\omega(k)-\omega(k'-\varepsilon\xi'/2)+\omega(k'+\varepsilon\xi'/2))}\right|\nonumber\\
		&\times\frac{1}{|Z/2-ia-i\varepsilon^{-1}(\omega(k)-\omega(k'+\varepsilon\xi'/2))|}d\xi'dk'dkda\label{limitplusdenomI1}\\
		&+\frac{\|\hat{F}\|_{\infty}^{2}\|\tilde{\Theta}_{F}\|_{\infty}^{2}}{\varepsilon}\int_{\R}\frac{1}{(Z/2)^2+a^2}\int_{\T^2\times \T_{2/\varepsilon}}\left|\frac{\hat{W}_{\varepsilon}(\xi',k')}{Z-i\varepsilon^{-1}(2\omega(k)-\omega(k'-\varepsilon\xi'/2)+\omega(k'+\varepsilon\xi'/2))}\right|\nonumber\\
		&\times\frac{1}{|Z/2+ia-i\varepsilon^{-1}(\omega(k)+\omega(k'-\varepsilon\eta'/2))|}d\eta' dk'dkda.\label{limitplusdenomI2}
	\end{align}

	Recalling that $|\xi'|<\frac{\varrho}{2^{100}\varepsilon}$ on $T^2_\varepsilon$. We choose $\varrho$ is sufficiently small, the limit \eqref{limitplusdenomI2} is $0$ as \[\varepsilon^{-1}|2\omega(k)-\omega(k'-\varepsilon\xi'/2)+\omega(k'+\varepsilon\xi'/2)|\approx |2\varepsilon^{-1}\omega(k)+\omega'(k')\xi'|\approx \varepsilon^{-1},\] and \[\varepsilon^{-1}(\omega(k)+\omega(k'+\varepsilon\xi'/2))\gtrsim \varepsilon^{-1}.\] We need to estimate \eqref{limitplusdenomI1}. To this end, we bound
	\begin{align*}
		\int_{\T}\frac{dk}{|\varepsilon Z/2-i\varepsilon a+i(\omega(k)-\omega(k'-\varepsilon\xi'/2))}&\le\sup_{E\in\R}\int_{\T}\frac{dk}{|\varepsilon Z/2-i(\omega(k)-E)|}\\
		&= 2\sup_{E\in\R}\int_{0}^{\omM}\frac{du}{|\varepsilon\Real(Z)/2-i(u-E)||\omega'(\varphi_{+}(u))|}\\
		&\lesssim \sup_{E\in\R}\int_{0}^{1}\frac{du}{|\varepsilon-i(u-E)||\omega'(\varphi_{+}(\omM u))|}\\
		&\lesssim \sup_{E\in[0,1]}\int_{0}^{1}\frac{du}{|\varepsilon-i(u-E)||\omega'(\varphi_{+}(\omM u))|}.
		%\label{estimateplusdenomI1}
	\end{align*}
	The first estimate is obtained by changing $\varepsilon a+\omega(k'-\varepsilon\xi'/2)-\varepsilon\Imagine(Z)$ into $E$, the second is because $\omega(k)$ is even so we only need to consider the positive branch. Note that $\varphi_{+}$ is defined in Assumption \ref{inverome}. Combining \eqref{minapprox} and \eqref{maxapprox}, for sufficiently small $\varepsilon$, we have
\begin{align*}
	\sup_{E\in[0,1]}&\int_{0}^{1}\frac{du}{|\varepsilon+(u-E)||\omega'(\varphi_{+}(\omM u))|}\\
	&\approx \sup_{E\in[0,1]}\int_{0}^{\varepsilon}\frac{du}{(\varepsilon+|u-E|)\sqrt{u}}+\int_{\varepsilon}^{1-\varepsilon}\frac{du}{(\varepsilon+|u-E|)\sqrt{\varepsilon}}+\int_{1-\varepsilon}^{1}\frac{du}{(\varepsilon+|u-E|)\sqrt{1-u}}\\
	&\lesssim \int_{0}^{\varepsilon}\frac{du}{\varepsilon\sqrt{u}}+\int_{\varepsilon}^{1-\varepsilon}\frac{du}{u\sqrt{\varepsilon}}+\int_{1-\varepsilon}^{1}\frac{du}{\varepsilon\sqrt{1-u}}\lesssim \varepsilon^{-1/2}\log(\varepsilon^{-1}).
	%\label{estimateplusdenomeI11}
\end{align*}
Combining \eqref{limitplusdenomI1},  with the above estimates, we find
\begin{align*}
	\int_{\R}&\frac{\|\hat{F}\|_{\infty}^{2}\|\tilde{\Theta}_{F}\|_{\infty}^{2}}{(Z/2)^2+a^2}\int_{T^2_{\varepsilon}}\frac{\varepsilon^{1/2}\log(\varepsilon^{-1})|\hat{W}_{\varepsilon}(\xi',k')|}{|2\omega(k)-\omega(k'-\varepsilon\xi'/2)+\omega(k'+\varepsilon\xi'/2)|}d\xi' dk' da\\
	&\lesssim \int_{T^2_{\varepsilon}}\frac{\varepsilon^{1/2}\log(\varepsilon^{-1})d\xi'dk'}{\langle \xi' \rangle^{3+\kappa}}\to 0.
	%\label{estimateplusdenomeI12}
\end{align*}

\subsection{Proof of Lemmas \ref{lemI1-4} and \ref{limim0}}
\label{proapp}
We next prove \eqref{approxI1} - \eqref{approxI4} for the case $+$. The domain $T^{3}_{\varepsilon,+}$ is divided into $4$ smaller domains, defined by
\begin{align*}
	T^{3}_{\varepsilon,+,(\pm)_{1},(\pm)_{2}}:=\left\{(k,\xi',\xi'')\in T^{3}_{\varepsilon,+}:(\pm)_{1}k>0,(\pm)_{2}(k-\varepsilon\xi'/2+\varepsilon\xi''/2)>0\right\}.
	%\label{T3divide}
\end{align*}
Due to the symmetry of the domains we will only provide the proofs for $T^{3}_{\varepsilon,+,+,+}$. the integrals $\cI_{\varepsilon,+},\cI^{(1)}_{\varepsilon,+},\cI^{(2)}_{\varepsilon,+},\cI^{(3)}_{\varepsilon,+}$ on this domain are denoted by $\cI,\cI^{(1)},\cI^{(2)},\cI^{(3)}$.
\subsubsection{Proof of \eqref{approxI1}}
We change the variables $\omega(k)\to u_{0},\omega(k-\varepsilon\xi'/2+\varepsilon\xi''/2)\to u_{1}, \omega(k+\varepsilon\xi'/2+\varepsilon\xi''/2)\to u_{2}$ and denote $(T^3_{\varepsilon})$ as the domain $T^3_{\varepsilon,+,+,+}$ under this change. Computing the difference gives
	\begin{align}
		&\cI(Z,\varepsilon,F)-\cI^{(1)}(Z,\varepsilon,F)\nonumber\\
		&= \frac{Zi}{8\pi \varepsilon^{2}}\int_{\R\times\T_{2/\varepsilon}}\frac{1}{(Z/2)^2+a^2}\int_{(T^3_{\varepsilon})}\frac{\hat{W}_{\varepsilon}(\varepsilon^{-1}(\varphi_{+}(u_{2})-\varphi_{+}(u_1)),(1/2)(\varphi_{+}(u_{2})+\varphi_{+}(u_1)))}{Z/2-ia-i\varepsilon^{-1}(u_{0}-u_{2})}\nonumber\\
		&\times\frac{\Lambda^{(\varepsilon)}_{+}(u_{2},u_{0},a)|\cL(F)(\varepsilon Z/2-i\varepsilon a-iu_{0})\tilde{\Theta}_{F}(\varepsilon Z/2-i\varepsilon a-iu_{0})|^{2}}{Z/2+ia+i\varepsilon^{-1}(u_{0}-u_{1})}\nonumber\\
		&\times\frac{\hat{O}^{*}(\xi,\varphi_{+}(u_{0})+\varepsilon\xi/2)}{Z+iD_{\varepsilon}^{+}(\varphi_{+}(u_{0},\xi))}\prod_{j=1}^{2}\frac{1}{\omega'(\varphi_{+}(u_{j}))}du_{0}du_{1}du_{2}d\xi da,
		\label{diffI1change}
	\end{align}
	where
\begin{align*}
	\Lambda{+}^{(\varepsilon)}(u',u,a)&:= \frac{\varepsilon^{-1}D\varphi_{+}(u',u)}{\lambda/2-ia-i\varepsilon^{-1}(u-u')-i\varepsilon^{-1}\omega'(\varphi_{+}(w))D\varphi_{+}(u',u)},\\
	%\label{Delta+def}\\
	\text{with }D\varphi_{+}(u',u)&:= \varphi_{+}(u)-\varphi_{+}(u')-\varphi'_{+}(u)(u-u').%\label{Dphi+}
\end{align*}
Our aim is to estimate \eqref{diffI1change} by $0$. Since $\cL(F)$ and $\tilde{\Theta}_{F}$ are bounded, we can ignore them in this proof.

If there exists $c_{u_1,u_2}>0$ such that $|u_{1}-u_{2}|>c_{u_1,u_2}$, then $\hat{W}_{\varepsilon}(\varepsilon^{-1}(\varphi_{+}(u_{2})-\varphi_{+}(u_{1})),\cdot)\lesssim \langle \varepsilon^{-1}\rangle^{-3-\kappa}$.  Hence, we only need to consider $|u_{1}-u_{2}|<\varrho_{2}$, where $\varrho_{2}$ is a small number we can choose. Thanks to \eqref{Trunass}, we also only consider $u_{1}, u_{2} \in [\omm+c(\rho),\omM-c(\rho)]$ for some $c(\rho)>0$; otherwise, the term vanishes. Let's choose $\varrho_{2}$ small enough so that $u_{2}\in I(2\varrho_{2})$ with $I(\varrho)=[\omm+\varrho,\omM-\varrho]$. This implies $u_{1}\in I(\varrho_{2})$. On this domain, $\left|\frac{1}{\omega'(\varphi_{+}(u_{j}))}\right|\lesssim 1,j=1,2$. Hence, we can also ignore these terms in our estimate.

Next, we further divide the domain of $u_{0}$ into $I(\varrho_3)$ as defined above and $I'(\varrho_3)=[\omm,\omM]\setminus I(\varrho_3)$, where $\varrho_3$ is a small number we choose later. When $u_{0}\in I'(\varrho_3)$, if we choose $\varrho_3$ small enough, say $\varrho_3<\varrho_{2}/2$, then $|u_{0}-u_{1}|>\varrho_{2}/2,|u_{0}-u_{2}|>\varrho_{2}$. On this domain, we use a simple estimate for $\Lambda^{(\varepsilon)}$
\begin{align}
	\left|\Lambda_{+}^{(\varepsilon)}(u_{2},u_{0},a)\right|&= \left|\frac{\varepsilon^{-1}\left( \varphi_{+}(u_{0})-\varphi_{+}(u_{2})-\varphi'_{+}(u_{0})(u_{0}-u_{2}) \right)}{Z/2-ia-i\varepsilon^{-1}\omega'(\varphi_{+}(u_{0}))\left( \varphi_{+}(u_{0})-\varphi_{+}(u_{2}) \right)}\right|\nonumber\\
	&\lesssim\varepsilon^{-1}\left( 1+\frac{1}{|\omega'(\varphi_{+}(u_{0}))|} \right).
	\label{Lambdaepsilonestimate}
\end{align}
From \eqref{diffI1change}, in $I'(\varrho_3)$ domain, applying \ref{WbracketR}, we use \eqref{asquareestimate} to get
\begin{align*}
	&\frac{1}{\varepsilon^{3}}\int_{I'(\varrho_3)}\left( 1+\frac{1}{\omega'(\varphi_{+}(u_{0}))} \right)\int_{I(\varrho_{2})}\int_{I(2\varrho_{2})}\int_{\R}\\
	&\quad\frac{\langle C\varepsilon^{-1}(u_{2}-u_{1})\rangle^{-3-\kappa}}{(1+|a+\varepsilon^{-1}(u_{0}-u_{2})|)(1+|a+\varepsilon^{-1}(u_{0}-u_{1})|)}\frac{1}{1+a^2}dadu_{1}du_{2}du_0\\
	&\lesssim\frac{1}{\varepsilon^{3}}\sum_{j=1}^{2}\int_{I'(\varrho_3)}\left( 1+\frac{1}{\omega'(\varphi_{+}(u_{0}))} \right)\int_{I(\varrho_{2})}\int_{I(2\varrho_{2})}\int_{\R}\frac{\langle C\varepsilon^{-1}(u_{2}-u_{1})\rangle^{-3-\kappa}}{(1+|a+\varepsilon^{-1}(u_{0}-u_{j})|)^2}\frac{1}{1+a^2}dadu_{1}du_{2}du_0\\
	&\lesssim\frac{1}{\varepsilon^{3}}\sum_{j=1}^{2}\int_{I'(\varrho_3)}\left( 1+\frac{1}{\omega'(\varphi_{+}(u_{0}))} \right)\int_{I(\varrho_{2})}\int_{I(2\varrho_{2})}\frac{\langle C\varepsilon^{-1}(u_{2}-u_{1})\rangle^{-3-\kappa}}{1+\varepsilon^{-2}(u_{0}-u_{j})^2}du_{1}du_{2}du_0\\
	&\lesssim\frac{1}{\varepsilon\varrho_{2}^{2}}\int_{I'(\varrho_3)}\left( 1+\frac{1}{\omega'(\varphi_{+}(u_{0}))} \right)du_{0}\int_{I(\varrho_{2})}\int_{I(2\varrho_{2})}\langle\varepsilon^{-1}(u_{2}-u_{1})\rangle^{-3-\kappa}du_{1}du_{2}\\
	&\lesssim\frac{1}{\varrho_{2}^{2}}\int_{I'(\varrho_3)}\left( 1+\frac{1}{\omega'(\varphi_{+}(u_{0}))} \right)du_{0}\\
	&\lesssim\frac{1}{\varrho_{2}^{2}}\left( \int_{I'(\varrho_3)}du_{0}+\int_{TI_{\varrho_3}}dk \right)\lesssim \rho(\varrho_3).
	%\label{Iprimeestimate}
\end{align*}
Here, $TI_{\varrho_3}\subset\T$ such that $\omega(k)\in I'(\varrho_3)$. This means that $\rho(\varrho_3)$ is small when $\varrho_3$ is small and the proof of the estimateon domain $I'(\rho_3)$ is complete.

We now work with $I(\varrho_3)$. On this domain, we have $\inf_{u_{0}\in I(\varrho_3)}\omega'(\varphi_{+}(u_{0}))>0$. Instead of the \eqref{Lambdaepsilonestimate}, a stronger estimate is used in this case. Before that, we make some changes of variables, changing $u_{1}$ into $\varepsilon^{-1}(u_{1}-u_{0})$, $u_{2}$ into $\varepsilon^{-1}(u_{2}-u_{0})$ and $a$ into $a-\varepsilon^{-1}(u_{2}-u_{0})$. Under these changes, we denote $I_{\varepsilon}(\varrho_{2})$ as the image of $I(\varrho_{2})\times I(2\varrho_{2})$, and we change the way of writing $\Lambda^{(\varepsilon)}_{+}$ into
\begin{align*}
	\tilde{\Lambda}^{(\varepsilon)}_{+}(u',u,a)&:= \frac{\varepsilon^{-1}D^{\varepsilon}\varphi_{+}(u',u)}{Z/2-ia-\varepsilon^{-1}\omega'(\varphi_{+}(w))D^{\varepsilon}\varphi_{+}(u',u)},\\%\label{tildeLambdadef}\\
	D^{\varepsilon}\omega_{+}(u',u)&:= -\int_{u}^{u+\varepsilon u'}(\omega'_{+}(v)-\omega'_{+}(u))dv=\varphi_{+}(u)+\varphi'_{+}(u)\varepsilon u'-\varphi_{+}(u+\varepsilon u').%\label{Desilonvarphidef}
\end{align*}
With these new notations, our estimate for \eqref{diffI1change} in the domain $I(\varrho_3)$ is
\begin{align}
	\int_{I(\varrho_3)}\int_{I_{\varepsilon}(\varrho_{2})}\int_{\R}\frac{\langle u_{2}-u_{1}\rangle^{-3-\kappa}|\tilde{\Delta}^{(\varepsilon)}_{+}(u_{2},u_{0},a)|}{(1+|a+u_{2}-u_{1}|)(1+(a+u_{2})^2)(1+|a|)}dadu_1du_2du_0.
	\label{estimateIvarrho3}
\end{align}
The domain of $u_{2}$ is split into
\begin{align*}
	(T^3_{\varepsilon})(\varrho_4):=\left\{u_{2}:|u_{2}|\le\varrho_4/\varepsilon\right\},
	%\label{T3epsilonrho}
\end{align*}
and $(T^3_\varepsilon)^{c}(\varrho_4)$, its complement, where $\varrho_4$ is small number to be chosen later. On $(T^3_{\varepsilon})(\varrho_4)$, we use two estimates
\begin{align}
	\int_{\R}\frac{\langle u\rangle^{-3-\kappa}}{1+|a+u|}du&\lesssim\frac{1}{1+|a|}, u=u_{2}-u_{1},\label{intduestimate}\\
	|\tilde{\Lambda}_{+}^{(\varepsilon)}(u_{2},u_{0},a)|&\lesssim|\varepsilon^{-1}D^{\varepsilon}\varphi_{+}(u_{2},u_{0})|\lesssim \varepsilon u^{2}_{2}.\nonumber%\label{tildeLambdaepsilonestimate}
\end{align}
Applying them on \eqref{estimateIvarrho3} with domain $(T^3_{\varepsilon})(\varrho_4)$, we estimate
\begin{align*}
	\int_{I(\varrho_3)}&\int_{-\varrho_4/\varepsilon}^{\varrho_4/\varepsilon}\int_{\R}\frac{\varepsilon u_{2}^{2}}{(1+(a+u_{2})^2)(1+a^{2})}dadu_2du_0\\
	&\lesssim\int_{I(\varrho_3)}\int_{-\varrho_4/\varepsilon}^{\varrho_4/\varepsilon}\frac{\varepsilon u_{2}^{2}}{1+u_{2}^{2}}du_2du_0\lesssim \varrho_4.
	%\label{estimateT3epsilonrho}
\end{align*}
For $(T^3_{\varepsilon})^{c}(\varrho_4)$, the domain is split again. We consider the following domains
\begin{align}
	\left\{(a,u_{2}):|a+u_{2}|\le|u_{2}|^{3/4}\right\}\text{ and }\left\{(a,u_{2}):|a+u_{2}|>|u_{2}|^{3/4}\right\}.
	\label{splittingT3epsilonc}
\end{align}
For the latter domain, we can again use \eqref{intduestimate} with the estimate
\begin{align*}
	\int_{\R}\frac{da}{(1+|a+b|)(1+a^2)}\lesssim\frac{1}{1+|b|}.
	%\label{aabsestimate}
\end{align*}
We estimate \eqref{estimateIvarrho3} as
\begin{align*}
	\int_{I(\varrho_3)}&\int_{\left\{(a,u_{2}):|a+u_{2}|>|u_{2}|^{3/4}\right\}}\frac{|\tilde{\Lambda}_{+}^{(\varepsilon)}(u_{2},u_{0},a)|}{(1+(a+u_{2})^2)(1+a^{2})}dadu_2du_0\\
	&\lesssim\int_{I(\varrho_3)}\int_{\varrho_4/\varepsilon}^{C_{5}/\varepsilon}\frac{|\varepsilon^{-1}D^{\varepsilon}\varphi_{+}(u_{2},u_{0})|}{(1+u_{2}^{3/2})(1+|\varepsilon^{-1}\omega'(\varphi_{+}(u_{0}))D^{\varepsilon}\varphi_{+}(u_{2},u_{0})|)}du_2du_0\\
	&\lesssim \int_{I(\varrho_3)}\int_{\varrho_4/\varepsilon}^{C_{5}/\varepsilon}\frac{1}{1+w_{2}^{3/2}}du_2du_0.
	%\label{T3epsilon2estimate}
\end{align*}
This ensures that we can use Dominated Convergence Theorem in this case. Furthermore, we have
\begin{align*}
	\varepsilon^{-1}D^{\varepsilon}\varphi_{+}(w_{2},w_{0})&=\varepsilon^{-1}(\varphi_{+}(u_{0})-\varphi_{+}(u_{0}+\varepsilon u_{2}))+\varphi_{+}'(u_{0})u_{2}\\
	&= u_{2}(\varphi_{+}'(u_{0})-\varphi_{+}'(u_{0}'))\approx 0.\text{($u_{0}'\in(u_{0},u_{0}+\varepsilon u_{2})$)}
	%\label{D^epsilonapprox}
\end{align*}
Thus, in the second domain of \eqref{splittingT3epsilonc}, the integral \eqref{estimateIvarrho3} is approximated by $0$.

Now, we consider the first domain of \eqref{splittingT3epsilonc}. We have the terms
\begin{align*}
	\phi_{\varepsilon}(u_{2},u_{0})&:= \varepsilon^{-1}\omega'(\varphi_{+}(w_{0}))D^{\varepsilon}\varphi_{+}(u_{2},u_{0})=u_{2}-\frac{\varphi_{+}(u_{0}+\varepsilon u_{2})-\varphi_{+}(u_{0})}{\varepsilon\varphi_{+}'(u_{0})},\\
	\Phi(u_{2},u_{0})&:= u_{2}-\frac{\varphi_{+}(u_{0}+u_{2})-\varphi_{+}(u_{0})}{\varphi'_{+}(u_{0})}.
	%\label{phiepsilon}
\end{align*}
We then estimate
\begin{align*}
	\int_{I(\varrho_3)}&\int_{\left\{(a,u_{2}):|a+u_{2}|\le|u_{2}|^{3/4}\right\}}\frac{|\tilde{\Lambda}_{+}^{(\varepsilon)}(u_{2},u_{0},a)|}{(1+(a+u_{2})^2)(1+a^{2})}dadu_2du_0\\
	&\lesssim\int_{I(\varrho_3)}\int_{\varrho_4/\varepsilon}^{C_5/\varepsilon}\int_{-u_{2}-u_{2}^{3/4}}^{-u_{2}+u_{2}^{3/4}}\frac{\varepsilon u_{2}^{2}}{(1+(a+u_{2})^2)(1+a^{2})(1+|a-\phi_{\varepsilon}|)}dadu_2du_0\\
	&\lesssim\int_{I(\varrho_3)}\int_{\varrho_4/\varepsilon}^{C_5/\varepsilon}\int_{-u_{2}-u_{2}^{3/4}-\phi_{\varepsilon}}^{-u_{2}+u_{2}^{3/4}-\phi_{\varepsilon}}\frac{\varepsilon u_{2}^{2}}{(1+u_{2}^2)(1+|a|)}dadu_2du_0\\
	&\lesssim\int_{I(\varrho_3)}\int_{\varrho_4/\varepsilon}^{C_5/\varepsilon}\frac{\varepsilon u_{2}^{2}}{1+u_{2}^2}\log\left( \frac{u_{2}+u_{2}^{3/4}+\phi_{\varepsilon}-1}{u_{2}-u_{2}^{3/4}+\phi_{\varepsilon}-1} \right)du_2du_0\\
	&\lesssim\int_{I(\varrho_3)}\int_{\varrho_4}^{C_5}\frac{u_{2}^{2}}{\varepsilon^{2}+u_{2}^2}\log\left( \frac{u_{2}+\varepsilon^{1/4}u_{2}^{3/4}+\Phi(u_{2},u_{0})-\varepsilon}{u_{2}-\varepsilon^{1/4}u_{2}^{3/4}+\Phi(u_{2},u_{0})-\varepsilon} \right)du_2du_0\overset{\varepsilon\to0}{\longrightarrow}0.
	%\label{T3epsilonestimate}
\end{align*}

Before ending our proof,   we remark   that the above estimate is valid if $u_{2}>0$  but the estimates for $u_{2}<0$ can be done similarly. Therefore, the proof is finished.

\subsubsection{Proof of \eqref{approxI2}}
	This proof is similar to the proof of \eqref{approxI1}, in which the roles of $u_{1}$ and $u_{2}$ switch.

\subsubsection{Proof of \eqref{approxI3}}
We state another lemma, that we will later use for \eqref{LestimateI21}.
\begin{lemma}
	If $\omega'(k_{0})=0$, we have
\begin{align*}
	\lim_{\varrho'\to0}\limsup_{\varepsilon\to0}\sup_{a\in(-\varrho',\varrho')}|\cL(F)(\varepsilon-i[a+\omega(k_{0})])\tilde{\Theta}_{F}(\varepsilon-i[a+\omega(k_{0})])|=0.
	%\label{limitk0}
\end{align*}
	\label{additionproperties}
\end{lemma}
\begin{proof}[Proof of lemma \ref{additionproperties}]

	We first prove that
\begin{align}
	\lim_{\varrho'\to0}\liminf_{\varepsilon\to0}\inf_{a\in(-\varrho',\varrho')}|\cL(C_\omega)(\varepsilon-i[a+\omega(k_{0})])|&= +\infty,\label{oldresultlimComega}\\
	\text{ or }\lim_{\varrho'\to0}\limsup_{\varepsilon\to0}\sup_{a\in(-\varrho',\varrho')}|\tilde{\Theta}(\varepsilon-i[a+\omega(k_{0})])|&= 0.\label{oldresultlimtheta}
\end{align}
Indeed, it is enough to show that for each $M$ there are $\rho',\varepsilon$ such that
\begin{equation*}
	\left|\int_0^{1/2}\frac{dh}{i\varepsilon+a+\omega(k_0)-\omega(h)}\right|>M.
\end{equation*}
The absolute value of the imaginary part of integral is
\begin{align*}
	\int_0^{1/2}\frac{\varepsilon dh}{\varepsilon^2+(a+\omega(k_0)-\omega(h))^2}&= \int_{\omm}^{\omM}\frac{\varepsilon \varphi_+'(w)dw}{\varepsilon^2+(a+\omega(k_0)-w)^2}\\
	&\gtrsim \int_{B(\omega(k_0),\rho')}\frac{\varepsilon dw}{(\rho')^{1/2}(\varepsilon^2+4(\rho')^2)}\\
	&\gtrsim \frac{\varepsilon(\rho')^{1/2}}{\varepsilon^2+(\rho')^2}.
\end{align*}
We take $\rho\gtrsim\varepsilon\gtrsim M^{-2}$. Thus, \eqref{oldresultlimComega} follows.

We recall \eqref{LaplaceFourier}, which gives $\lim_{\varepsilon\to0}\cL(F)(\varepsilon-i[a+\omega(k_{0})])=\hat{F}^{*}(\frac{a+\omega(k_{0})}{2\pi})$. Because of Assumption \ref{inequalitycondition1}, $\hat{F}^{*}(\frac{\omega(k_{0})}{2\pi})\ne0$, then $\inf_{a\in(-\varrho',\varrho')}|\cL(F)(\varepsilon-i[a+\omega(k_{0})])|>0$ when   $\varepsilon,\varrho'$ are small enough. Therefore,
\begin{align*}
	\lim_{\varrho'\to0}\liminf_{\varepsilon\to0}\inf_{a\in(-\varrho',\varrho')}|\cL(F)(\varepsilon-i[a+\omega(k_{0})])\cL(C_\omega)(\varepsilon-i[a+\omega(k_{0})])|&= +\infty\\%\label{newresultlimComega}
	\text{ or }\lim_{\varrho'\to0}\limsup_{\varepsilon\to0}\sup_{a\in(-\varrho',\varrho')}|\tilde{\Theta}_{F}(\varepsilon-i[a+\omega(k_{0})])|&= 0.%\label{newresultlimtheta}
\end{align*}
\end{proof}
	We first write down the difference
	\begin{align*}
		\cI^{(2)}(Z,\varepsilon,F)-\cI^{(3)}(Z,\varepsilon,F)&= \frac{Z}{32\pi}\int_{\R\times\T_{2/\varepsilon}}\frac{1}{(Z/2)^2+a^2}\int_{T^{3}_{\varepsilon,+,+,+}}\frac{\hat{W}_{\varepsilon}(\xi',k+\varepsilon\xi''/2)}{Z/2-ia+i\omega'(k)(\xi'+\xi'')/2}\\
		&\times\frac{D\theta_{F,\varepsilon}(a,k)\hat{O}^{*}(\xi,k+\varepsilon\xi/2)}{(Z/2+ia+i\omega'(k)(\xi'-\xi'')/2)(Z+iD_{\varepsilon}^{+}(k,\xi))}d\xi''d\xi'dkd\xi da,%\label{DiffI3}
	\end{align*}
	where $D\theta_{F,\varepsilon}(a,k)=\left|\cL(F)(\varepsilon Z/2-i\varepsilon a-i\omega(k))\tilde{\Theta}(\varepsilon Z/2-i\varepsilon a-i\omega(k))\right|^{2}-\left|\hat{F}(\omega(k)/(2\pi))\theta_{F}(k)\right|^{2}$.

	By \eqref{Trunass}, if $k+\varepsilon\xi''\in [\T](\varrho)$, then both $\cI^{(2)}$ and $\cI^{(3)}$ are zero on this domain. We only need to consider $k+\varepsilon\xi''\in [\T]^{c}(\varrho)$.

	If $k\in [\T]^{c}(\varrho_5)$ for a small constant $\varrho_5>0$, then $\inf\omega'(k)>0$. We have the following estimate
	\begin{align}
		\int_{\R}\frac{1}{(Z/2)^{2}+a^2}&\int_{[\T]^{c}(\varrho_5)}\int_{\R^{2}}\frac{\langle\xi'\rangle^{-3-\kappa}|D\theta_{F,\varepsilon}(a,k)|}{(Z/2+i\omega'(k)\xi'/2)^{2}+(a+\omega'(k)\xi''/2)^{2}}d\xi''d\xi'dkda\nonumber\\
		&\lesssim\int_{\R}\frac{1}{1+a^2}\int_{[\T]^{c}(\varrho_5)}\int_{\R^{2}}\frac{\langle\xi'\rangle^{-3-\kappa}|D\theta_{F,\varepsilon}(a,k)|}{1+(a+\omega'(k)\xi''/2)^{2}}d\xi''d\xi'dkda\nonumber\\
		&\lesssim\int_{[\T]^{c}(\varrho_5)}\int_{\R}\frac{1}{|\omega'(k)|}\frac{C_6}{1+(\xi'')^{2}}d\xi''dk\lesssim\frac{C_6}{\inf \omega'(k)}.
		\label{LcestimateDiffI3}
	\end{align}
	Here, $C_6$ is a bound for $D\theta_{F,\varepsilon}$ when $\varepsilon$ is small enough. Therefore, we can use Dominated Convergence Theorem to prove that  the difference goes to $0$.

	We consider $k\in [\T](\varrho_5)$, we can see that $\cI^{(3)}$ on this domain is $0$. Indeed, we can make the same estimation like in $[\T]^{c}$,
	\begin{align}
		\int_{\R}\frac{1}{(Z/2)^{2}+a^2}&\int_{[\T](\varrho_5)}\int_{\R^{2}}\frac{\langle\xi'\rangle^{-3-\kappa}\left|\hat{F}(\omega(k)/(2\pi))\theta_{F}(k)\right|^{2}}{(\lambda/2+i\omega'(k)\xi'/2)^{2}+(a+\omega'(k)\xi''/2)^{2}}d\xi''d\xi'dkda\nonumber\\
	&\lesssim\int_{\R}\frac{1}{1+a^2}\int_{[\T](\varrho_5)}\int_{\R^{2}}\frac{\langle\xi'\rangle^{-3-\kappa}\left|\hat{F}(\omega(k)/(2\pi))\theta_{F}(k)\right|^{2}}{1+(a+\omega'(k)\xi''/2)^{2}}d\xi''d\xi'dkda\nonumber\\
	&\lesssim\int_{[\T](\varrho_5)}\int_{\R}\frac{1}{|\omega'(k)|}\frac{\left|\hat{F}(\omega(k)/(2\pi))\theta_{F}(k)\right|^{2}}{1+(\xi'')^{2}}d\xi''dk\lesssim\varrho_5.
		\label{LestimateI3}
	\end{align}
	The last estimate comes from Lemma \ref{thetaFprop}, which leads to $\left|\hat{F}(\omega(k)/(2\pi))\theta_{F}(k)\right|^{2}\lesssim |\omega'(k)|$.

	For a small number $\varrho'$ and a small $\varepsilon$, by Lemma \ref{additionproperties}, there is small number $\rho(\varrho')$ such that
	\begin{align}
		\int_{\R}\frac{1}{(Z/2)^{2}+a^2}&\int_{[\T](\varrho_5)}\int_{(\varrho-\varrho_5)/\varepsilon\le|\xi''|\le2/\varepsilon}\frac{\left|\cL(F)(\varepsilon Z/2-i\varepsilon a-i\omega(k))\tilde{\Theta}(\varepsilon Z/2-i\varepsilon a-i\omega(k))\right|^{2}}{(Z/2+i\omega'(k)\xi'/2)^{2}+(a+\omega'(k)\xi''/2)^{2}}d\xi''dkda\nonumber\\
		&\lesssim\int_{-\varrho'/\varepsilon}^{\varrho'/\varepsilon}\frac{1}{1+a^2}\int_{[\T](\varrho_5)}\int_{(\varrho-\varrho_5)/\varepsilon\le|\xi''|\le2/\varepsilon}\frac{\rho(\varrho')}{1+(a+\omega'(k)\xi''/2)^{2}}d\xi''dkda\label{LestimateI21}\\
		&\quad+\int_{|a|>\varrho/\varepsilon}\frac{1}{1+a^2}\int_{[\T](\varrho_5)}\int_{(\varrho-\varrho_5)/\varepsilon\le|\xi''|\le2/\varepsilon}\frac{1}{1+(a+\omega'(k)\xi''/2)^{2}}d\xi''dkda.
		\label{LestimateI22}
	\end{align}
	If $\varrho_5$ is small enough, then we have $|\omega'(k)||\xi''|\le\frac{\varrho}{\varepsilon}\le|a|$. Hence, $|a+\omega'(k)\xi''/2|>|a/2|$. We obtain an estimate for \eqref{LestimateI22}
\begin{align*}
	\frac{\varrho_5}{\varepsilon}\int_{|a|>\frac{\varrho}{\varepsilon}}\frac{da}{1+a^4}\lesssim \varrho_5\int_{|a|>\varrho}\frac{da}{\varepsilon+\varepsilon^{-3}a^{4}}\lesssim \varrho_5\varepsilon^{3}.
	%\label{LestimateI22cont}
\end{align*}
	For \eqref{LestimateI21}, we have
	\begin{align*}
		\rho(\varrho')\int_{[\T](\varrho_5)}&\int_{|\omega'(k)|(\varrho-\varrho_5)/\varepsilon\le|\xi''|\le2|\omega'(k)|/\varepsilon}\frac{1}{|\omega'(k)|}\frac{1}{1+(\xi'')^2}d\xi''dk\\
		&\lesssim \rho(\varrho')\int_{[\T](\varrho_5)}\frac{\arctan(2|\omega'(k)|/\varepsilon)-\arctan(|\omega'(k)|(\varrho-\varrho_5)/\varepsilon)}{|\omega'(k)|}dk\\
		&\lesssim \rho(\varrho')\int_{[\T](\varrho_5)}\frac{\arctan\left( \frac{(2-\varrho+\varrho_5)|\omega'(k)|}{\varepsilon}\left( 1+\frac{2(\varrho-\varrho_5)|\omega'(k)|^2}{\varepsilon^{2}} \right)^{-1} \right)}{|\omega'(k)|}dk\\
		&\lesssim \rho(\varrho')\int_{[\T](\varrho_5)}\frac{(2-\varrho+\varrho_5)}{\varepsilon}\left( 1+\frac{2(\varrho-\varrho_5)|\omega'(k)|^2}{\varepsilon^{2}} \right)^{-1}dk\\
		&\lesssim \rho(\varrho')\int_{\omega'([\T](\varrho_5))/\varepsilon}\frac{dk}{1+k^2}\lesssim\rho(\varrho')\int_{\R}\frac{dk}{1+k^2}\lesssim\rho(\varrho')
		%\label{LestiamteI21cont}
	\end{align*}

\subsubsection{Proof of \eqref{approxI4}}
	As $\hat{O}^{*}$ is a Schwartz function, we approximate $\hat{O}^{*}(\eta,k+\varepsilon\xi/2)$ by $\hat{O}^{*}(\xi,k)$ and $D_{\varepsilon}^{+}(k,\xi)$ by $\omega'(k)\xi$.
	Hence, we approximate $\cI^{(3)}(Z,\varepsilon,F)$ with
	\begin{align*}
		\tilde{\cI}^{(3)}(Z,\varepsilon,F)&= \frac{Z}{32\pi}\int_{\R\times\T_{2/\varepsilon}}\frac{1}{(Z/2)^2+a^2}\int_{T^{3}_{\varepsilon,+,+,+}}\frac{\hat{W}_{\varepsilon}(\xi',k+\varepsilon\xi''/2)}{Z/2-ia+i\omega'(k)(\xi'+\xi'')/2}\\
		&\times\frac{|\hat{F}(\omega(k)/(2\pi))\theta_{F}(k)|^{2}\hat{O}^{*}(\xi,k)}{(Z/2+ia+i\omega'(k)(\xi'-\xi'')/2)(Z+i\omega'(k)\xi)}d\xi''d\xi'dkd\xi da.
		%\label{tildeI3}
	\end{align*}
	We change the variable $k$ into $k-\varepsilon\xi''/2$ and denote by $[T^3_\varepsilon]$ the image of $T^{3}_{\varepsilon,+,+,+}$ under this change. The integral $\tilde{\cI}^{(3)}$ becomes
	\begin{align*}
		\tilde{\cI}^{(3)}(Z,\varepsilon,F)&= \frac{Z}{32\pi}\int_{\R\times\T_{2/\varepsilon}}\frac{1}{(Z/2)^2+a^2}\int_{[T^3_\varepsilon]}\frac{\hat{W}_{\varepsilon}(\xi',k)}{Z/2-ia+i\omega'(k-\varepsilon\xi''/2)(\xi'+\xi'')/2}\\
		&\times\frac{|\hat{F}(\omega(k-\varepsilon\xi''/2)/(2\pi))\theta_{F}(k-\varepsilon\xi''/2)|^{2}\hat{O}^{*}(\xi,k-\varepsilon\xi''/2)}{(Z/2+ia+i\omega'(k-\varepsilon\xi''/2)(\xi'-\xi'')/2)(Z+i\omega'(k-\varepsilon\xi''/2)\xi)}d\xi''d\xi'dkd\xi da
		%\label{tildeI3change}
	\end{align*}
	Once again, we use the regularity of $\hat{O}^{*}$ to estimate $\hat{O}^{*}(\xi,k-\varepsilon\xi''/2)$ by $\hat{O}^{*}(\xi,k)$ and $\omega'(k-\varepsilon\xi''/2)\eta$ by $\omega'(k)\xi$. Now, we approximate $\tilde{\cI}^{(3)}$ with
	\begin{align*}
		\tilde{\tilde{\cI}}^{(3)}(Z,\varepsilon,F)&= \frac{Z}{32\pi}\int_{\R\times\T_{2/\varepsilon}}\frac{1}{(Z/2)^2+a^2}\int_{[T^3_\varepsilon]}\frac{\hat{W}_{\varepsilon}(\xi',k)}{Z/2-ia+i\omega'(k-\varepsilon\xi''/2)(\xi'+\xi'')/2}\\
		&\times\frac{|\hat{F}(\omega(k-\varepsilon\xi''/2)/(2\pi))\theta_{F}(k-\varepsilon\xi''/2)|^{2}\hat{O}^{*}(\xi,k)}{(Z/2+ia+i\omega'(k-\varepsilon\xi''/2)(\xi'-\xi'')/2)(Z+i\omega'(k)\xi)}d\xi''d\xi'dkd\xi da.
		%\label{tildetildeI3}
	\end{align*}
	Then, we will approximate $\tilde{\tilde{\cI}}^{(3)}$ with 
	\begin{align*}
		\tilde{\tilde{\tilde{\cI}}}^{(3)}(Z,\varepsilon,F)&= \frac{Z}{32\pi}\int_{\R\times\T_{2/\varepsilon}}\frac{1}{(Z/2)^2+a^2}\int_{[T^3_\varepsilon]}\frac{\hat{W}_{\varepsilon}(\xi',k)}{Z/2-ia+i\omega'(k)(\xi'+\xi'')/2}\\
		&\times\frac{|\hat{F}(\omega(k)/(2\pi))\theta_{F}(k)|^{2}\hat{O}^{*}(\xi,k)}{(Z/2+ia+i\omega'(k)(\xi'-\xi'')/2)(Z+i\omega'(k)\xi)}d\xi''d\xi'dkd\xi da.
		%\label{tilde3I3}
	\end{align*}
	We denote
	\begin{align*}
		d_{\varepsilon}(k,\xi',\xi'')&:= \frac{|\hat{F}(\omega(k-\varepsilon\xi''/2)/(2\pi))\theta_{F}(k-\varepsilon\xi''/2)|^{2}}{(Z/2-ia+i\omega'(k-\varepsilon\xi''/2)(\xi'+\xi'')/2)(Z/2+ia+i\omega'(k-\varepsilon\xi''/2)(\xi'-\xi'')/2)}\\
		&\quad-\frac{|\hat{F}(\omega(k)/(2\pi))\theta_{F}(k)|^{2}}{(Z/2-ia+i\omega'(k)(\xi'+\xi'')/2)(Z/2+ia+i\omega'(k)(\xi'-\xi'')/2)}.
		%\label{depsilon}
	\end{align*}
	This difference converges to $0$ as $\varepsilon \to 0$ and
	\begin{align*}
		|d_{\varepsilon}|&\le\frac{|\hat{F}(\omega(k-\varepsilon\xi''/2)/(2\pi))\theta_{F}(k-\varepsilon\xi''/2)|^{2}}{(Z/2+i\omega'(k-\varepsilon\xi''/2)\xi'/2)^{2}+(a+\omega'(k-\varepsilon\xi''/2)\xi''/2)^{2}}\\
		&\quad+\frac{|\hat{F}(\omega(k)/(2\pi))\theta_{F}(k)|^{2}}{(Z/2+i\omega'(k)\xi'/2)^{2}+(a+\omega'(k)\xi''/2)^{2}}.
		%\label{depsilonestimate}
	\end{align*}
	To reach the desired result, we can repeat the techniques used in \eqref{LcestimateDiffI3} and \eqref{LestimateI3}. We consider $k\in [\T]^{c}(\varrho-\varrho/2^{100})$, because $\hat{W}_{\varepsilon}(\xi',k)$ vanishes otherwise. We consider two domains $k-\varepsilon\xi''\in [\T]^{c}(\varrho_6)$ and $k-\varepsilon\xi''\in [\T](\varrho_6)$. The techniques in \eqref{LcestimateDiffI3} is reused to bound the terms concerning the first domain by Dominated Convergence Theorem. The techniques in  \eqref{LestimateI3} is reused to choose the small parameter $\varrho_6$.

	For the final step in this proof, we approximate $\tilde{\tilde{\tilde{\cI}}}^{(3)}$ with $\cI^{(4)}$. This can is done by using Dominated Convergence Theorem to turn $\T_{2/\varepsilon}$ into $\R$ and $[T^3_\varepsilon]$ into $\T\times\R^2$ in the limit of $\varepsilon \to 0$.

\subsubsection{Proof of Lemma \ref{limim0}}
	We compute $\Imagine S_{\varepsilon}^{II}$ using similar computations used to obtain \eqref{sumreal2cosnew} and \eqref{sumreal2sinnew}
\begin{align*}
	&\Imagine S^{II}_{\varepsilon}(Z,k,F)\\
	&= \varepsilon^{2}Z\int_{0}^{t}e^{-\varepsilon Zt}\left\langle\left( \int_{0}^{t}\cos(\omega(k)s)(\Theta_{F}\star F\star \Psi)(s)ds \right) \left(\int_{0}^{t}\sin(\omega(k)s)(\Theta_{F}\star F\star \Psi)(s)ds\right)\right\rangle_{\mu_{\varepsilon}}dt.
	%\label{Imagine2}
\end{align*}

Similar computations used to obtain \eqref{approxI1} - \eqref{approxI4} can be reused and we get
\begin{align}
	\lim_{\varepsilon\to0}&\int_{\R\times\T}\Imagine S_{\varepsilon}^{II}(Z,k,F)\frac{\hat{O}^{*}(\xi,k+\varepsilon\xi/2)}{Z+iD_{\varepsilon}^{+}(k,\xi)}d\xi dk\label{Imaginelimit}\\
	&= -\int_{\T\times\R^{2}}\frac{|\hat{F}(\omega(k)/(2\pi))|^{2}|\theta_{F}(k)|^{2}\hat{W}_{\varepsilon}(\xi',k)\hat{O}(\xi,k)}{16i|v_g(k)|(Z+i\omega(k)\xi')(Z+i\omega'(k)\xi)}d\xi'd\xi dk\nonumber\\
	&-\int_{\T\times\R^{2}}\frac{|\hat{F}(\omega(k)/(2\pi))|^{2}|\theta_{F}(k)|^{2}\hat{W}_{\varepsilon}(\xi',-k)\hat{O}(\xi,k)}{16i|v_g(k)|(Z-i\omega(k)\xi')(Z+i\omega'(k)\xi)}d\xi'd\xi dk.  \nonumber
\end{align}
On the other hand, the limit when $\varepsilon \to 0$ of $\int_{\R\times\T}\Imagine S_{\varepsilon}^{II}(Z,k,F)\frac{\hat{O}^{*}(\xi,k-\varepsilon\xi/2)}{Z+iD_{\varepsilon}^{-}(k,\xi)}$ gives exactly the same result with \eqref{Imaginelimit}. Thus, the limit of the difference is $0$.

\subsection{An Ito Analysis}
\subsubsection{Proof of Lemma \ref{ligrowimpul}}
\label{proligrowimpul}
Consider the SDE (stochastic differential equation): 
\begin{align*}
	d\hat{\psi}(t,k)&=  -i\omega(k)\hat{\psi}(t,k)dt-i\nu \alpha_0(t)dt+iF(t)dt+\sqrt{2\nu T}idR(t).
	%\label{dpsi2}
\end{align*}
Taking the functional $g(\hat{\psi},\hat{\psi}^*)= \| \psi \|^2_{L^2(\mathbb{T})}$, and applying the Ito-formula we get
%(We refer to \cite{karatzas2012brownian} Chapter 5, 
%\cite{da2014stochastic} Chapter 4, for wellposedness and applicability of Ito formula):
%\textcolor{blue}{Making this rigorous requires a "correct" reference: this is obvious in finite-dimensional spaces, we can either discretize and pass to the limit, or directly use infinite-dimensional ito formula e.g. Kuksin book}): 
\begin{align*}
    d\|\hat{\psi}(t) \|^2_{L^2(\mathbb{T})} = & 
    \left[\int_{\mathbb{T}} \hat{\psi}^*(t,k) \left(-i \omega(k) 
    \hat{\psi}(k,t) - i \nu \alpha_0(t) + iF(t) \right)\right]dt    \\
    &+ \left[\int_{\mathbb{T}} \hat{\psi}(k,t) \left(i \omega(k) 
    \hat{\psi}^*(k,t) + i \nu \alpha_0(t) - iF(t) \right) \right] dt   \\
    &- \left[i \sqrt{2 \nu T} \int_{\mathbb{T}}
    (\hat{\psi}(k,t)- \hat{\psi}^*(k,t))
    \right] d R(t)   \\
    &+ \left[ 2 \nu T \right] dt.
\end{align*}
Using the definition of $\alpha_0=\int Im( \psi)$ we get %(\textcolor{blue}{check for $2$ and signs}): 
\begin{align*}
    d\|\hat{\psi}(t) \|^2_{L^2(\mathbb{T})} = 
    \left[-2 \nu \alpha_0^2(t) + 2 \alpha_0(t) F(t)
    + 2 \nu T \right] dt + \sqrt{2 \nu T} \alpha_0(t) dR(t).
\end{align*}
Upon a time rescaling $t \to t/\varepsilon$, using scaling 
properties of the Wiener process we get: 
\begin{equation*}
     d\|\hat{\psi}^{(\varepsilon)}(t) \|^2_{L^2(\mathbb{T})} = 
    \left[-\frac{2 \nu}{\varepsilon}  
    [\alpha_0^{(\varepsilon)}]^2(t) + \frac{2}{\varepsilon} 
    \alpha_0^{(\varepsilon)}(t) F^{(\varepsilon)}(t)
    + \frac{2 \nu T}{\varepsilon} \right] dt + \sqrt{\frac{2 
    \nu T}{\varepsilon}} \alpha_0(t) dR(t).
\end{equation*}
We would like bounds of the form: 
\begin{equation} \label{eq:boundtme}
    \sup_{\varepsilon \in (0,1]} \varepsilon \mathbb{E}_{\varepsilon}\| \hat{\psi}^{(\varepsilon)}(t)\|_{L^2(\mathbb{T})}^2 \leq    \sup_{\varepsilon \in (0,1]} 
    \varepsilon 
    \mathbb{E}_{\varepsilon} \| \hat{\psi}^{(\varepsilon)}(0)\|_{L^2(\mathbb{T})}^2+Ct,
\end{equation}
for a constant $C>0$. 
For this, it is sufficient to have bounds of the form: 
\begin{equation*}
    \sup_{\varepsilon \in (0,1]} \frac1{\nu}\int_0^{t/\varepsilon} \varepsilon F^2(s) ds \leq C't, 
\end{equation*}
%\textcolor{red}{We should check if such an assumption is compatible with other assumptions and scaling concerning $F$}
as stated in Assumption \ref{integralF^2bound}.

\subsubsection{Proof of Lemma \ref{linearfeed}}
\label{prolinearfeed}
We consider
\begin{align*}
	d\hat{\psi}(t,k)&= -i\omega(k)\hat{\psi}(t,k)dt+iF\star \alpha_{0}(t)dt+i\sqrt{2\nu T}dR(t).
	%\label{newdpsi2}
\end{align*}
Thanks the It\^{o} formula we get
\begin{align*}
    d\|\hat{\psi}(t) \|^2_{L^2(\mathbb{T})} = & 
    \left[\int_{\mathbb{T}} \hat{\psi}^*(t,k) \left(-i \omega(k) 
    \hat{\psi}(k,t) +iF\star \alpha_{0}(t) \right)+\int_{\mathbb{T}} \hat{\psi}(k,t) \left(i \omega(k) 
    \hat{\psi}^*(k,t) - iF\star \alpha_{0}(t) \right)\right]dt \\
    & -\left[i \sqrt{2 \nu T} \int_{\mathbb{T}}
    (\hat{\psi}(k,t)- \hat{\psi}^*(k,t))
    \right] d R(t) + 
      2 \nu T  dt\\
    = &\left[ 2F\star \alpha_{0}(t)\alpha_{0}(t)+2\nu T \right]dt+\sqrt{2\nu T}\alpha_{0}(t)dR(t).
\end{align*}
Rescaling  time by $1/\varepsilon$ and using scaling properties of $R(t)$ yields: 
\begin{equation*}
     d\|\hat{\psi}^{(\varepsilon)}(t) \|^2_{L^2(\mathbb{T})} = 
    \left[\frac{2}{\varepsilon}  
	    F\star \alpha_{0}(t/\varepsilon)\alpha_0(t/\varepsilon) 
    + \frac{2 \nu T}{\varepsilon} \right] dt + \sqrt{\frac{2 
    \nu T}{\varepsilon}} \alpha_0(t/\varepsilon) dR(t).
\end{equation*}
By solving  \eqref{alpha_0Comeganew}, $\alpha_0$ is given by
\begin{align*}
	\alpha_{0}(t/\varepsilon)=\Theta_{F}\star \Psi(t/\varepsilon)+\sqrt{2\nu T}\Theta_{F}\star C_\omega\star dR(t/\varepsilon).
\end{align*}
Combining above expressions, and applying the It\^{o} formula we get
\begin{align*}
\varepsilon \mathbb{E}_{\varepsilon}\| \hat{\psi}^{(\varepsilon)}(t)\|_{L^2(\mathbb{T})}^2&=
\varepsilon \mathbb{E}_{\varepsilon}\| \hat{\psi}^{(\varepsilon)}(0)\|_{L^2(\mathbb{T})}^2 + 2 \nu Tt  + 
2\int_0^t  \E_{\varepsilon}\left[F\star \Theta_{F}\star \Psi(s/\varepsilon)\Theta_{F}\star \Psi(s/\varepsilon) \right]ds  \nonumber\\
&\quad+  4\nu T\int_{0}^{t}\int_{0}^{s/\varepsilon}F\star \Theta_{F}\star C_\omega(s/\varepsilon-s_1)\Theta_{F}\star C_\omega(s/\varepsilon-s_1)ds_1ds\\
&\quad+2 \nu T\sqrt{\varepsilon}\E_{\varepsilon}\left[ \int_{0}^{t}\int_{0}^{s/\varepsilon}\Theta_{F}\star C_\omega(s/\varepsilon-s_1)dR(s_1)dR(s) \right].
\end{align*}
In the above expression, we used the fact that terms with single $dR$ has expectation zero thanks to properties of the Wiener process. 
Moreover, the second term is obtained by $\E\left[ dw(s)dw(s') \right]=\delta(s-s')$ (in the distribution sense). Using the mentioned covariance structure we also get
\begin{align*}
	\sqrt{\varepsilon}\E_{\varepsilon}\left[ \int_{0}^{t}\int_{0}^{s/\varepsilon}\Theta_{F}\star C_\omega(s/\varepsilon-s_1)dR(s_1)dR(s) \right]&= \E_{\varepsilon}\left[ \int_{0}^{t}\int_{0}^{s/\varepsilon}\Theta_{F}\star C_\omega(s/\varepsilon-s_1)dR(s_1/\varepsilon)dR(s) \right]  \\ \leq \Theta_F \star C_\omega(0) t=Ct
 %\nonumber\\
%	&= \E_{\varepsilon}\left[ \int_{0}^{t}dw(s)\int_{0}^{s}dw(s_1)g_{F}\star J(s/\varepsilon-s_1/\varepsilon) \right]\nonumber\\
%	&= \int_{0}^{t}g_{F}\star J(0)ds=t.
\end{align*}

To obtain a linear in-time growth bound as in \eqref{eq:boundtme}, it is sufficient to bound the following expression uniformly in time 
\begin{align}
	\underbrace{\E_{\varepsilon}\left[ F\star \Theta_{F}\star \Psi(t/\varepsilon)\Theta_{F}\star \Psi(t/\varepsilon) \right]}_{(a)}+ 2\nu T\underbrace{\int_{0}^{t/\varepsilon}F\star \Theta_{F}\star C_\omega(s)\Theta_{F}\star C_\omega(s)ds}_{(b)}.\label{boundforgrowth}
\end{align}
To this end, we approximate the terms using the Laplace transform. For term $(a)$ in \eqref{boundforgrowth},
\begin{align}
	\E_{\varepsilon}&\left[ F\star \Theta_{F}\star \Psi(t/\varepsilon)\Theta_{F}\star \Psi(t/\varepsilon) \right]\nonumber\\
%	&= \E_{\varepsilon}\left[ \int_{\T}dk\Imagine\left( \hat{\psi}(k)F\star g_{F}\star e^{-i\omega(k)\cdot}(t/\varepsilon) \right)\int_{\T}dk\Imagine\left( \hat{\psi}(k)g_{F}\star e^{-i\omega(k)\cdot}(t/\varepsilon) \right) \right] ,
% \nonumber\\
	&= \frac{1}{4}\E_{\varepsilon}\left[ \int_{\T} \hat{\psi}(k)F\star \Theta_{F}\star e^{-i\omega(k)\cdot}(t/\varepsilon)dk \int_{\T}\hat{\psi}^{*}(k)\Theta_{F}\star e^{i\omega(k)\cdot}(t/\varepsilon)dk \right] + c.c.
 %\nonumber\\
%	&\qquad +\frac{1}{4}\E_{\varepsilon}\left[ \int_{\T}dk \hat{\psi}^{*}(k)F\star g_{F}\star e^{i\omega(k)\cdot}(t/\varepsilon) \int_{\T}dk\hat{\psi}(k)g_{F}\star e^{-i\omega(k)\cdot}(t/\varepsilon) \right].
	\label{Psiproduct}
\end{align}
where c.c. stands for the complex conjugate. Thanks to the assumption of $F$, this term will be negative when $\varepsilon$ is small enough. Indeed, taking a Laplace transform and by some manipulations   we have
\begin{align}
	\cL&\left( \E_{\varepsilon}\left[ \int_{\T} \hat{\psi}(k)F\star \Theta_{F}\star e^{-i\omega(k)\cdot}(t/\varepsilon)dk \int_{\T}\hat{\psi}^{*}(k)\Theta_{F}\star e^{i\omega(k)\cdot}(t/\varepsilon)dk \right] \right)(Z)\nonumber\\
	%&= \int_{\T^2}dkd\ell\E_{\varepsilon}[\hat{\psi}(k)\hat{\psi}(\ell)]\frac{\varepsilon}{\lambda}\cL\left( F\star g_{F}\star e^{-i\omega(k)\cdot}(t)g_{F}\star e^{i\omega(\ell)\cdot}(t) \right)(\varepsilon\lambda)\nonumber\\
	%&= \int_{\T^2}dkd\ell\E_{\varepsilon}[\hat{\psi}(k)\hat{\psi}(\ell)]\frac{1}{\lambda}\cL\left( F\star g_{F}\star e^{i\omega(k)\cdot}(t/\varepsilon)g_{F}\star e^{i\omega(\ell)\cdot}(t/\varepsilon) \right)(\lambda)\nonumber\\
	&= \int_{\T^2}\E_{\varepsilon}[\hat{\psi}(k)\hat{\psi}^*(h)]\frac{1}{Z}\cL\left( \int_{0}^{t/\varepsilon} e^{i\omega(k)s}F\star \Theta_{F}(s)e^{-i\omega(k)t/\varepsilon}ds\int_{0}^{t/\varepsilon}e^{-i\omega(h)s'}\Theta_{F}(s') e^{i\omega(\ell)t/\varepsilon}ds' \right)(Z)dkdh\nonumber\\
	&= \int_{\T^{2}}\E_{\varepsilon}[\hat{\psi}(k)\hat{\psi}^*(h)]\frac{\varepsilon^2}{2Z\pi i}\lim_{L\to\infty}\int_{c-iL}^{c+iL}\cL\left( \int_{0}^{t} e^{i\omega(k)s}F\star \Theta_{F}(s)ds \right)(\varepsilon\zeta)\nonumber\\
	&\hspace{0.2\textwidth}\times\cL\left( \int_{0}^{t}e^{-i\omega(h)s'}\Theta_{F}(s')ds'  \right)(\varepsilon(Z-\zeta)+i(\omega(k)-\omega(h)))d\zeta dkdh\nonumber\\
	&= \int_{\T^{2}}\frac{\E_{\varepsilon}[\hat{\psi}(k)\hat{\psi}^*(h)]}{2Z\pi i}\lim_{L\to\infty}\int_{c-iL}^{c+iL}\frac{\cL(F)(\varepsilon\zeta-i\omega(k))\tilde{\Theta}_{F}(\varepsilon\zeta-i\omega(k))}{\zeta}\nonumber\\
	&\hspace{0.2\textwidth}\times\frac{\tilde{\Theta}_{F}(\varepsilon(Z-\zeta)+i\omega(k))}{(Z-\zeta)+i\varepsilon^{-1}(\omega(k)-\omega(h))}d\zeta dk dh.
	\label{Psiproductcal}
\end{align}
We change the variable from $k$ into $k+\varepsilon\xi/2$, $h$ into $k-\varepsilon\xi/2$, \eqref{Psiproductcal} becomes
\begin{align*}
	&\int_{\T\times\T_{2/\varepsilon}}\frac{\hat{W}_{\varepsilon}(\xi,k)}{Z\pi i}\lim_{L\to\infty}\int_{c-iL}^{c+iL}\frac{\cL(F)(\varepsilon\zeta-i\omega(k))\tilde{\Theta}_{F}(\varepsilon\zeta-i\omega(k))\tilde{\Theta}_{F}(\varepsilon(Z-\zeta)+i\omega(k))}{\zeta(Z-\zeta+iD_{\varepsilon}(k,\xi))}d\zeta d\xi dk\\
	&\overset{\varepsilon\to0}{\longrightarrow}\int_{\T\times\R}\frac{\hat{W}_{0}(\xi,k)}{Z\pi i}\lim_{L\to\infty}\int_{c-iL}^{c+iL}\frac{\hat{F}^{*}(\omega(k)/(2\pi))|\theta_{F}(k)|^{2}}{\zeta(Z-\zeta+i\omega'(k)\xi)}d\zeta d\xi dk\\
	&= \int_{\T\times\R}\frac{2\hat{W}_{0}(\xi,k)}{Z(Z+i\omega'(k)\xi)}\hat{F}^{*}(\omega(k)/(2\pi))|\theta_{F}(k)|^{2}d\xi dk.
	%\label{Psiproductlimit}
\end{align*}
The limit can be obtained using Dominated Convergence Theorem. Indeed, $\cL(F), \tilde{\Theta}_{F}$ are both bounded when $\varepsilon$ is small; $Z-\zeta+iD_{\varepsilon}(k,\xi)$ is bounded by $\Real(Z)-c$; $\lim_{L\to\infty}\int_{c-iL}^{c+iL}1/\sigma d\sigma$ is bounded using \eqref{integralfraction}; finally, Assumption \ref{WbracketR} gives a dominating function $\langle\xi\rangle^{-3-\kappa}$ for $1_{[-1/\varepsilon,1/\varepsilon]}\hat{W}_{0}(\xi,k)$. Then, a similar estimation can be made for the complex conjugate term (c.c.). By using Parseval identity like in \eqref{thermal4limitshowing}, \eqref{Psiproduct} is approximated by
\begin{align*}
	\int_{\T\times\R}W_{0}(x,k)\Real(\hat{F}(\omega(k)/(2\pi)))|\theta_{F}(k)|^{2}1_{[0,v_g(k)]}(x)dxdk.
\end{align*}
This term is negative due to Assumption \ref{inequalitycondition1}. Hence, when $\varepsilon$ is small enough, $\E_{\varepsilon}[F\star \Theta_{F}\star\Psi(t/\varepsilon)\Theta_{F}\star\Psi(t/\varepsilon)]$ is negative.

In the same way, we now deal with term $(b)$ in \eqref{boundforgrowth}. We write
\begin{align*}
	F\star \Theta_{F}\star C_\omega(s)\Theta_{F}\star C_\omega(s)
 %&= \frac{1}{4}\left( \int_{\T}dk F\star g_{F}\star e^{-i\omega(k)\cdot}(s)+\int_{\T}dk F\star g_{F}\star e^{i\omega(k)\cdot}(s) \right)%\nonumber\\
	%&\qquad\times\left( \int_{\T}dk g_{F}\star e^{-i\omega(k)\cdot}(s)+\int_{\T}dk g_{F}\star e^{i\omega(k)\cdot}(s) \right)\nonumber\\
	= \frac{1}{4}\sum_{\sigma_{1},\sigma_{2}\in\{+,-\}}\int_{\T^{2}} F\star \Theta_{F}\star e^{\sigma_{1}i\omega(k)\cdot}(s)\Theta_{F}\star e^{\sigma_{2}i\omega(h)\cdot}(s)dkdh.
\end{align*}

We compute the Laplace transform for each sign; as before we do the same manipulation
\begin{align}
	&\cL\left( \int_{0}^{t/\varepsilon}\int_{\T^{2}} F\star \Theta_{F}\star e^{\sigma_{1}i\omega(k)\cdot}(s)\Theta_{F}\star e^{\sigma_{2}i\omega(h)\cdot}(s)dkdh \right)(Z)\nonumber\\
	%&= \int_{\T^2}dkd\ell\frac{\varepsilon}{\lambda}\cL\left( F\star g_{F}\star e^{\sigma_{1}i\omega(k)\cdot}(t)g_{F}\star e^{\sigma_{2}i\omega(\ell)\cdot}(t) \right)(\varepsilon\lambda)\nonumber\\
	%&= \int_{\T^2}dkd\ell\frac{1}{\lambda}\cL\left( F\star g_{F}\star e^{\sigma_{1}i\omega(k)\cdot}(t/\varepsilon)g_{F}\star e^{\sigma_{2}i\omega(\ell)\cdot}(t/\varepsilon) \right)(\lambda)\nonumber\\
	&= \int_{\T^2}\frac{1}{Z}\cL\left( \int_{0}^{t/\varepsilon} e^{-\sigma_{1}i\omega(k)s}F\star \Theta_{F}(s)dse^{\sigma_{1}i\omega(k)t/\varepsilon}\int_{0}^{t/\varepsilon}e^{-\sigma_2i\omega(h)s'}\Theta_{F}(ds') e^{\sigma_{2}i\omega(h)t/\varepsilon} \right)(Z)dkdh\nonumber\\
	&= \int_{\T^{2}}\frac{\varepsilon^2}{ 2Z\pi i}\lim_{L\to\infty}\int_{c-iL}^{c+iL}\cL\left( \int_{0}^{t} e^{-\sigma_{1}i\omega(k)s}F\star \Theta_{F}(s)ds \right)(\varepsilon\zeta)\nonumber\\
	&\hspace{0.2\textwidth}\times\cL\left( \int_{0}^{t}e^{-\sigma_2i\omega(h)s'}\Theta_{F}(ds')  \right)(\varepsilon(Z-\zeta)-i(\sigma_{1}\omega(k)+\sigma_{2}\omega(h)))d\zeta dkdh\nonumber\\
	&= \int_{\T^{2}}\frac{1}{2Z\pi i}\lim_{L\to\infty}\int_{c-iL}^{c+iL}\frac{\cL(F)(\varepsilon\zeta+\sigma_{1}i\omega(k))\tilde{\Theta}_{F}(\varepsilon\zeta+\sigma_{1}i\omega(k))\tilde{\Theta}_{F}(\varepsilon\zeta-\sigma_{1}i\omega(k))}{\zeta(Z-\zeta-i\varepsilon^{-1}(\sigma_{1}\omega(k)+\sigma_{2}\omega(h)))}d\zeta dkdh\label{signgeneral}
\end{align}
The case $\sigma_{1}$ has the same sign as $\sigma_{2}$, then we can easily see that the limit is $0$. We estimate the case $\sigma_{1}=+,\sigma_{2}=-$. We change $k,h$ into $k\pm\frac{\varepsilon\xi}{2}$, then \eqref{signgeneral} becomes
\begin{align}
	\int_{\T\times\T_{2/\varepsilon}}\frac{\varepsilon}{2Z\pi i}\lim_{L\to\infty}\int_{c-iL}^{c+iL}\frac{\cL(F)(\varepsilon\zeta+i\omega(k))\tilde{\Theta}_{F}(\varepsilon\zeta+i\omega(k))\tilde{\Theta}_{F}(\varepsilon(Z-\zeta)-i\omega(k))}{\zeta(Z-\zeta-iD_{\varepsilon}(k,\xi))}d\zeta dkdh.
	\label{termJchangevariable}
\end{align}
The Dominated Convergence Theorem does not work here as we don't have a dominating function $\langle\xi\rangle^{-3-\kappa}$, like the one in the first term. Though, the theorem can still be used to prove
\begin{align*}
	\lim_{\varepsilon\to0}\frac{1}{2\pi i}\lim_{L\to\infty}&\int_{c-iL}^{c+iL}\frac{\cL(F)(\varepsilon\zeta+i\omega(k))\tilde{\Theta}_{F}(\varepsilon\zeta+i\omega(k))\tilde{\Theta}_{F}(\varepsilon(Z-\zeta)-i\omega(k))}{\zeta(Z-\zeta-iD_{\varepsilon}(k,\xi))}d\zeta\\
	&= \frac{\hat{F}(\omega(k)/(2\pi))|\theta_{F}(k)|^{2}}{Z-i\omega'(k)\xi}.
\end{align*}

Hence, we have that
\begin{align*}
	\text{\eqref{termJchangevariable}}=\int_{\T\times\R}\frac{1_{[-1/\varepsilon,1/\varepsilon]}(\xi)\varepsilon}{Z(Z-i\omega'(k)\xi)}\hat{F}(\omega(k)/(2\pi))|\theta_{F}(k)|^{2}d\xi dk+\int_{\T\times\T_{2/\varepsilon}}\frac{\varepsilon }{Z}o(\varepsilon)d\xi dk.
\end{align*}
Considering the other case $\sigma_{1}=-,\sigma_{2}=+$, and we taking its inverse Laplace transform, one gets
\begin{align*}
	\int_{\T\times\R}\frac{\sin(2\pi x/\varepsilon)\varepsilon}{2\pi x}1_{[0,v_g(k)t]}(x)\Real(\hat{F}(\omega(k)/(2\pi)))|\theta_{F}(k)|^{2}dx dk+\int_{\T\times\T_{2/\varepsilon}}\varepsilon o(\varepsilon)dx dk.
\end{align*}
The above term converges to $0$. This means that the term \eqref{boundforgrowth} converges to a negative number as $\varepsilon\to0$, which concludes the desired linear estimate.

\bibliographystyle{plain}

\bibliography{references}

\end{document}